\def\published{y}
\documentclass[12pt,reqno]{amsart}
\def\draft{n}
\usepackage{
  fullpage,graphicx,dbnsymb,amsmath,amssymb,floatflt,picins,
  multicol,mathtools,stmaryrd,xcolor
}
\usepackage{hyperref}

\usepackage{hyperref}\hypersetup{colorlinks,
   linkcolor={green!50!black},
   citecolor={green!50!black},
   urlcolor=blue
}

\usepackage[all]{xy}

\if\draft y
  \usepackage[color]{showkeys}
  \usepackage{everypage,draftwatermark}
  \SetWatermarkScale{4}
  \SetWatermarkLightness{0.9}
\fi

\theoremstyle{plain}
\newtheorem{theorem}{Theorem}[section]

\newtheorem{proposition}[theorem]{Proposition}

\newtheorem{lemma}[theorem]{Lemma}

\newtheorem{claim}[theorem]{Claim}

\newtheorem{conjecture}[theorem]{Conjecture}

\theoremstyle{definition}
\newtheorem{definition}[theorem]{Definition}
\newtheorem{defwarn}[theorem]{Definition and Warning}

\newtheorem{problem}[theorem]{Problem}

\theoremstyle{remark}

\newtheorem{comments}[theorem]{Comments}
\newtheorem{example}[theorem]{Example}
\newtheorem{exercise}[theorem]{Exercise}

\newtheorem{remark}[theorem]{Remark}
\newtheorem{warning}[theorem]{Warning}

\newlength{\standardunitlength}
\catcode`\@=11
\long\def\@makecaption#1#2{%
    \vskip 10pt
    \setbox\@tempboxa\hbox{%\ifvoid\tinybox\else\box\tinybox\fi
      \small\sf{\bfcaptionfont #1. }\ignorespaces #2}%
    \ifdim \wd\@tempboxa >\captionwidth {%
        \rightskip=\@captionmargin\leftskip=\@captionmargin
        \unhbox\@tempboxa\par}%
      \else
        \hbox to\hsize{\hfil\box\@tempboxa\hfil}%
    \fi}
\font\bfcaptionfont=cmssbx10 scaled \magstephalf
\newdimen\@captionmargin\@captionmargin=2\parindent
\newdimen\captionwidth\captionwidth=\hsize
\catcode`\@=12

\newcommand{\ad}{\operatorname{ad}}

\newcommand{\Aut}{\operatorname{Aut}}

\newcommand{\gr}{\operatorname{gr}}

\newcommand{\sign}{\operatorname{sign}}

\newcommand{\tr}{\operatorname{tr}}

\def\qed{{\hfill\text{$\Box$}}}

\newlength{\globalparindent}
\newenvironment{myitemize}{
        \begin{list}{$\bullet$}{\setlength{\leftmargin}{16pt}
        \setlength{\labelwidth}{12pt}
        \setlength{\labelsep}{4pt}}
}{
        \end{list}
}

\def\arXiv#1{{\href{http://front.math.ucdavis.edu/#1}{arXiv:#1}}}

\def\bbQ{{\mathbb Q}}
\def\bbR{{\mathbb R}}
\def\bbZ{{\mathbb Z}}
\def\calA{{\mathcal A}}
\def\calB{{\mathcal B}}
\def\calAcf{{\mathcal A}^\text{\it cf}}

\def\calD{{\mathcal D}}
\def\calF{{\mathcal F}}
\def\calG{{\mathcal G}}
\def\calI{{\mathcal I}}
\def\calK{{\mathcal K}}
\def\calL{{\mathcal L}}

\def\calP{{\mathcal P}}
\def\calR{{\mathcal R}}

\def\calT{{\mathcal T}}
\def\calU{{\mathcal U}}
\def\fraka{{\mathfrak a}}
\def\frakg{{\mathfrak g}}
\def\tilE{{\tilde{E}}}

\def\ad{\operatorname{ad}}
\def\Aut{\operatorname{Aut}}
\def\Autop{\operatorname{Aut}^\text{op}}

\def\gr{\operatorname{gr}\,}
\def\IAM{\operatorname{\it IAM}}
\def\proj{{\operatorname{proj}\,}}

\def\attr{\operatorname{\mathfrak{tr}}}

\def\sder{\operatorname{\mathfrak{sder}}}
\def\tder{\operatorname{\mathfrak{tder}}}

\def\FA{{\mathit F\!A}}
\def\uB{{\mathit u\!B}}
\def\PuB{{\mathit P\!u\!B}}
\def\vB{{\mathit v\!B}}
\def\vT{{\mathit v\!T}}
\def\PvB{{\mathit P\!v\!B}}
\def\wB{{\mathit w\!B}}
\def\swB{{\mathit s\!w\!B}}
\def\uT{{\mathit u\!T}}
\def\wT{{\mathit w\!T}}

\def\PwB{{\mathit P\!w\!B}}
\def\sl{{\mathit sl}}

\def\aft{{\overrightarrow{4T}}}
\def\aAS{{\overrightarrow{AS}}}
\def\aSTU{{\overrightarrow{STU}}}
\def\aIHX{{\overrightarrow{IHX}}}
\def\Rs{{R$1^{\!s}$}}

\def\mathsize#1#2{{\begin{array}{c}\text{#1$\!#2\!$}\end{array}}}
\def\pstex#1{\begin{array}{c}
  \input figs/#1.pstex_t
  \if\draft y
    \smash{\makebox[0in]{\color{labelkey}{\tt figs/#1}}}
  \fi
\end{array}}

\def\draftcut{\if\draft y \newpage \fi}

\def\glos#1{\setlength{\fboxsep}{0pt}\colorbox{yellow}{$#1$}}
\def\glost#1{\setlength{\fboxsep}{0pt}\colorbox{yellow}{#1}}
\def\red{\color{red}}

\newcommand\ifpub[2]{\if\published y #1 \else #2 \fi}

\begin{document}
\newdimen\captionwidth\captionwidth=\hsize
\setcounter{secnumdepth}{4}

\title[Finite Type Invariants of w-Knotted Objects]{Finite Type Invariants
  of w-Knotted Objects I: w-Knots and the Alexander Polynomial}

\author{Dror~Bar-Natan and Zsuzsanna Dancso}
\address{
  Department of Mathematics\\
  University of Toronto\\
  Toronto Ontario M5S 2E4\\
  Canada
}
\email{drorbn@math.toronto.edu, zsuzsi@math.toronto.edu}
\urladdr{http://www.math.toronto.edu/~drorbn, http://www.math.toronto.edu/zsuzsi}

\date{first edition Sep.\ 27, 2013, this edition Jun.~30,~2015. The
\arXiv{1309.7155} edition may be older}


\thanks{This work was partially supported by NSERC grant RGPIN 262178.
  Electronic version, videos (wClips) and related files at~\cite{WKO},
  \url{http://www.math.toronto.edu/~drorbn/papers/WKO/}.
}

\begin{abstract}
  \ifpub{This is the first in a series of papers studying w-knots, and more generally,
w-knotted objects (w-braids, w-tangles,
etc.). These  are classes of knotted objects which are \underline{w}ider, but
\underline{w}eaker than their ``\underline{u}sual'' counterparts.}
{This is the first in a series of papers studying w-knots, and more generally,
w-knotted objects (w-braids, w-tangles,
etc.). These  are classes of knotted objects which are \underline{w}ider, but
\underline{w}eaker than their ``\underline{u}sual'' counterparts. To
get (say) w-knots from usual knots (or u-knots), one has to allow non-planar ``virtual''
knot diagrams, hence enlarging the base set of knots. But then one
imposes a new relation beyond the ordinary collection of Reidemeister moves,
called the ``overcrossings commute'' relation, making w-knotted
objects a bit weaker once again.}

The group of w-braids was studied (under the name
``\underline{w}elded braids'') by Fenn, Rimanyi and
Rourke~\cite{FennRimanyiRourke:BraidPermutation} and was shown to
be isomorphic to the McCool group~\cite{McCool:BasisConjugating}
of ``basis-conjugating'' automorphisms of a free group $F_n$ ---
the smallest subgroup of $\Aut(F_n)$ that contains both braids and
permutations. Brendle and Hatcher~\cite{BrendleHatcher:RingsAndWickets},
in work that traces back to Goldsmith~\cite{Goldsmith:MotionGroups},
have shown this group to be a group of movies of flying rings in
$\bbR^3$. Satoh~\cite{Satoh:RibbonTorusKnots} studied several classes of
w-knotted objects (under the name ``\underline{w}eakly-virtual'') and has
shown them to be closely related to certain classes of knotted surfaces
in $\bbR^4$. So w-knotted objects are algebraically and topologically 
interesting.

In this article we study finite type invariants of w-brainds and w-knots.  
Following Berceanu and
Papadima~\cite{BerceanuPapadima:BraidPermutation}, we construct
homomorphic universal finite type invariants of w-braids. 
We find that the universal finite type invariant of w-knots
is essentially the Alexander polynomial.

Much as the spaces $\calA$ of chord diagrams for ordinary knotted
objects are related to metrized Lie algebras, we find that the spaces
$\calA^w$ of ``arrow diagrams'' for w-knotted objects are related to
not-necessarily-metrized Lie algebras. Many questions concerning w-knotted
objects turn out to be equivalent to questions about Lie algebras, and
in later papers of this series we re-interpret Alekseev-Torossian's~\cite{AlekseevTorossian:KashiwaraVergne}
work on Drinfel'd associators and the Kashiwara-Vergne problem 
as a study of w-knotted trivalent graphs.

\ifpub{}{The true value of w-knots, though, is likely to emerge later, for we
expect them to serve as a warmup example for what we expect
will be even more interesting --- the study of virtual
knots, or v-knots. We expect v-knotted objects to provide the global
context whose projectivization (or ``associated graded structure'')
will be the Etingof-Kazhdan theory of deformation quantization of Lie
bialgebras~\cite{EtingofKazhdan:BialgebrasI}.}

\end{abstract}

\maketitle

\clearpage
{\small \tableofcontents}
%{\small \begin{multicols}{2} \tableofcontents \end{multicols}}
%{\small \twocolumn \tableofcontents \onecolumn}

\draftcut
\section{Introduction} \label{sec:intro}

\subsection{Dreams} \label{subsec:dreams} We
have a dream\footnote{Understanding
the authors' history and psychology ought never be necessary to
understand their papers, yet it may be helpful. Nothing material
in the rest of this paper relies on Section~\ref{subsec:dreams}.},
at least partially founded on reality, that many of the difficult
algebraic equations in mathematics, especially those that are
written in graded spaces, more especially those that are related in
one way or another to quantum groups,
and even more especially those related to the work of Etingof and
Kazhdan~\cite{EtingofKazhdan:BialgebrasI}, can be understood, and indeed,
would appear more natural, in terms of finite type invariants of various
topological objects.

We believe this is the case for Drinfel'd's theory
of associators~\cite{Drinfeld:QuasiHopf}, which can be
interpreted as a theory of well-behaved universal finite type
invariants of parenthesized tangles\footnote{``$q$-tangles''
in~\cite{LeMurakami:Universal}, ``non-associative tangles''
in~\cite{Bar-Natan:NAT}.}~\cite{LeMurakami:Universal, Bar-Natan:NAT},
and as a theory of universal finite type invariants
of knotted trivalent graphs~\cite{Dancso:KIforKTG}.

We believe this is the case for Drinfel'd's ``Grothendieck-Teichm\"uller
group''~\cite{Drinfeld:GalQQ}, which is better understood as a
group of automorphisms of a certain algebraic structure, also
related to universal finite type invariants of parenthesized
tangles~\cite{Bar-Natan:Associators}.

And we're optimistic, indeed we believe, that sooner or later the
work of Etingof and Kazhdan~\cite{EtingofKazhdan:BialgebrasI}
on quantization of Lie bialgebras will be re-interpreted as a
construction of a well-behaved universal finite type invariant of
virtual knots~\cite{Kauffman:VirtualKnotTheory} or of some other class
of virtually knotted objects. Some steps in that direction were taken
by Haviv~\cite{Haviv:DiagrammaticAnalogue}.

We have another dream, to construct a useful ``Algebraic Knot Theory''. As
at least a partial writeup exists~\cite{Bar-Natan:AKT-CFA},
we'll only state that an important ingredient necessary to
fulfil that dream would be a ``closed form''\footnote{The
phrase ``closed form'' in itself requires an explanation. See
Section~\ref{subsec:ClosedForm}. \label{foot:ClosedForm}} formula for an
associator, at least in some reduced sense. Formulae for associators or
reduced associators were in themselves the goal of several studies
undertaken for various other reasons~\cite{LeMurakami:HOMFLY,
Lieberum:gl11, Kurlin:CompressedAssociators, LeeP:ClosedForm}.

\draftcut \subsection{Stories}

Thus, the first named author, DBN, was absolutely delighted
when in January 2008 Anton Alekseev described to him his joint
work~\cite{AlekseevTorossian:KashiwaraVergne} with Charles
Torossian --- Anton told DBN that they found a relationship between the
Kashiwara-Vergne conjecture~\cite{KashiwaraVergne:Conjecture},
a cousin of the Duflo isomorphism (which DBN already knew to be
knot-theoretic~\cite{Bar-NatanLeThurston:TwoApplications}), and
associators taking values in a space called $\sder$, which he could
identify as ``tree-level Jacobi diagrams'', also a knot-theoretic
space related to the Milnor invariants~\cite{Bar-Natan:Homotopy,
HabeggerMasbaum:Milnor}. What's more, Anton told DBN that in certain
quotient spaces the Kashiwara-Vergne conjecture can be solved explicitly;
this should lead to some explicit associators!

So DBN spent the following several months trying to
understand~\cite{AlekseevTorossian:KashiwaraVergne} which eventually led to this sequence of papers. 
One main thing we learned is that
the Alekseev-Torossian paper, and with it the Kashiwara-Vergne (KV)
conjecture, fit very nicely with our first dream recalled above,
about interpreting algebra in terms of knot theory. Indeed much
of~\cite{AlekseevTorossian:KashiwaraVergne} can be reformulated as a
construction and a discussion of a well-behaved universal finite type
invariant\footnote{The notation $Z$ for universal finite type invarants
comes from the famous universal finite type invariant of classical links, the Kontsevich itegral.} $Z$
of a certain class of knotted objects (which we will call 
w-knotted), a certain natural quotient of the space of virtual knots
(more precisely, virtual trivalent tangles): this will be the subject of the second
paper in the series. It is also possible to provide a topological interpretation 
(and independent topological proof) of the \cite{AlekseevEnriquezTorossian:ExplicitSolutions}
formula for explicit solutions to the KV problem in terms of associators. This will be done in 
the third paper.
And our hopes remain high
that later we (or somebody else) will be able to exploit this relationship
in directions compatible with our second dream recalled above, on the
construction of an ``algebraic knot theory''.

The story, in fact, is prettier than we were hoping for, for it has the
following additional qualities:

\begin{myitemize}

\item w-Knotted objects are quite interesting in themselves: as
stated in the abstract, they are related to combinatorial group
theory via ``basis-conjugating'' automorphisms of a free group $F_n$,
to groups of movies of flying rings in $\bbR^3$, and more generally, to
certain classes of knotted surfaces in $\bbR^4$. The references include
\cite{Goldsmith:MotionGroups, McCool:BasisConjugating, FennRimanyiRourke:BraidPermutation, 
Satoh:RibbonTorusKnots, BrendleHatcher:RingsAndWickets}.

\item The ``chord diagrams'' for w-knotted objects (really, these are ``arrow
diagrams'') describe formulae for invariant tensors in spaces pertaining to
not-necessarily-metrized Lie algebras in much of the same way as ordinary
chord diagrams for ordinary knotted objects describe formulae for invariant
tensors in spaces pertaining to metrized Lie algebras. This observation is
bound to have further implications.

\item Arrow diagrams also describe the Feynman diagrams of topological BF
theory \cite{CattaneoCotta-RamusinoMartellini:Alexander, CCFM:BF34} and of a
certain class of Chern-Simons theories~\cite{Naot:BF}. Thus, it is likely that
our story is directly related to quantum field theory\footnote{Some
non-perturbative relations between BF theory and w-knots was discussed by
Baez, Wise and Crans~\cite{BaezWiseCrans:ExoticStatistics}.}.

\item The main objective of this paper is to prove that, when composed 
with the map from knots to w-knots, $Z$ becomes the
Alexander polynomial. For links, it becomes an invariant stronger than the
multi-variable Alexander polynomial which contains the multi-variable
Alexander polynomial as an easily identifiable reduction. 

\item On other
w-knotted objects $Z$ has easily identifiable reductions that can be
considered as ``Alexander polynomials'' with good behaviour relative
to various knot-theoretic operations --- cablings, compositions
of tangles, etc. There is also a certain specific reduction of $Z$
that can be considered as an ``ultimate Alexander polynomial'' ---
in the appropriate sense, it is the minimal extension of the Alexander
polynomial to other knotted objects which is well behaved under a whole
slew of knot theoretic operations, including the ones named above. See
\cite{Bar-NatanSelmani:MetaMonoids, Bar-Natan:KBH}.

\ifpub{\item The true value of w-knots, though, is likely to emerge later, for we
expect them to serve as a \underline{w}armup example for what we expect
will be even more interesting --- the study of \underline{v}irtual
knots, or v-knots. We expect v-knotted objects to provide the global
context whose projectivization (or ``associated graded structure'')
will be the Etingof-Kazhdan theory of deformation quantization of Lie
bialgebras~\cite{EtingofKazhdan:BialgebrasI}.}{}

\end{myitemize}

\begin{figure}
\[
  \def\uT{\parbox[t]{1.875in}{\small
    Ordinary (\underline{u}sual) knotted objects in 3D --- braids,
    knots, links, tangles, knotted graphs, etc.
  }}
  \def\vT{\parbox[t]{1.875in}{\small
    \underline{V}irtual knotted objects --- ``algebraic'' knotted objects,
    or ``not specifically embedded'' knotted objects; knots drawn on a
    surface, modulo stabilization.
  }}
  \def\wT{\parbox[t]{1.875in}{\small
    Ribbon knotted objects in 4D; ``flying rings''. Like v, but also with
    ``overcrossings commute''.
  }}
  \def\uC{\parbox[t]{1.875in}{\small
    Chord diagrams and Jacobi diagrams, modulo $4T$, $STU$, $IHX$, etc.
  }}
  \def\vC{\parbox[t]{1.875in}{\small
   Arrow diagrams and v-Jacobi diagrams, modulo $6T$ and various
   ``directed'' $STU$s and $IHX$s, etc.
  }}
  \def\wC{\parbox[t]{1.875in}{\small
    Like v, but also with ``tails commute''. Only ``two in one out''
    internal vertices.
  }}
  \def\uL{\parbox[t]{1.875in}{\small
    Finite dimensional metrized Lie algebras, representations,  and
    associated spaces.
  }}
  \def\vL{\parbox[t]{1.875in}{\small
    Finite dimensional Lie bi-algebras, representations,  and associated
    spaces.
  }}
  \def\wL{\parbox[t]{1.875in}{\small
    Finite dimensional co-commutative Lie bi-algebras (i.e.,
    $\frakg\ltimes\frakg^\ast$), representations,  and associated
    spaces.
  }}
  \def\uH{\parbox[t]{1.875in}{\small
    The Drinfel'd theory of associators.
  }}
  \def\vH{\parbox[t]{1.875in}{\small
    Likely, quantum groups and the Etingof-Kazhdan theory of quantization
    of Lie bi-algebras.
  }}
  \def\wH{\parbox[t]{1.875in}{\small
    The Kashiwara-Vergne-Alekseev-Torossian theory of convolutions on Lie
    groups and Lie algebras.
  }}
  \input figs/uvw.pstex_t
\]
\caption{The u-v-w Stories} \label{fig:uvw}
\end{figure}

\draftcut \subsection{The Bigger Picture} 
Parallel to the w-story run the possibly more significant u-story
and v-story. The u-story is about u-knots, or more generally,
u-knotted objects (braids, links, tangles, etc.), where ``u'' stands
for \underline{u}sual; hence the u-story is about classical knot
theory. The v-story is about v-knots, or more generally, v-knotted
objects, where ``v'' stands for \underline{v}irtual, in the sense of
Kauffman~\cite{Kauffman:VirtualKnotTheory}.

The u, v, and w-knotted objects, are quite different from each other. Yet they
can be told along similar lines --- first the knots (topology), then their
finite type invariants and their ``chord diagrams'' (combinatorics),
then those map into certain universal enveloping algebras and similar
spaces associated with various classes of Lie algebras (low algebra),
and finally, in order to construct a ``good'' universal finite type
invariant, in each case one has to confront a certain deeper algebraic
subject (high algebra). These stories are summarized in a table form
in Figure~\ref{fig:uvw}.

u-Knots map into v-knots, and v-knots map into w-knots\footnote{Though
the composition ``$u\to v\to w$'' is not $0$. In fact, the composed
map $u\to w$ is injective. u-Knots, for example, are determined by the
fundamental groups of their complements plus ``peripheral systems''
(or alternatively, by their ``quandles''~\cite{Joyce:TheKnotQuandle}),
and this information is easily recovered from the w-knot images of
u-knots. Similar considerations apply to other classes of u-knotted
objects.}. The other parts of our stories, the ``combinatorics'' and
``low algebra'' and ``high algebra'' rows of Figure~\ref{fig:uvw}, are
likewise related, and this relationship is a crucial part of our overall
theme. Thus, we cannot and will not tell the w-story in isolation, and
while it is central to this article, we will necessarily also include
some episodes from the u and v series.

\subsection{Plans} In this paper we study w-braids and w-knots; the main
result is Theorem \ref{thm:Alexander}, which states that the universal
finite type invariant of w-knots is essenially the Alexander polynomial.
However, starting with braids and taking a classical approach to 
finite type invariants, this paper also serves as a gentle introduction 
to the subsequent papers and in particular to \cite{Bar-NatanDancso:WKO2} 
where we will present a more algebraic point of view.
For more detailed information on the content 
consult the ``Section Summary'' paragraphs below and at the beginning
of each section. An ``odds and ends'' section and a glossary of notation 
follows the main sections.

\def\summarybraids{This section is largely a compilation of existing
literature, though we also introduce the language of arrow diagrams that
we use throughout the rest of the paper. In Sections~\ref{subsec:VirtualBraids}
and~\ref{subsec:wBraids} we define v-braids and then w-braids and
survey their relationship with basis-conjugating automorphisms of free
groups and with ``the group of (horizontal) flying rings in $\bbR^3$''
(really, a group of knotted tubes in $\bbR^4$). In Section ~\ref{subsec:FT4Braids}
we play the usual game of introducing finite type invariants, weight
systems, chord diagrams (arrow diagrams, for this case), and 4T-like
relations. In Section~\ref{subsec:wBraidExpansion} we define and construct
a universal finite type invariant $Z$ for w-braids --- it turns out
that the only algebraic tool we need to use is the formal exponential
function $\exp(a):=\sum a^n/n!$. In Section~\ref{subsec:bcomments} we study
some good algebraic properties of $Z$, its injectivity, and its
uniqueness, and we conclude with the slight modifications needed for the
study of non-horizontal flying rings.}

\def\summaryknots{In Section~\ref{subsec:VirtualKnots} we define v-knots and
w-knots (long v-knots and long w-knots, to be precise) and discuss a map
$v\to w$. In Section~\ref{subsec:FTforvwKnots} we determine the space of ``chord
diagrams'' for w-knots to be the space $\calA^w(\uparrow)$ of arrow
diagrams modulo $\aft$ and TC relations and in Section~\ref{subsec:SomeDimensions}
we compute some relevant dimensions. In Section~\ref{subsec:Jacobi} we show
that $\calA^w(\uparrow)$ can be re-interpreted as a space of trivalent
graphs modulo STU- and IHX-like relations, and is therefore related
to Lie algebras (Section~\ref{subsec:LieAlgebras}). This allows us to
completely determine $\calA^w(\uparrow)$.  With no difficulty 
in Section~\ref{subsec:Z4Knots} we construct a universal finite type
invariant for w-knots. With a bit of further difficulty we show in
Section~\ref{subsec:Alexander} that it is essentially equal to the Alexander
polynomial.}

%\noindent
%{\small \begin{multicols}{2}

{\bf Section~\ref{sec:w-braids}, w-Braids.} %(page~\pageref{sec:w-braids})
\summarybraids

{\bf Section~\ref{sec:w-knots}, w-Knots.} %(page~\pageref{sec:w-knots})
\summaryknots

%$\end{multicols}}

\subsection{Acknowledgements} We wish to thank the Anonymous Referee, Anton Alekseev, Jana
Archibald, Scott Carter, Karene Chu, Iva Halacheva, Joel Kamnitzer,
Lou Kauffman, Peter Lee, Louis Leung, Jean-Baptiste Meilhan, Dylan Thurston, Daniel Tubbenhauer and Lucy Zhang
for comments and suggestions.

\draftcut
\section{w-Braids} \label{sec:w-braids}

\begin{quote} \small {\bf Section Summary. }
  \summarybraids
\end{quote}

\subsection{Preliminary: Virtual Braids, or v-Braids.}
\label{subsec:VirtualBraids}
Our main object of study for this section, w-braids, are best
viewed as ``virtual braids''~\cite{Bardakov:VirtualAndUniversal,
KauffmanLambropoulou:VirtualBraids, BardakovBellingeri:VirtualBraids},
or v-braids, modulo one additional relation; hence, we start with v-braids.

It is simplest to define v-braids in
terms of generators and relations, either algebraically or pictorially.
This can be done in at least two ways -- the easier-at-first but
philosophically less satisfying ``planar'' way, and the harder-to-digest 
but morally more correct ``abstract'' way.\footnote{Compare with a
similar choice that exists in the definition of manifolds, as either
appropriate subsets of some ambient Euclidean spaces (module some
equivalences) or as abstract gluings of coordinate patches (modulo some
other equivalences). Here in the ``planar'' approach of
Section~\ref{subsubsec:Planar} we consider v-braids
as ``planar'' objects, and in the ``abstract approach'' of
Section~\ref{subsubsec:Abstract} they are just ``gluings'' of abstract
``crossings'', not drawn anywhere in particular.}

\subsubsection{The ``Planar'' Way} \label{subsubsec:Planar} For a
natural number $n$ set $\glos{\vB_n}$ to be the group generated by
symbols $\glos{\sigma_i}$ ($1\leq i\leq n-1$), called ``crossings''
and graphically represented by an overcrossing $\overcrossing$ ``between
strand $i$ and strand $i+1$'' (with inverse $\undercrossing$)\footnote{We
sometimes refer to $\overcrossing$ as a ``positive crossing'' and to
$\undercrossing$ as a ``negative crossing''.}, and $\glos{s_i}$, called
``virtual crossings'' and graphically represented by a non-crossing,
$\virtualcrossing$, also ``between strand $i$ and strand $i+1$'',
subject to the following relations:

\begin{myitemize}

\item The subgroup of $\vB_n$ generated by the virtual crossings $s_i$
is the symmetric group $\glos{S_n}$, and the $s_i$'s correspond to the
transpositions $(i,i+1)$. That is, we have
\begin{equation} \label{eq:sRelations}
  s_i^2=1,
  \qquad s_is_{i+1}s_i = s_{i+1}s_is_{i+1},
  \qquad\text{and if $|i-j|>1$, then}
  \qquad s_is_j=s_js_i.
\end{equation}
In pictures, this is
\begin{equation} \label{eq:sRels}
  \def\i{{$i$}}
  \def\ip{{$i\!+\!1$}}
  \def\ipp{{$i\!+\!2$}}
  \def\j{{$j$}}
  \def\jp{{$j\!+\!1$}}
  \input figs/sRels.pstex_t
\end{equation}
Note that we read our braids from bottom to top, and that all relations (and most pitcures in this paper) are local:
the braids may be bigger than shown but the parts not shown remain the same throughout a relation.

\item The subgroup of $\vB_n$ generated by the crossings $\sigma_i$'s is
the usual braid group $\glos{\uB_n}$, and $\sigma_i$ corresponds to the ``braiding
of strand $i$ over strand $i+1$''. That is, we have
\begin{equation} \label{eq:R3}
  \sigma_i\sigma_{i+1}\sigma_i
    = \sigma_{i+1}\sigma_i\sigma_{i+1},
  \qquad\text{and if $|i-j|>1$ then}
  \qquad \sigma_i\sigma_j=\sigma_j\sigma_i.
\end{equation}
In pictures, dropping the indices, this is
\begin{equation} \label{eq:sigmaRels} \begin{picture}(0,0)%
\includegraphics{figs/sigmaRels.pstex}%
\end{picture}%
%
%  pstex_opts: -m 0.8 
%
\setlength{\unitlength}{3158sp}%
\begingroup\makeatletter\ifx\SetFigFont\undefined%
\gdef\SetFigFont#1#2#3#4#5{%
  \reset@font\fontsize{#1}{#2pt}%
  \fontfamily{#3}\fontseries{#4}\fontshape{#5}%
  \selectfont}%
\fi\endgroup%
\begin{picture}(5154,924)(1639,-673)
\put(5551,-286){\makebox(0,0)[b]{\smash{{\SetFigFont{10}{12.0}{\rmdefault}{\mddefault}{\updefault}{\color[rgb]{0,0,0}$=$}%
}}}}
\put(2401,-286){\makebox(0,0)[b]{\smash{{\SetFigFont{10}{12.0}{\rmdefault}{\mddefault}{\updefault}{\color[rgb]{0,0,0}$=$}%
}}}}
\end{picture}%
 \end{equation}
The first of these relations is the ``Reidemeister 3 move''\footnote{The
Reidemeister 2 move is the relations $\sigma_i\sigma_i^{-1}=1$ which is
part of the definition of a group. There is no Reidemeister 1 move in
the theory of braids.} of knot theory. The second is sometimes called
``locality in space''~\cite{Bar-Natan:NAT}.

\item Some ``mixed relations'', that is,
\begin{equation} \label{eq:MixedRelations}
  s_i\sigma^{\pm 1}_{i+1}s_i = s_{i+1}\sigma^{\pm 1}_is_{i+1},
  \qquad\text{and if $|i-j|>1$, then}
  \qquad s_i\sigma_j=\sigma_js_i.
\end{equation}
In pictures, this is
\begin{equation} \label{eq:MixedRels}
  \input figs/MixedRels.pstex_t
\end{equation}

\end{myitemize}

\begin{remark} \label{rem:Skeleton} The ``skeleton'' of a v-braid $B$
is the set of strands appearing in it, retaining the association
between their beginning and ends but ignoring all the crossing
information. More precisely, it is the permutation induced by
tracing along $B$, and even more precisely it is the image of $B$
via the ``skeleton morphism'' $\glos{\varsigma}\colon\vB_n\to S_n$ defined
by $\varsigma(\sigma_i)=\varsigma(s_i)=s_i$ (or pictorially, by
$\varsigma(\overcrossing)=\varsigma(\virtualcrossing)=\virtualcrossing$).
Thus, the symmetric group $S_n$ is both a subgroup and a quotient group
of $\vB_n$.
\end{remark}

Like there are pure braids to accompany braids, there are pure virtual
braids as well:

\begin{definition} A pure v-braid is a v-braid whose skeleton is the
identity permutation; the group $\glos{\PvB_n}$ of all pure v-braids is
simply the kernel of the skeleton morphism $\varsigma\colon\vB_n\to S_n$.
\end{definition}

We note the short exact sequence of group homomorphisms
\begin{equation} \label{eq:ExcatSeqForPvB}
  1\longrightarrow\PvB_n\xhookrightarrow{\quad}\vB_n
  \overset{\varsigma}{\longrightarrow}S_n
  \longrightarrow 1.
\end{equation}
This short exact sequence splits, with the splitting given by the inclusion
$S_n\hookrightarrow\vB_n$ mentioned above~\eqref{eq:sRelations}.
Therefore, we have that
\begin{equation} \label{eq:vBSemiDirect}
  \vB_n=\PvB_n\rtimes S_n.
\end{equation}

\subsubsection{The ``Abstract'' Way} \label{subsubsec:Abstract}
The relations~\eqref{eq:sRels} and~\eqref{eq:MixedRels} that govern
the behaviour of virtual crossings precisely say that virtual crossings
really are ``virtual'' --- if a piece of strand is routed within a braid
so that there are only virtual crossings around it, it can be rerouted
in any other ``virtual only'' way, provided the ends remain fixed
(this is Kauffman's ``detour move''~\cite{Kauffman:VirtualKnotTheory,
KauffmanLambropoulou:VirtualBraids}). Since a v-braid $B$ is independent
of the routing of virtual pieces of strand, we may as well never supply
this routing information.

\parpic[r]{$\input figs/PvBExample.pstex_t$}
Thus, for example, a perfectly fair verbal description of the
(pure!) v-braid on the right is ``strand 1 goes over strand 3 by a
positive crossing then likewise positively over strand 2 then negatively
over 3 then 2 goes positively over 1''. We don't need to specify how
strand 1 got to be near strand 3 so it can go over it --- it got there
by means of virtual crossings, and it doesn't matter how. Hence we arrive
at the following ``abstract'' presentation of $\PvB_n$ and $\vB_n$:

\begin{proposition} (E.g.~\cite[Theorems 1 and 2]{Bardakov:VirtualAndUniversal})
\begin{enumerate}
\item The group $\PvB_n$ of pure v-braids
is isomorphic to the group generated by symbols $\glos{\sigma_{ij}}$ for
$1\leq i\neq j\leq n$ (meaning ``strand $i$ crosses over strand $j$
at a positive crossing''\footnote{The inverse, $\sigma_{ij}^{-1}$,
is ``strand $i$ crosses over strand $j$ at a negative crossing''}),
subject to the third Reidemeister move and to locality in space (compare
with~\eqref{eq:R3} and~\eqref{eq:sigmaRels}):

\begin{align*}
  \sigma_{ij}\sigma_{ik}\sigma_{jk} &= \sigma_{jk}\sigma_{ik}\sigma_{ij}
    & \text{whenever}\qquad & |\{i,j,k\}|=3, \\
  \sigma_{ij}\sigma_{kl} &= \sigma_{kl}\sigma_{ij}
     & \text{whenever}\qquad & |\{i,j,k,l\}|=4.
\end{align*}
\item If $\tau\in S_n$, then with the action
$\sigma_{ij}^\tau:=\sigma_{\tau i,\tau j}$ we recover the semi-direct
product decomposition $\vB_n=\PvB_n\rtimes S_n$. \qed
\end{enumerate}
\end{proposition}

\draftcut \subsection{On to w-Braids} \label{subsec:wBraids}
To define w-braids, we break the symmetry between overcrossings and
undercrossings by imposing one of the ``forbidden moves'' in virtual knot
theory, but not the other:
\begin{equation} \label{eq:OvercrossingsCommute}
  \sigma_i\sigma_{i+1}s_i = s_{i+1}\sigma_i\sigma_{i+1},
  \qquad\text{yet}\qquad
  s_i\sigma_{i+1}\sigma_i \neq \sigma_{i+1}\sigma_is_{i+1}.
\end{equation}
Alternatively,
\[ \sigma_{ij}\sigma_{ik} = \sigma_{ik}\sigma_{ij},
  \qquad\text{yet}\qquad
  \sigma_{ik}\sigma_{jk} \neq \sigma_{jk}\sigma_{ik}.
\]
In pictures, this is
\begin{equation} \label{eq:OC} \input figs/OCUC.pstex_t \end{equation}

The relation we have just imposed may be called the ``unforbidden
relation'', or, perhaps more appropriately, the ``overcrossings commute''
relation, abbreviated \glost{OC}. Ignoring the non-crossings\footnote{Why this is
appropriate was explained in the previous section.} $\virtualcrossing$,
the OC relation says that it is the same if strand $i$ first crosses
over strand $j$ and then over strand $k$, or if it first crosses over
strand $k$ and then over strand $j$. The ``undercrossings commute''
relation \glost{UC}, the one we do not impose
in~\eqref{eq:OvercrossingsCommute}, would say the same except with
``under'' replacing ``over''.

\begin{definition} The group of w-braids is $\glos{\wB_n}:=\vB_n/OC$. Note
that $\varsigma$ descends to $\wB_n$, and hence we can define the group $\glos{\PwB_n}$ of
pure w-braids to be the kernel of the map $\varsigma\colon\wB_n\to S_n$. We
still have a split exact sequence as at~\eqref{eq:ExcatSeqForPvB} and
a thus, a semi-direct product decomposition $\wB_n=\PwB_n\rtimes S_n$.
\end{definition}

\begin{exercise} Show that the OC relation is equivalent to the relation
\[
  \sigma_i^{-1}s_{i+1}\sigma_i = \sigma_{i+1}s_i\sigma_{i+1}^{-1}
  \qquad\text{or}\qquad
  \parbox[m]{1.5in}{\begin{picture}(0,0)%
\includegraphics{figs/AltOC.pstex}%
\end{picture}%
%
%  pstex_opts: -m 0.8 
%
\setlength{\unitlength}{3158sp}%
\begingroup\makeatletter\ifx\SetFigFont\undefined%
\gdef\SetFigFont#1#2#3#4#5{%
  \reset@font\fontsize{#1}{#2pt}%
  \fontfamily{#3}\fontseries{#4}\fontshape{#5}%
  \selectfont}%
\fi\endgroup%
\begin{picture}(1524,924)(1639,-673)
\put(2401,-286){\makebox(0,0)[b]{\smash{{\SetFigFont{10}{12.0}{\rmdefault}{\mddefault}{\updefault}{\color[rgb]{0,0,0}$=$}%
}}}}
\end{picture}%
 }
\]
\end{exercise}

While for most of this paper the pictorial / algebraic definition of w-braids
(and other w-knotted objects) will suffice, we ought describe at least
briefly a few further interpretations of $\wB_n$:

\subsubsection{The Group of Flying Rings} \label{subsubsec:FlyingRings}
Let \glos{X_n} be the space of all placements of $n$ numbered disjoint
geometric circles in $\bbR^3$, such that all circles are parallel to
the $xy$ plane. Such placements will be called horizontal\footnote{
For the group of non-horizontal flying rings see
Section \ref{subsubsec:NonHorRings}.}. A horizontal
placement is determined by the centres in $\bbR^3$ of the $n$ circles
and by $n$ radii, so $\dim X_n=3n+n=4n$. The permutation group $S_n$
acts on $X_n$ by permuting the circles, and one may think of the quotient
$\glos{\tilde{X}_n}:=X_n/S_n$ as the space of all horizontal
placements of $n$ unmarked circles in $\bbR^3$. The fundamental group
$\pi_1(\tilde{X}_n)$ is a group of paths traced by $n$ disjoint horizontal
circles (modulo homotopy), so it is fair to think of it as ``the group
of flying rings''.

\begin{theorem} The group of pure w-braids $\PwB_n$ is isomorphic to the group
of flying rings $\pi_1(X_n)$. The group $\wB_n$ is isomorphic to the group
of unmarked flying rings $\pi_1(\tilde{X}_n)$.
\end{theorem}

For the proof of this theorem,
see~\cite{Goldsmith:MotionGroups, Satoh:RibbonTorusKnots} and
especially~\cite[Proposition~3.3]{BrendleHatcher:RingsAndWickets}. Here
we will contend ourselves with pictures describing the images of the
generators of $\wB_n$ in $\pi_1(\tilde{X}_n)$ and a few comments:

\[ \input figs/FlyingRings.pstex_t \]

Thus, we map the permutation $s_i$ to the movie clip in which ring
number $i$ trades its place with ring number $i+1$ by having the two
flying around each other. This acrobatic feat is performed in $\bbR^3$
and it does not matter if ring number $i$ goes ``above'' or ``below''
or ``left'' or ``right'' of ring number $i+1$ when they trade places,
as all of these possibilities are homotopic. More interestingly, we
map the braiding $\sigma_i$ to the movie clip in which ring $i+1$
shrinks a bit and flies through ring $i$. It is a worthwhile
exercise for the reader to verify that the relations in the definition
of $\wB_n$ become homotopies of movie clips. Of these relations it
is most interesting to see why the ``overcrossings commute'' relation
$\sigma_i\sigma_{i+1}s_i = s_{i+1}\sigma_i\sigma_{i+1}$ holds, yet the
``undercrossings commute'' relation $\sigma^{-1}_i\sigma^{-1}_{i+1}s_i =
s_{i+1}\sigma^{-1}_i\sigma^{-1}_{i+1}$ doesn't.

\begin{exercise}\label{ex:swBn}
To be perfectly precise, we have to
specify the fly-through direction. In our notation, $\sigma_i$ means
that the ring corresponding to the strand going under (in the local picture 
for $\sigma_i$) approaches from below the bigger
ring representing the strand going over, then flies through it and exists
above.  For $\sigma_i^{-1}$ we are ``playing the movie backwards'', i.e.,
the ring of the strand going under comes from above and exits below the ring
of the ``over'' strand.

Let ``the signed $w$ braid group'', $\swB_n$, be the group
of horizontal flying rings where both fly-through 
directions are allowed. This introduces a ``sign'' for
each crossing $\sigma_i$:
\begin{center}
 \input figs/FlyingRings2.pstex_t
\end{center}
In other words, $\swB_n$ is generated by $s_i$, 
$\sigma_{i+}$ and $\sigma_{i-}$, for $i=1,...,n-1$. Check that in $\swB_n$
$\sigma_{i-}=s_i\sigma_{i+}^{-1}s_i$, and this, along with
the other obvious relations implies $\swB_n \cong \wB_n$.
\end{exercise}

\subsubsection{Certain Ribbon Tubes in $\bbR^4$} \label{subsubsec:ribbon}
With
time as the added dimension, a flying ring in $\bbR^3$ traces a tube
(an annulus) in $\bbR^4$, as shown in the picture below:
\[ \input figs/RibbonTubes.pstex_t \]
Note that we adopt here the drawing conventions of Carter and
Saito~\cite{CarterSaito:KnottedSurfaces} --- we draw surfaces as if they
were projected from $\bbR^4$ to $\bbR^3$, and we cut them open whenever
they are ``hidden'' by something with a higher fourth coordinate.

Note
also that the tubes we get in $\bbR^4$ always bound natural 3D
``solids'' --- their ``insides'', in the pictures above. These solids
are disjoint in the case of $s_i$ and have a very specific kind
of intersection in the case of $\sigma_i$ --- these are transverse
intersections with no triple points, and their inverse images are a
meridional disk on the ``thin'' solid tube and an interior disk on the
``thick'' one. By analogy with the case of ribbon knots and ribbon
singularities in $\bbR^3$ (e.g.~\cite[Chapter V]{Kauffman:OnKnots}) and
following Satoh~\cite{Satoh:RibbonTorusKnots}, we call this kind if
intersections of solids in $\bbR^4$ ``ribbon singularities'' and thus, our
tubes in $\bbR^4$ are always ``ribbon tubes''.

\subsubsection{Basis Conjugating Automorphisms of $F_n$}
\label{subsubsec:McCool}
Let
$\glos{F_n}$ be the free (non-abelian) group with generators
$\glos{\xi_1,\ldots,\xi_n}$. Artin's theorem (Theorems 15 and 16
of~\cite{Artin:TheoryOfBraids}) says that the (usual) braid
group $\uB_n$ (equivalently, the subgroup of $\wB_n$ generated by
the $\sigma_i$'s) has a faithful right action on $F_n$. In other
words, $\uB_n$ is isomorphic to a subgroup $H$ of $\Autop(F_n)$
(the group of automorphisms of $F_n$ with opposite multiplication, i.e.,
$\psi_1\psi_2:=\psi_2\circ\psi_1$). Precisely, using $(\xi,
B)\mapsto\xi\glos{\sslash}B$ to denote the right action of $\Autop(F_n)$ on
$F_n$, the subgroup $H$ consists of those automorphisms $B\colon F_n\to F_n$
of $F_n$ that satisfy the following two conditions:
\begin{enumerate}
\item $B$ maps any generator $\xi_i$ to a
  conjugate of a generator (possibly different). That is, there is a
  permutation $\beta\in S_n$ and elements $a_i\in F_n$ so that, for every $i$,
  \begin{equation} \label{eq:BasisConjugating}
    \xi_i \sslash B = a_i^{-1}\xi_{\beta (i)}a_i.
  \end{equation}
\item $B$ fixes the ordered product of the generators of $F_n$,
  \[ \xi_1\xi_2\cdots \xi_n \sslash B = \xi_1\xi_2\cdots \xi_n. \]
\end{enumerate}

McCool's theorem\footnote{Stricktly speaking, the main theorem
of~\cite{McCool:BasisConjugating} is about $\PwB_n$, yet it can easily be
restated for $\wB_n$.}~\cite{McCool:BasisConjugating} says that almost the same statement 
holds true\footnote{Though see Warning~\ref{warn:NoArtin}.} for the bigger group $\wB_n$:
namely, $\wB_n$ is isomorphic to the subgroup of $\Autop(F_n)$ consisting of
automorphisms satisfying only the first condition above.
So $\wB_n$ is precisely the group
of ``basis-conjugating'' automorphisms of the free group $F_n$,
the group of those automorphisms which map any ``basis element''
in $\{\xi_1,\ldots,\xi_n\}$ to a conjugate of a (possibly different)
basis element.

The relevant action is explicitly defined on the generators of $\wB_n$
and $F_n$ as follows (we state how each generator of $\wB_n$ acts on 
each generator of $F_n$, in each case omitting the generators of $F_n$ 
which are fixed under the action):
\begin{equation} \label{eq:ExplicitPsi}
  (\xi_i, \xi_{i+1})\sslash s_i = (\xi_{i+1}, \xi_i),
  \qquad
  (\xi_i, \xi_{i+1})\sslash \sigma_i = (\xi_{i+1},
    \xi_{i+1}\xi_i\xi_{i+1}^{-1}),
  \qquad
  \xi_j\sslash \sigma_{ij} = \xi_i\xi_j\xi_i^{-1}.
\end{equation} 
It is a worthwhile exercise to verify that $\sslash$ respects the
relations in the definition of $\wB_n$ and that the permutation $\beta$
in~\eqref{eq:BasisConjugating} is the skeleton $\varsigma(B)$.

There is a more conceptual description of $\sslash$, in terms of the
structure of $\wB_{n+1}$. Consider the inclusions
\begin{equation} \label{eq:inclusions}
  \wB_n \xhookrightarrow{\iota} \wB_{n+1} \xhookleftarrow{i_u} F_n.
\end{equation}

Here $\glos{\iota}$ is the inclusion of $\wB_n$ into $\wB_{n+1}$ by adding an
inert $(n+1)$st strand (it is injective as it has a well-defined
one sided inverse -- the deletion of the $(n+1)$st strand). 

\parpic[r]{$\input figs/xi.pstex_t$}
The inclusion $\glos{i_u}$ of the free group $F_n$ into $\wB_{n+1}$ is defined by
$i_u(\xi_i):=\sigma_{i,n+1}$.  The image $i_u(F_n)\subset\wB_{n+1}$ is the
set of all w-braids whose first $n$ strands are straight and vertical,
and whose $(n+1)$-st strand wanders among the first $n$ strands mostly
virtually (i.e., mostly using virtual crossings), occasionally slipping
under one of those $n$ strands, but never going over anything. 
It is easier to see that
this is indeed injective using the
``flying rings'' picture of Section~\ref{subsubsec:FlyingRings}. The image
$i_u(F_n)\subset\wB_{n+1}$ can be interpreted as the fundamental group
of the complement in $\bbR^3$ of $n$ stationary rings (which is indeed
$F_n$) --- in $i_u(F_n)$ the only ring in motion is the last, and it only
goes under, or ``through'', other rings, so it can be replaced by a point
object whose path is an element of the fundamental group. The injectivity
of $i_u$ follows from this geometric picture.

\parpic[r]{$\input figs/Bgamma.pstex_t$} \picskip{4}
One may explicitly verify that $i_u(F_n)$ is normalized by $\iota(\wB_n)$
in $\wB_{n+1}$ (that is, the set $i_u(F_n)$ is preserved by conjugation
by elements of $\iota(\wB_n)$). Thus, the following definition (also shown
as a picture on the right) makes sense, for \linebreak $B\in\wB_n\subset\wB_{n+1}$
and for $\gamma\in F_n\subset\wB_{n+1}$:

\begin{equation} \label{eq:ConceptualPsi}
  \gamma\sslash B := i_u^{-1}(B^{-1}\gamma B)
\end{equation}

It is a worthwhile exercise to recover the explicit formulae
in~\eqref{eq:ExplicitPsi} from the above definition.

\begin{warning} \label{warn:NoArtin} People familiar with the Artin story for
ordinary braids should be warned that even though $\wB_n$ acts on $F_n$ and
the action is induced from the inclusions in~\eqref{eq:inclusions} in much
of the same way as the Artin action is induced by inclusions $\uB_n
\xhookrightarrow{\iota} \uB_{n+1} \xhookleftarrow{i} F_n$, there are also
some differences, and some further warnings apply:
\begin{myitemize}
\item In the ordinary Artin story, $i(F_n)$ is the set of braids in
$\uB_{n+1}$ whose first $n$ strands are unbraided (that is, whose image in
$\uB_n$ via ``dropping the last strand'' is the identity). This is not true
for w-braids. For w-braids, in $i_u(F_n)$ the last strand always goes
``under'' all other strands (or just virtually crosses them), but never
``over''.
\item Thus, unlike the isomorphism $\PuB_{n+1}\cong \PuB_n\ltimes F_n$,
it is not true that $\PwB_{n+1}$ is isomorphic to $\PwB_n\ltimes F_n$.
\item The OC relation imposed in $\wB$ breaks
the symmetry between overcrossings and undercrossings. Thus, let
$i_o\colon F_n\to\wB_n$ be the ``opposite'' of $i_u$, mapping into braids in
which the last strand is always ``over'' or virtual. Then $i_o$ is not
injective (its image is in fact abelian) and its image is not normalized
by $\iota(\wB_n)$. So there is no ``second'' action of $\wB_n$ on $F_n$
defined using $i_o$.
\item For v-braids, both $i_u$ and $i_o$ are injective and there are two
actions of $\vB_n$ on $F_n$ --- one defined by first projecting into
w-braids, and the other defined by first projecting into v-braids modulo
``undercrossings commute''. Yet v-braids contain more information than
these two actions can see. The ``Kishino'' v-braid below, for example,
is visibly trivial if either overcrossings or undercrossings are made to
commute, yet by computing its Kauffman bracket we know it is non-trivial
as a v-braid~\cite[``The Kishino Braid'']{WKO}:
\[
  \begin{picture}(0,0)%
\includegraphics{figs/KishinoBraid.pstex}%
\end{picture}%
%
%  pstex_opts: -m 0.7 
%
\setlength{\unitlength}{2763sp}%
\begingroup\makeatletter\ifx\SetFigFont\undefined%
\gdef\SetFigFont#1#2#3#4#5{%
  \reset@font\fontsize{#1}{#2pt}%
  \fontfamily{#3}\fontseries{#4}\fontshape{#5}%
  \selectfont}%
\fi\endgroup%
\begin{picture}(7074,1027)(-11,-176)
\put(976,-136){\makebox(0,0)[b]{\smash{{\SetFigFont{7}{8.4}{\rmdefault}{\mddefault}{\updefault}{\color[rgb]{0,0,0}$a$}%
}}}}
\put(6076,-136){\makebox(0,0)[b]{\smash{{\SetFigFont{7}{8.4}{\rmdefault}{\mddefault}{\updefault}{\color[rgb]{0,0,0}$b$}%
}}}}
\end{picture}%
 \quad \left(\parbox{1.6in}{\footnotesize
    The commutator $ab^{-1}a^{-1}b$ of v-braids $a,b$ annihilated by OC/UC,
    respectively, with a minor cancellation.
  }\right)
\]
\end{myitemize}
\end{warning}

\begin{problem} \label{prob:wCombing}
Are $\PvB_n$ and $\PwB_n$ semi-direct products of free groups?
For $\PuB_n$, this
is the well-known ``combing of braids'' and it follows from $\PuB_n\cong
\PuB_{n-1}\ltimes F_{n-1}$ and induction.
\end{problem}

\begin{remark} \label{rem:GutierrezKrstic}
Note that Guti\'errez and Krsti\'c~\cite{GutierrezKrstic:NormalForms}
have found ``normal forms'' for the elements of $\PwB_n$, yet they
do not decide whether $\PwB_n$ is ``automatic'' in the sense
of~\cite{Epstein:WordProcessing}.
\end{remark}

\draftcut \subsection{Finite Type Invariants of v-Braids and w-Braids}
\label{subsec:FT4Braids}

Just as we had two definitions for v-braids (and thus, for w-braids)
in Section~\ref{subsec:VirtualBraids}, we will give two 
equivalent developments of the theory of finite type invariants of
v-braids and w-braids --- a pictorial/topological version in
Section~\ref{subsubsec:FTPictorial}, and a more abstract algebraic version
in Section~\ref{subsubsec:FTAlgebraic}.

\subsubsection{Finite Type Invariants, the Pictorial Approach}
\label{subsubsec:FTPictorial}

In
the standard theory of finite type invariants of knots (also known as
Vassiliev or Goussarov-Vassiliev invariants)~\cite{Vassiliev:CohKnot, Goussarov:nEquivalence,
Bar-Natan:OnVassiliev, Bar-Natan:EMP}
one progresses from the definition of finite type via iterated
differences to chord diagrams and weight systems, to $4T$ (and other)
relations, to the definition of universal finite type invariants,
and beyond. The exact same progression (with different objects
playing similar roles, and sometimes, when yet insufficiently
studied, with the last step or two missing) is also seen in the theories
of finite type invariants of braids~\cite{Bar-Natan:Braids},
3-manifolds~\cite{Ohtsuki:IntegralHomology, LeMurakamiOhtsuki:Universal,
Le:UniversalIHS}, virtual knots~\cite{GoussarovPolyakViro:VirtualKnots,
Polyak:ArrowDiagrams} and of several other classes of objects. We thus
assume that the reader has familiarity with these basic ideas, and we
only indicate briefly how they are implemented in the case of v-braids
and w-braids.

\begin{figure}
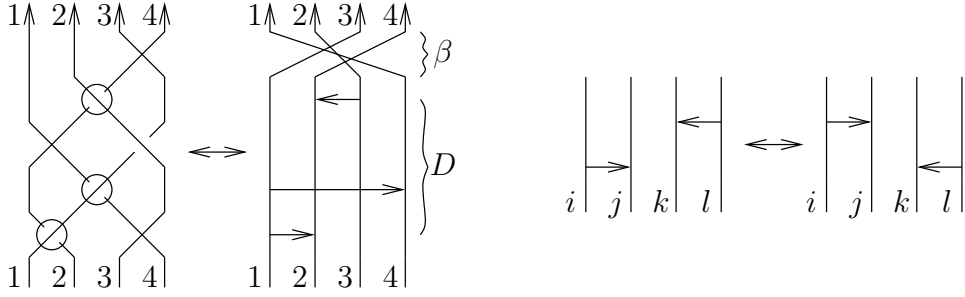

\[ \input figs/Dvh1.pstex_t \]
  \caption{
    On the left, a 3-singular v-braid and its corresponding 3-arrow
    diagram. A self-explanatory algebraic notation for this arrow
    diagram is $(\glos{a_{12}a_{41}a_{23}},\,3421)$. 
    Note that we regard arrow diagrams as graph-theoretic
    objects, and hence, the two arrow diagrams on the right, whose
    underlying graphs are the same, are regarded as equal. In
    algebraic notation this means that we always impose the relation
    $a_{ij}a_{kl}=a_{kl}a_{ij}$ when the indices $i$, $j$, $k$, and $l$
    are all distinct.
  } \label{fig:Dvh1}
\end{figure}

Much
like the formula $\doublepoint\to\overcrossing-\undercrossing$ of
the Vassiliev-Goussarov fame, given a v-braid invariant $\glos{V}\colon
\vB_n\to A$ valued in some abelian group $A$, we extend it to
``singular'' v-braids, i.e., braids that contain ``semi-virtual crossings''
like $\glos{\semivirtualover}$ and $\glos{\semivirtualunder}$ using the
formulae $V(\semivirtualover):=V(\overcrossing)-V(\virtualcrossing)$
and $V(\semivirtualunder):=V(\undercrossing)-V(\virtualcrossing)$
(see~\cite{GoussarovPolyakViro:VirtualKnots, Polyak:ArrowDiagrams,
Bar-NatanHalachevaLeungRoukema:v-Dims}). We
say that ``$V$ is of type $m$'' if its extension vanishes on singular
v-braids having more than $m$ semi-virtual crossings. Up to invariants
of lower type, an invariant of type $m$ is determined by its ``weight
system'', which is a functional $W=\glos{W_m}(V)$ defined on ``$m$-singular
v-braids modulo $\overcrossing=\virtualcrossing=\undercrossing$''. Let
us denote the vector space of all formal linear combinations of such
equivalence classes by $\glos{\calG_m}\calD^v_n$. Much as $m$-singular knots
modulo $\overcrossing=\undercrossing$ can be identified with chord
diagrams, the basis elements of $\calG_m\calD^v_n$ can be identified with
pairs $(D,\beta)$, where $D$ is a horizontal arrow diagram and $\beta$
is a ``skeleton permutation'', see Figure~\ref{fig:Dvh1}.

We assemble the spaces $\calG_m\calD^v_n$ together to form a single
graded space, $\glos{\calD^v_n}:=\oplus_{m=0}^\infty\calG_m\calD^v_n$. Note that
throughout this paper, whenever we write an infinite direct sum,
we automatically complete it. Thus, in $\calD^v_n$ we allow infinite sums
with one term in each homogeneous piece $\calG_m\calD^v_n$, in particular,
exponential-like sums will be heavily used.

\begin{figure}
\[ \input figs/6T.pstex_t \]
\[
  a_{ij}a_{ik}+a_{ij}a_{jk}+a_{ik}a_{jk}
  = a_{ik}a_{ij}+a_{jk}a_{ij}+a_{jk}a_{ik}
\]
\[
  \text{or}\qquad
  [a_{ij}, a_{ik}] + [a_{ij}, a_{jk}] + [a_{ik}, a_{jk}] = 0
\]
\caption{The $6T$ relation. Standard knot theoretic conventions apply ---
only the relevant parts of each diagram is shown; in reality each diagram
may have further vertical strands and horizontal arrows, provided the
extras are the same in all 6 diagrams. Also, the vertical strands are in no
particular order --- other valid $6T$ relations are obtained when those
strands are permuted in other ways.} \label{fig:6T}
\end{figure}

\begin{figure}
\[ \begin{array}{ccc}
  \input figs/TC.pstex_t & \qquad & \input figs/4TArrow.pstex_t \\
  a_{ij}a_{ik} = a_{ik}a_{ij} &&
  a_{ij}a_{jk} + a_{ik}a_{jk} = a_{jk}a_{ij} + a_{jk}a_{ik} \\
  \text{or} \quad [a_{ij}, a_{ik}] = 0 &&
  \text{or} \quad [a_{ij} + a_{ik}, a_{jk}] = 0
\end{array} \]
\caption{The TC and the $\protect\aft$ relations.} \label{fig:TCand4T}
\end{figure}

In the standard finite-type theory for knots, weight
systems always satisfy the $4T$ relation, and are therefore
functionals on $\calA:=\calD/4T$. Likewise, in the case of
v-braids, weight systems satisfy the ``$\glos{6T}$ relation''
of~\cite{GoussarovPolyakViro:VirtualKnots, Polyak:ArrowDiagrams,
Bar-NatanHalachevaLeungRoukema:v-Dims},
shown in Figure~\ref{fig:6T}, and are therefore functionals on
$\glos{\calA^v_n}:=\calD^v_n/6T$. In the case of w-braids, the OC 
relation~\eqref{eq:OvercrossingsCommute} implies the
``tails commute'' (\glost{TC}) relation on the level of arrow diagrams,
and in the presence of the TC relation, two of the terms in the $6T$
relation drop out, and what remains is the ``$\glos{\aft}$'' relation. These
relations are shown in Figure~\ref{fig:TCand4T}. Thus, weight systems
of finite type invariants of w-braids are linear functionals on
$\glos{\calA^w_n}:=\calD^v_n/TC,\aft$.

The next question that arises is whether we have already found {\em
all} the relations that weight systems always satisfy. More precisely,
given a degree $m$ linear functional on $\calA^v_n=\calD^v_n/6T$ (or
on $\calA^w_n=\calD^v_n/TC,\aft$), is it always the weight system
of some type $m$ invariant $V$ of v-braids (or w-braids)? As in every
other theory of finite type invariants, the answer to this question
is affirmative if and only if there exists a ``universal finite type
invariant'' (or simply, an ``expansion'') of v-braids (or w-braids):

\begin{definition} \label{def:vwbraidexpansion} An expansion
for v-braids (or w-braids) is an invariant $Z\colon \vB_n\to\calA^v_n$
(or $Z\colon \wB_n\to\calA^w_n$) satisfying the following ``universality
condition'':
\begin{itemize}
\item If $B$ is an $m$-singular v-braid (or w-braid) and
$D\in\calG_m\calD^v_n$ is its underlying arrow diagram as in
Figure~\ref{fig:Dvh1}, then
\[ Z(B)=D+(\text{terms of degree\,}>m). \]
\end{itemize}
\end{definition}

Indeed if $Z$ is an expansion and $W\in\calG_m\calA^\star$,\footnote{$\calA^\star$
here denotes either $\calA^v_n$ or $\calA^w_n$, or in fact, any of many
similar spaces that we will discuss later on.} the universality condition
implies that $W\circ Z$ is a finite type invariant whose weight system
is $W$. To go the other way, if $(D_i)$ is a basis of $\calA$ consisting
of homogeneous elements, if $(W_i)$ is the dual basis of $\calA^\star$ and
$(V_i)$ are finite type invariants whose weight systems are the $W_i$'s,
then $Z(B):=\sum_iD_iV_i(B)$ defines an expansion.

In general, constructing a universal finite type invariant is a
hard task. For knots, one uses either the Kontsevich integral or
perturbative Chern-Simons theory (also known as ``configuration
space integrals''~\cite{BottTaubes:SelfLinking} or ``tinker-toy
towers''~\cite{Thurston:IntegralExpressions}) or the rather fancy
algebraic theory of ``Drinfel'd associators'' (a summary of all those
approaches is at~\cite{Bar-NatanStoimenow:Fundamental}). For homology
spheres, this is the ``LMO invariant''~\cite{LeMurakamiOhtsuki:Universal,
Le:UniversalIHS} (also the ``\AA{}rhus
integral''~\cite{Bar-NatanGaroufalidisRozanskyThurston:Aarhus}). For
v-braids, we still don't know if an expansion exists. In contrast, as we shall see
below, the construction of an expansion for w-braids is quite easy.

\subsubsection{Finite Type Invariants, the Algebraic Approach}
\label{subsubsec:FTAlgebraic}
For any group $G$, one can form the group algebra ${\mathbb Q}G$ by
allowing formal linear combinations of group elements and extending
multiplication linearly, where $\mathbb Q$ is the field of rational
numbers\footnote{The definitions in this subsection make sense over $\bbZ$
as well, but the main result of the next subsection requires a field of
characteristic $0$. For simplicity of notation we stick with $\bbQ$.}.
The {\it augmentation ideal} is the ideal generated by
differences, or equivalently, the set of linear combinations of group
elements whose coefficients sum to zero:
\[ \glos{\calI} := \left\{\sum_{i=1}^k a_ig_i\colon
  a_i \in {\mathbb Q}, g_i \in G, \sum_{i=1}^k a_i=0\right\}.
\]
Powers of the augmentation ideal provide a filtration of the group
algebra. We denote by $\glos{\calA(G)}:= \bigoplus_{m\geq 0} \calI^m/\calI^{m+1}$
the associated graded space corresponding to this filtration.

\begin{definition}\label{def:grpexpansion} An
expansion for the group
$G$ is a map $Z\colon G \to \calA(G)$, such that the linear extension
$Z\colon  {\mathbb Q}G \to \calA(G)$ is filtration preserving and
the induced map $$\gr Z\colon  (\gr {\mathbb Q}G=\calA(G)) \to (\gr
\calA(G)=\calA(G))$$ is the identity. An equivalent way to phrase this
is that the degree $m$ piece of $Z$ restricted to $\calI^m$ is the
projection onto $\calI^m/\calI^{m+1}$.

\begin{exercise}\label{ex:BraidsAlgApproach}
Verify that for the groups $\PvB_n$ and $\PwB_n$ the m-th power of the
augmentation ideal coincides with the span of all resolutions of
$m$-singular $v$- or $w$-braids (by a resolution we mean the formal
linear combination where each semivirtual crossing is replaced by
the appropriate difference of a virtual and a regular crossing, as in Figure \ref{fig:Dvh1}). Then
check that the notion of expansion defined above is the same as that of
Definition \ref{def:vwbraidexpansion}, restricted to pure braids.
\end{exercise}

Finally, note the functorial nature of the construction above. What we have 
described is a functor from the category of groups to the category of graded
algebras, called {\em projectivization}\footnote{We use this name to
distinguish the associated graded with respect to this particular filtration,
which will be a repeating theme in \cite{Bar-NatanDancso:WKO2}.}
$\proj\colon  Grp \to GrAlg$, which assigns to each group
$G$ the graded algebra $\calA(G)$. To each homomorphism $\phi\colon  G \to H$,
$\proj$ assigns
the induced map $$\gr \phi\colon  (\calA(G)=\gr {\mathbb Q}G) \to (\calA(H)= \gr {\mathbb Q}H).$$
 
\end{definition}

\draftcut \subsection{Expansions for w-Braids}\label{subsec:wBraidExpansion}

The space $\calA^w_n$ of arrow diagrams on $n$ strands is an associative
algebra in an obvious manner: if the permutations underlying two arrow
diagrams are the identity permutations, then we simply juxtapose the diagrams.
Otherwise we ``slide'' arrows through permutations in the obvious manner
--- if $\tau$ is a permutation, we declare that $\tau a_{(\tau i)(\tau
j)}=a_{ij}\tau$. Instead of seeking an expansion $\wB_n\to\calA^w_n$, we
set the bar a little higher and seek a ``homomorphic expansion'':

\begin{definition} \label{def:Universallity} A homomorphic expansion
$Z\colon \wB_n\to\calA^w_n$ is an expansion that carries products in $\wB_n$
to products in $\calA^w_n$.
\end{definition}

To find a homomorphic expansion, we just need to define it
on the generators of $\wB_n$ and verify that it satisfies the
relations defining $\wB_n$ and the universality condition.
Following~\cite[Section~5.3]{BerceanuPapadima:BraidPermutation}
and~\cite[Section~8.1]{AlekseevTorossian:KashiwaraVergne} we set
$Z(\virtualcrossing)=\virtualcrossing$ (that is, a transposition in
$\wB_n$ gets mapped to the same transposition in $\calA^w_n$, adding no
arrows) and $Z(\overcrossing)=\exp(\rightarrowdiagram)\virtualcrossing$.
(Reacall that we work in the degree completion.)
This last formula is important so deserves to be magnified, explained
and replaced by some new notation:

\begin{equation} \label{eq:reservoir}
  Z\left(\!\mathsize{\Huge}{\overcrossing}\!\right)\! =
  \exp\left(\!\mathsize{\Huge}{\rightarrowdiagram}\!\right)
    \cdot\mathsize{\Huge}{\virtualcrossing}
  = \input figs/ZIsExp.pstex_t+\ldots =: \input figs/ArrowReservoir.pstex_t.
\end{equation}

Thus the new notation $\overset{e^a}{\longrightarrow}$ stands
for an ``exponential reservoir'' of parallel arrows, much like
$e^a=1+a+aa/2+aaa/3!+\ldots$ is a ``reservoir'' of $a$'s. With
the obvious interpretation for $\overset{e^{-a}}{\longrightarrow}$
(that is, the $-$ sign indicates that the terms should have alternating signs,
as in $e^{-a}=1-a+a^2/2-a^3/3!+\ldots$), the second Reidemeister move
$\overcrossing\undercrossing=1$ forces that we set
\[ Z\left(\mathsize{\Huge}{\undercrossing}\right) =
  \mathsize{\Huge}{\virtualcrossing}
  \cdot\exp\left(-\mathsize{\Huge}{\rightarrowdiagram}\right)
  = \begin{picture}(0,0)%
\includegraphics{figs/NegReservoir1.pstex}%
\end{picture}%
%
%  pstex_opts: -m 0.75 
%
\setlength{\unitlength}{2960sp}%
\begingroup\makeatletter\ifx\SetFigFont\undefined%
\gdef\SetFigFont#1#2#3#4#5{%
  \reset@font\fontsize{#1}{#2pt}%
  \fontfamily{#3}\fontseries{#4}\fontshape{#5}%
  \selectfont}%
\fi\endgroup%
\begin{picture}(739,1374)(1281,-523)
\put(1651,239){\makebox(0,0)[b]{\smash{{\SetFigFont{9}{10.8}{\rmdefault}{\mddefault}{\updefault}{\color[rgb]{0,0,0}$e^{-a}$}%
}}}}
\end{picture}%
 = \input figs/NegReservoir2.pstex_t.
\]

\begin{theorem} \label{thm:RInvariance} The above formulae define
an invariant $Z\colon \wB_n\to\calA^w_n$ (that is, $Z$ satisfies all the
defining relations of $\wB_n$). The resulting $Z$ is a homomorphic
expansion (that is, it satisfies the universality property of
Definition~\ref{def:Universallity}).
\end{theorem}

\begin{proof} Following~\cite{BerceanuPapadima:BraidPermutation,
AlekseevTorossian:KashiwaraVergne}: for the invariance of $Z$, the
only interesting relations to check are the Reidemeister 3 relation
of~\eqref{eq:sigmaRels} and the OC relation
of~\eqref{eq:OC}. For Reidemeister 3, we have
\[ \input figs/R3Left.pstex_t
  = e^{a_{12}}e^{a_{13}}e^{a_{23}}\tau
  \overset{1}{=} e^{a_{12}+a_{13}}e^{a_{23}}\tau
  \overset{2}{=} e^{a_{12}+a_{13}+a_{23}}\tau,
\]
where $\tau$ is the permutation $321$ and equality 1 holds because
$[a_{12},a_{13}]=0$ by a TC relation and equality 2 holds
because $[a_{12}+a_{13}, a_{23}]=0$ by a $\aft$ relation.
Likewise, again using TC and $\aft$,
\[ \input figs/R3Right.pstex_t
  = e^{a_{23}}e^{a_{13}}e^{a_{12}}\tau
  = e^{a_{23}}e^{a_{13}+a_{12}}\tau
  = e^{a_{23}+a_{13}+a_{12}}\tau,
\]
and so Reidemeister 3 holds. An even simpler proof using just the TC 
relation shows that the OC relation also holds.
Finally, since $Z$ is homomorphic, it is enough to check the universality
property at degree $1$, where it is very easy:
\[ Z\left(\mathsize{\Huge}{\semivirtualover}\right) =
  \exp\left(\mathsize{\Huge}{\rightarrowdiagram}\right)
    \cdot\mathsize{\Huge}{\virtualcrossing}
    - \mathsize{\Huge}{\virtualcrossing}
  = \mathsize{\Huge}{\rightarrowdiagram}\cdot\mathsize{\Huge}{\virtualcrossing}
    + (\text{terms of degree\,}>1),
\]
and a similar computation manages the $\semivirtualunder$ case. \qed
\end{proof}

\begin{remark} \label{rem:YangBaxter} Note that the main ingredient
of the above proof was to show that \linebreak $\glos{R}:=Z(\sigma_{12})=e^{a_{12}}$
satisfies the famed Yang-Baxter equation,
\[ R^{12}R^{13}R^{23} = R^{23}R^{13}R^{12}, \]
where $R^{ij}$ means ``place $R$ on strands $i$ and $j$''.
\end{remark}

\draftcut
\subsection{Some Further Comments} \label{subsec:bcomments}
\subsubsection{Compatibility with Braid Operations}
  \label{subsubsec:BraidCompatibility}
As with
any new gadget, we would like to know how compatible the expansion
$Z$ of the previous section is with the gadgets we already have; namely,
with various operations that are available on w-braids and with the action
of w-braids on the free group $F_n$, see Section~\ref{subsubsec:McCool}.

\parpic[r]{$\xymatrix{
  \wB_n \ar[r]^\theta \ar[d]_Z   & \wB_n \ar[d]^Z    \\
  \calA^w_n \ar[r]_\theta        & \calA^w_n
    \ar@{}[ul]|{\text{\huge$\circlearrowleft$}}
}$}
\paragraph{$Z$ is Compatible with Braid Inversion} \label{par:theta}
Let $\theta$ denote both the ``braid inversion'' operation
$\glos{\theta}\colon \wB_n\to\wB_n$ defined by $B\mapsto B^{-1}$ and
the ``antipode'' anti-automorphism $\theta\colon \calA^w_n\to\calA^w_n$
defined by mapping permutations to their inverses and arrows to their
negatives (that is, $a_{ij}\mapsto-a_{ij}$). Then the diagram on the
right commutes.

\pagebreak[2]

\parpic[r]{$\xymatrix{
  \wB_n \ar[r]^<>(0.5)\Delta \ar[d]_Z   & \wB_n\times\wB_n \ar[d]^{Z\times Z} \\
  \calA^w_n \ar[r]_<>(0.5)\Delta        & \calA^w_n\otimes\calA^w_n
    \ar@{}[ul]|{\text{\huge$\circlearrowleft$}}
}$}
\paragraph{Braid Cloning and the Group-Like Property} \label{par:Delta}
Let $\glos{\Delta}$ denote both the ``braid cloning''
operation $\Delta\colon \wB_n\to\wB_n\times\wB_n$ defined by
$B\mapsto (B,B)$ and the ``co-product'' algebra morphism \linebreak
$\Delta\colon \calA^w_n\to\calA^w_n\otimes\calA^w_n$ defined by cloning
permutations (that is, $\tau\mapsto\tau\otimes\tau$) and by treating
arrows as primitives (that is, \linebreak $a_{ij}\mapsto a_{ij}\otimes 1+1\otimes
a_{ij}$). Then the diagram on the right commutes. In formulae, this is
$\Delta(Z(B))=Z(B)\otimes Z(B)$, which is the statement ``$Z(B)$ is
group-like''.

\parpic[r]{$\xymatrix{
  \wB_n \ar[r]^<>(0.5)\iota \ar[d]_Z	& \wB_{n+1} \ar[d]^Z	\\
  \calA^w_n \ar[r]_<>(0.5)\iota	& \calA^w_{n+1}
    \ar@{}[ul]|{\text{\huge$\circlearrowleft$}}
}$}
\paragraph{Strand Insertions} \label{par:iota}
Let $\iota\colon \wB_n\to\wB_{n+1}$ be an operation of ``inert strand
insertion''. Given $B\in\wB_n$, the resulting $\iota B\in\wB_{n+1}$
will be $B$ with one strand $S$ added at some location chosen in
advance --- to the left of all existing strands, or to the right, or
starting from between the 3rd and the 4th strand of $B$ and ending
between the 6th and the 7th strand of $B$; when adding $S$, add it
``inert'', so that all crossings on it are virtual (this is well
defined). There is a corresponding inert strand addition operation
$\iota\colon \calA^w_n\to\calA^w_{n+1}$, obtained by adding a strand at the
same location as for the original $\iota$ and adding no arrows. It is
easy to check that $Z$ is compatible with $\iota$; namely, that the
diagram on the right is commutative.

\parpic[r]{$\xymatrix{
  \wB_n \ar[r]^<>(0.5){d_k} \ar[d]_Z   & \wB_{n-1} \ar[d]^Z    \\
  \calA^w_n \ar[r]_<>(0.5){d_k}        & \calA^w_{n-1}
    \ar@{}[ul]|{\text{\huge$\circlearrowleft$}}
}$}
\paragraph{Strand Deletions} \label{par:deletions} Given $1 \leq k \leq n$, 
let $\glos{d_k}\colon \wB_n\to\wB_{n-1}$ be the operation of
``removing the $k$th strand''.  This operation induces a homonymous
operation $d_k\colon \calA^w_n\to\calA^w_{n-1}$: if $D\in\calA^w_n$ is an
arrow diagram, then $d_kD$ is $D$ with its $k$th strand removed if no arrows
in $D$ start or end on the $k$th strand, and it is $0$ otherwise. It
is easy to check that $Z$ is compatible with $d_k$; namely, that the
diagram on the right is commutative.\footnote{In \cite{Bar-NatanDancso:WKO2} we'll 
say that ``$d_k\colon \wB_n\to\wB_{n-1}$''
is an algebraic structure made of two spaces ($\wB_n$ and $\wB_{n-1}$),
two binary operations (braid composition in $\wB_n$ and in $\wB_{n-1}$),
and one unary operation, $d_k$. After projectivization we get the
algebraic structure $d_k\colon \calA^w_n\to\calA^w_{n-1}$ with $d_k$
as described above, and an alternative way of stating our assertion is
to say that $Z$ is a morphism of algebraic structures. A similar remark
applies (sometimes in the negative form) to the other operations discussed
in this section.}

\parpic[r]{$\xymatrix{
  F_n \ar@{}[r]|{\mathsize{\Huge}{\actsonright}} \ar[d]_Z & \wB_n \ar[d]^Z \\
  \FA_n \ar@{}[r]|{\mathsize{\Huge}{\actsonright}} & \calA^w_n 
    \ar@{}[ul]|{\text{\huge$\circlearrowleft$}}
}$}
\paragraph{Compatibility with the Action on $F_n$} \label{par:action}
Let $\glos{\FA_n}$ denote the (degree-completed) free, associative (but
not commutative) algebra on the generators $\glos{x_1,\dots,x_n}$. Then
there is an ``expansion'' $Z\colon F_n\to \FA_n$ defined by $\xi_i\mapsto
e^{x_i}$ (see~\cite{Lin:Expansions} and the related ``Magnus Expansion''
of~\cite{MagnusKarrasSolitar:CGT}). Also, there is a right action\footnote{In the language of
\cite{Bar-NatanDancso:WKO2}, we will say that $\FA_n=\proj F_n$ and
that when the actions involved are regarded as instances of the algebraic
structure ``one monoid acting on another'', we have that \linebreak
$\left(\FA_n\actsonright\calA^w_n\right)=\proj\left(F_n\actsonright
\wB_n\right)$.} of
$\calA^w_n$ on $\FA_n$ defined on generators by $x_i\tau=x_{\tau i}$
for $\tau\in S_n$ and by $x_ja_{ij}=[x_i,x_j]$ and $x_ka_{ij}=0$ for
$k\neq j$ and extended by the Leibniz rule to the rest of $\FA_n$ and
then multiplicatively to the rest of $\calA^w_n$.

\begin{exercise}  Use the definition of the action in
\eqref{eq:ConceptualPsi} and the commutative diagrams of paragraphs
\ref{par:theta}, \ref{par:Delta} and~\ref{par:iota} to show that the
diagram of paragraph~\ref{par:action} is also commutative.
\end{exercise}

\pagebreak[2]

\parpic[r]{$\begin{array}{c}
  \input figs/StrandDoubling.pstex_t \\
  \xymatrix{
    \wB_n \ar[r]^<>(0.5){u_k} \ar[d]_Z   & \wB_{n+1} \ar[d]^Z    \\
    \calA^w_n \ar[r]_<>(0.5){u_k}        & \calA^w_{n+1}
      \ar@{}[ul]|{\text{\huge$\not\circlearrowleft$}}
  }
\end{array}$}
\paragraph{Unzipping a Strand} \label{par:unzip} Given $k$ between $1$ and
$n$, let $\glos{u_k}\colon \wB_n\to\wB_{n+1}$ the operation of ``unzipping
the $k$th strand'', briefly defined on the right\footnote{Unzipping
a knotted zipper turns a single band into two parallel ones. This
operation is also known as ``strand doubling'', but for compatibility with
operations that will be introduced later, we prefer ``unzipping''.}. The
induced operation $u_k\colon \calA^w_n\to\calA^w_{n+1}$ is also shown on
the right --- if an arrow starts (or ends) on the strand being doubled,
it is replaced by a sum of two arrows that start (or end) on either
of the two ``daughter strands'' (and this is performed for each arrow
independently; so if there are $t$ arrows touching the $k$th strands in
a diagram $D$, then $u_kD$ will be a sum of $2^t$ diagrams).

In some sense, much of this current series of papers as well as
the works of Kashiwara and Vergne~\cite{KashiwaraVergne:Conjecture}
and Alekseev and Torossian~\cite{AlekseevTorossian:KashiwaraVergne}
are about coming to grips with the fact that $Z$ is {\bf not}
compatible with $u_k$ (that the diagram on the right is {\bf not}
commutative). Indeed, let $x:=a_{13}$ and $y:=a_{23}$ be as on the
right, and let $s$ be the permutation $21$ and $\tau$ the permutation
$231$.  Then $d_1Z(\overcrossing)=d_1(e^{a_{12}}s)=e^{x+y}\tau$
while $Z(d_1\overcrossing)=e^ye^x\tau$. So the failure of $d_1$
and $Z$ to commute is the ill-behaviour of the exponential function
when its arguments are not commuting, which is measured by the
BCH formula, central to both~\cite{KashiwaraVergne:Conjecture}
and~\cite{AlekseevTorossian:KashiwaraVergne}.

\subsubsection{Power and Injectivity} The following theorem is due to
Berceanu and
Papadima~\cite[Theorem~5.4]{BerceanuPapadima:BraidPermutation}; a variant of
this theorem are also true for ordinary braids~\cite{Kohno:deRham, Bar-Natan:Homotopy,
HabeggerMasbaum:Milnor}, and can be proven by similar means:

\begin{theorem} $Z\colon \wB_n\to\calA^w_n$ is injective. In other words, finite
type invariants separate w-braids.
\end{theorem}

\begin{proof}The statement follows immediately from the faithfulness of
the action $F_n\actsonright\wB_n$, from the compatibility
of $Z$ with this action, and from the injectivity
of $Z\colon F_n\to\FA_n$ (the latter is well known, see
e.g.~\cite[Theorem~5.6]{MagnusKarrasSolitar:CGT}\footnote{Though notice
that we use $\xi_i\mapsto e^{x_i}$ whereas
\cite[Theorem~5.6]{MagnusKarrasSolitar:CGT} uses $\xi_i\mapsto 1+x_i$. The
\cite{MagnusKarrasSolitar:CGT} injectivity proof holds in our case
just as well.} and~\cite{Lin:Expansions}). Indeed, if $B_1$ and $B_2$
are w-braids and $Z(B_1)=Z(B_2)$, then $Z(\xi)Z(B_1)=Z(\xi)Z(B_2)$ for
any $\xi\in F_n$. Therefore $\forall\xi\, Z(\xi\sslash B_1)=Z(\xi\sslash
B_2)$, therefore $\forall\xi\,\xi\sslash B_1=\xi\sslash B_2$, therefore
$B_1=B_2$. \qed
\end{proof}

\begin{remark} Apart from the easy fact that $\calA^w_n$ can be computed
degree by degree in exponential time, we do not know a simple formula for
the dimension of the degree $m$ piece of $\calA^w_n$ or a natural basis of
that space. This compares unfavourably with the situation for ordinary
braids (see e.g.~\cite{Bar-Natan:Braids}). Also compare with
Problem~\ref{prob:wCombing} and with Remark~\ref{rem:GutierrezKrstic}.
\end{remark}

\subsubsection{Uniqueness} There is certainly not a unique expansion for
w-braids --- if $Z_1$ is an expansion and and $P$ is any degree-increasing
linear map $\calA^w\to\calA^w$ (a ``pollution'' map), then $Z_2:=(I+P)\circ
Z_1$ is also an expansion, where $I\colon \calA^w\to\calA^w$ is the
identity. But that's all, and if we require a bit more, even that
freedom disappears.

\begin{proposition} If $Z_{1,2}\colon \wB_n\to\calA^w_n$ are expansions then
there exists a degree-increasing linear map $P\colon \calA^w\to\calA^w$ so
that $Z_2=(I+P)\circ Z_1$.
\end{proposition}

\begin{proof} (Sketch).  Let $\widehat{\wB_n}$ be the unipotent completion
of $\wB_n$. That is, let $\bbQ\wB_n$ be the algebra of formal linear
combinations of w-braids, let $\calI$ be the ideal in $\bbQ\wB_n$ be the
ideal generated by $\semivirtualover=\overcrossing-\virtualcrossing$ and by
$\semivirtualunder=\virtualcrossing-\undercrossing$, and set
\[ \widehat{\wB_n}:=
  \underleftarrow{\lim}_{m\to\infty} \bbQ\wB_n \left/\calI^m\right..
\]
$\widehat{\wB_n}$ is filtered with
$\calF_m\widehat{\wB_n}:=\underleftarrow{\lim}_{m'>m} \calI^m
\left/\calI^{m'}\right..$ An ``expansion'' can be re-interpreted as an
``isomorphism of $\widehat{\wB_n}$ and $\calA^w_n$ as filtered vector
spaces''. Always, any two isomorphisms differ by an automorphism of the
target space, and that's the essence of $I+P$. \qed
\end{proof}

\begin{proposition} If $Z_{1,2}\colon \wB_n\to\calA^w_n$ are homomorphic
expansions that commute with braid cloning (Paragraph~\ref{par:Delta}) and
with strand insertion (Paragraph~\ref{par:iota}), then \linebreak $Z_1=Z_2$.
\end{proposition}

\begin{proof} (Sketch). A homomorphic expansion that commutes with strand
insertions is determined by its values on the generators $\overcrossing$,
$\undercrossing$ and $\virtualcrossing$ of $\wB_2$. Commutativity
with braid cloning (see Paragraph \ref{par:Delta}) implies that these values must be, up to permuting
the strands, group like: that is, the exponentials of primitives. But
the only primitives are $a_{12}$ and $a_{21}$, and one may easily
verify that there is only one way to arrange these so that $Z$
will respect $\virtualcrossing^2=\overcrossing\cdot\undercrossing=1$ and
$\semivirtualover\mapsto\rightarrowdiagram+(\text{higher degree terms})$. \qed
\end{proof}

\subsubsection{The Group of Non-Horizontal Flying Rings}
\label{subsubsec:NonHorRings}
Let
$\glos{Y_n}$ denote the space of all placements of $n$ numbered disjoint
oriented unlinked geometric circles in $\bbR^3$. Such a placement
is determined by the centres in $\bbR^3$ of the circles, the radii,
and a unit normal vector for each circle pointing in the positive
direction, so $\dim Y_n=3n+n+3n=7n$.  $S_n \ltimes \bbZ_2^n$ acts on
$Y_n$ by permuting the circles and mapping each circle to its image in
either an orientation-preserving or an orientation-reversing way. Let
$\glos{\tilde{Y}_n}$ denote the quotient $Y_n/S_n \ltimes \bbZ_2^n$.
The fundamental group $\pi_1(\tilde{Y}_n)$ can be thought of as the
``group of flippable flying rings''.  Without loss of generality, we
can assume that the basepoint is chosen to be a horizontal placement.
We want to study the relationship of this group to $\wB_n$.

\begin{theorem}
$\pi_1(\tilde{Y}_n)$ is a $\bbZ_2^n$-extension of $\wB_n$, generated
by $s_i$, $\sigma_{i}$ ($1\leq i \leq n-1)$, and $\glos{w_i}$ (``flips''),
for $1\leq i \leq n$; with the relations as above, and in addition:
\[
 w_i^2=1; \qquad w_iw_j=w_jw_i; \qquad w_js_i=s_iw_j \quad \text{when } \quad i\neq j, j+1;
\]
\[
 w_is_i=s_iw_{i+1}; \qquad w_{i+1}s_i=s_iw_i; \quad w_j\sigma_{i}=\sigma_{i}w_j \quad \text{if } \quad j \neq i, i+1;
\]
\[
w_{i+1}\sigma_{i}=\sigma_{i}w_{i}; \quad \text{yet } \quad
w_i\sigma_{i}=s_i\sigma_i^{-1}s_iw_{i+1}.
\]
\end{theorem}
The two most interesting flip relations in pictures:
\begin{equation}\label{eq:FlipRels}
  \raisebox{-10mm}{\input figs/FlipRels.pstex_t}
\end{equation}

\parpic[r]{\input{figs/FlippingRing.pstex_t}}
Instead of a proof, we provide some heuristics.
Since each circle starts out in a horizontal position and returns 
to a horizontal position, there is an integer number of 
``flips'' they do in between, these are the generators $w_i$, as 
shown on the right.

The first relation says that a double flip is homotopic to doing nothing.
Technically, there are two different directions of flips, and they are the 
same via this (non-obvious but true) relation. The rest of the first line is 
obvious: flips of different rings commute, and if
two rings fly around each other while another one flips, the order of these
events can be switched by homotopy. The second line says that, if two rings trade
places with no interaction while one flips, then the order of these events can be 
switched as well. However, we have to re-number the flip to conform to the
strand labelling convention.

The only subtle point is how flips interact with crossings. First of all,
if one ring flies through another while a third one flips, the order clearly
does not matter. If a ring flies through another and also flips, the 
order can be switched. However, if ring $A$ flips and then ring $B$ flies 
through it, this is homotopic to ring $B$ flying through ring $A$
from the other direction and then ring $A$ flipping. In other words, commuting
$\sigma_i$ with $w_i$ changes the ``sign of the crossing'' in the sense of
Exercise \ref{ex:swBn}. This gives the last, and the only truly non-commutative flip 
relation.

\parpic[r]{\begin{picture}(0,0)%
\includegraphics{figs/Wen.pstex}%
\end{picture}%
%
%  pstex_opts: -m 1.25 
%
\setlength{\unitlength}{4934sp}%
\begingroup\makeatletter\ifx\SetFigFont\undefined%
\gdef\SetFigFont#1#2#3#4#5{%
  \reset@font\fontsize{#1}{#2pt}%
  \fontfamily{#3}\fontseries{#4}\fontshape{#5}%
  \selectfont}%
\fi\endgroup%
\begin{picture}(545,1337)(-382,-464)
\end{picture}%
}
To explain why the flip is denoted by $w$, let us consider the alternative
description by ribbon tubes. A flipping ring traces a so called 
wen\footnote{The term wen was coined by Kanenobu and Shima in 
\cite{KanenobuShima:TwoFiltrationsR2K}}
in $\bbR^4$. A wen is a Klein bottle cut along a meridian circle, 
as shown. The wen is embedded in $\bbR^4$.

Finally, note that $\pi_1Y_n$ is exactly the pure $w$-braid group
$\PwB_n$: since each ring has to return to its original position and
orientation, each does an even number of flips.  The flips (or wens)
can all be moved to the bottoms of the braid diagram strands (to the
bottoms of the tubes, to the beginning of words), at a possible cost, as
specified by Equation~\eqref{eq:FlipRels}. Once together, they pairwise
cancel each other.  As a result, this group can be thought of as not
containing wens at all.

\subsubsection{The Relationship with u-Braids} \label{subsubsec:RelWithu} 
For
the sake of ignoring strand permutations, we restrict our
attention to pure braids.

\parpic[r]{$\xymatrix{
  \PuB \ar@{.>}[r]^{Z^u} \ar[d]^a & \calA^u \ar[d]^\alpha \\
  \PwB \ar[r]^{Z^w} & \calA^w
}$}
By Section \ref{subsubsec:FTAlgebraic}, for any expansion $Z^u\colon
\PuB_n \to \calA^u_n$ (where $\PuB_n$ is the ``usual'' braid group
and $\calA^u_n$ is the algebra of horizontal chord diagrams on $n$
strands), there is a square of maps as shown on the right. Here $Z^w$
is the expansion constructed in Section~\ref{subsec:wBraidExpansion},
the left vertical map $\glos{a}$ is the composition of the inclusion
and projection maps $\PuB_n \to \PvB_n \to \PwB_n$. The map $\glos{\alpha}$
is the induced map by the functoriality of projectivization, as noted
after Exercise \ref{ex:BraidsAlgApproach}.  The reader can verify that
$\alpha$ maps each chord to the sum of its two possible directed versions.

Note that this square is {\it not} commutative for any choice of $Z^u$ even
in degree 2: the image of a crossing under $Z^w$ is outside the image 
of $\alpha$.

\parpic[r]{\input{figs/uwsquare2.pstex_t}}
More specifically, for any choice $c$ of a ``parenthesization'' of $n$
points, the KZ-construction / Kontsevich integral (see for example
\cite{Bar-Natan:NAT}) returns an expansion $Z_c^u$ of $u$-braids. We shall 
see in \cite{Bar-NatanDancso:WKO2} (Proposition 4.15 there) 
that for any choice of $c$, the two
compositions $\alpha \circ Z_c^u$ and $Z^w \circ a$ are ``conjugate in a
bigger space'': there is a map $i$ from $\calA^w$ to a larger space of
``non-horizontal arrow diagrams'', and in this space the images of the
above composites are conjugate.  However, we are not certain that $i$
is an injection, and whether the conjugation leaves the $i$-image of
$\calA^w$ invariant, and so we do not know if the two compositions differ
merely by an outer automorphism of $\calA^w$.

\draftcut
\section{w-Knots} \label{sec:w-knots}

\begin{quote} \small {\bf Section Summary. }
  \summaryknots
\end{quote}

{\bf Knots are the wrong objects for study in knot theory,}
v-knots are the wrong objects for study in the theory of v-knotted
objects and w-knots are the wrong objects for study in the theory of
w-knotted objects. Studying uvw-knots on their own is the parallel of
studying cakes, cookies and pastries as they come out of the bakery --- we sure
want to make them our own, but the theory of desserts is more about the
ingredients and how they are put together than about the end products. In
algebraic knot theory this reflects through the fact that knots are
not finitely generated in any sense (hence, they must be made of some
more basic ingredients), and through the fact that there are very few
operations defined on knots (connected sums and satellite operations
being the main exceptions), and thus, most interesting properties of
knots are transcendental, or non-algebraic, when viewed from within the
algebra of knots and operations on knots~\cite{Bar-Natan:AKT-CFA}.

The right objects for study in knot theory, or v-knot theory or w-knot
theory, are thus the ingredients that make up knots and that permit a
richer algebraic structure. These are braids, studied in the previous
section, and even more so tangles and tangled graphs, studied in 
\cite{Bar-NatanDancso:WKO2}.  Yet tradition has its place and the sweets are
tempting, and we can introduce and apply some of the tools we will
use in the deeper and healthier study of w-tangles and w-tangled foams
in the limited, but tasty, arena of the baked goods of knot theory,
the knots themselves.

\draftcut \subsection{v-Knots and w-Knots} \label{subsec:VirtualKnots}
v-Knots may be understood either as knots drawn on surfaces modulo the
addition or removal of empty handles~\cite{Kuperberg:VirtualLink} 
or as ``Gauss diagrams'' (see Remark~\ref{rem:GD}),
or simply ``unembedded but wired together'' crossings modulo the
Reidemeister moves (\cite{Kauffman:VirtualKnotTheory, Roukema:GPV}
and Section 2 of \cite{Bar-NatanDancso:WKO2}). But right now we forgo the
topological and the abstract and give only the ``planar'' (and somewhat
less philosophically satisfying) definition of v-knots.

\begin{figure}[h]
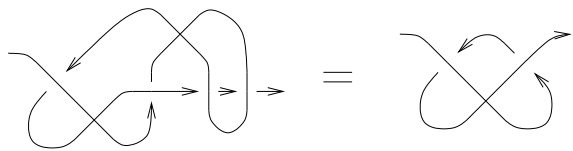

\[ \input figs/VKnot.pstex_t \]
\caption{
  A long v-knot diagram with 2 virtual crossings, 2 positive crossings and
  2 negative crossings. A positive-negative pair can easily be cancelled
  using R2, and then a virtual crossing can be cancelled using VR1, and it
  seems that the rest cannot be simplified any further.
} \label{fig:VKnot}
\end{figure}

\begin{figure}
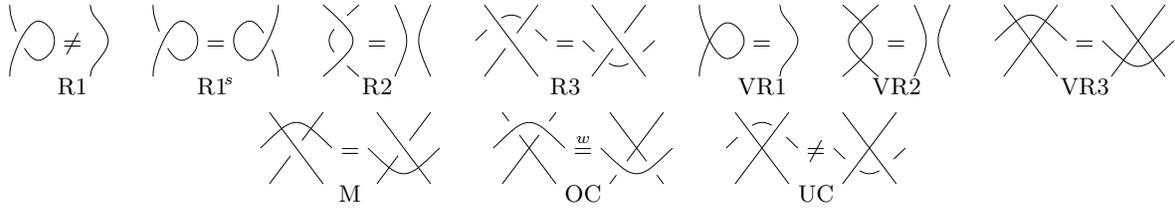

\[ \input figs/VKnotRels.pstex_t \]
\caption{
  The relations defining v-knots and w-knots, along with two relations that
  are {\em not} imposed.
} \label{fig:VKnotRels}
\end{figure}

\begin{definition} A ``long v-knot diagram'' is an arc smoothly
drawn in the plane from $-\infty$ to $+\infty$, with finitely many
self-intersections, divided into ``virtual crossings'' $\virtualcrossing$,
overcrossings $\overcrossing$ (a.k.a.~positive crossings), and undercrossings $\undercrossing$ (a.k.a.~negative crossings);
and regarded up to planar isotopy. A picture is worth more than a more
formal definition, and one appears in Figure~\ref{fig:VKnot}. A ``long
v-knot'' is an equivalence class of long v-knot diagrams, modulo the
equivalence generated by the Reidemeister $1^{\!s}$, 2 and 3 moves
(\glost{\Rs}, \glost{R2} and \glost{R3})\footnote{
\Rs\ is the ``spun'' version of R1 --- kinks can
be spun around, but not removed outright. See
Figure~\ref{fig:VKnotRels}.}, the virtual Reidemeister 1 through 3 moves
(\glost{VR1}, \glost{VR2}, \glost{VR3}), and by the mixed relations
(\glost{M}); all these are shown in Figure~\ref{fig:VKnotRels}. Finally,
``long w-knots'' are obtained from long v-knots by also dividing
by the overcrossings commute (OC) relations, also shown in
Figure~\ref{fig:VKnotRels}.  Note that we never mod out by the
Reidemeister 1 (\glost{R1}) move nor by the undercrossings commute
relation (UC).
\end{definition}

\begin{defwarn} A ``circular v-knot'' is like a long v-knot, except
parametrized by a circle rather than by a long line. Unlike the case of
usual knots, circular v-knots are {\bf not} equivalent to long v-knots \cite{Kauffman:VirtualKnotTheory}.
The same applies to w-knots.
\end{defwarn}

\begin{defwarn} Long v-knots form a monoid using the concatenation operation
$\#$. Unlike the case of usual knots, the resulting monoid is {\bf not} abelian \cite{Kauffman:VirtualKnotTheory}.
The same applies to w-knots.
\end{defwarn}

\begin{remark} \label{rem:GD} A
``Gauss diagram'' is a straight ``skeleton
line'' along with signed directed chords (signed ``arrows'') marked along
it (more at~\cite{Kauffman:VirtualKnotTheory,
GoussarovPolyakViro:VirtualKnots}). Gauss diagrams are in an obvious
bijection with long v-knot diagrams; the skeleton line of a Gauss diagram
corresponds to the parameter space of the v-knot, and the arrows
correspond to the crossings, with each arrow heading from the upper strand
to the lower strand, marked by the sign of the relevant crossing:
\[ \input figs/GDExample.pstex_t \]
One may also describe the relations in Figure~\ref{fig:VKnotRels} as well
as circular v-knots and other types of v-knots (as we will encounter later)
in terms of Gauss diagrams with varying skeletons.
\end{remark}

\begin{figure}
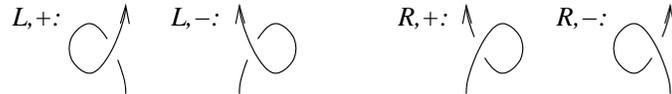

\[ \input figs/Kinks.pstex_t \]
\caption{
  The positive and negative under-then-over kinks (left), and the positive
  and negative over-then-under kinks (right). In each pair the
  negative kink is the $\#$-inverse of the positive kink, where $\#$
  denotes the concatenation operation.
\label{fig:Kinks}}
\end{figure}

\begin{remark}\label{rem:Framing} 
Since we do not mod out by R1, it is perhaps more
appropriate to call our class of v/w-knots ``framed long v/w-knots'',
but since we care more about framed v/w-knots than about unframed ones,
we reserve the unqualified name for the framed case, and when we do wish to
mod out by R1 we will explicitly write ``unframed long v/w-knots''.

Recall that in the case of ``usual knots'', or u-knots, dropping the R1
relation altogether also results in a $\bbZ^2$-extension of unframed
knot theory, where the two factors of $\bbZ$ are framing and rotation
number. If one wants to talk about ``true'' framed knots, one mods out
by the spun Reidemeister 1 relation (\Rs\ of Figure~\ref{fig:VKnotRels}),
which preserves the blackboard framing but does not preserve the rotation
number. We take the analogous approach here, including the \Rs\ relation
-- but not R1 -- in the v and w cases.

This said, note that the monoid of long v-knots is just a central extension
by $\bbZ$ of the monoid of unframed long v-knots, and so studying the
framed case is not very different from studying the unframed case. Indeed
the four ``kinks'' of Figure~\ref{fig:Kinks} generate a central $\bbZ$ within
long v-knots, and it is not hard to show that the sequence
\begin{equation} \label{eq:FramedAndUnframed}
   1\longrightarrow
   \bbZ \longrightarrow
   \{\text{long v-knots}\} \longrightarrow
   \{\text{unframed long v-knots}\} \longrightarrow 1
\end{equation}
is split and exact. The same can be said for w-knots.
\end{remark}

\begin{exercise} \label{ex:sl} Show that a splitting of the
sequence~\eqref{eq:FramedAndUnframed} is given by the ``self-linking''
invariant $\glos{\sl}\colon \{\text{long v-knots}\}\to\bbZ$ defined by
\[
  \sl(K):=\sum_{\text{crossings}\atop x\text{ in }K}\sign x ,
\]
where $K$ is a v-knot diagram, and the sign of a crossing $x$ is defined
so as to agree with the signs in Figure~\ref{fig:Kinks}.
\end{exercise}

\begin{remark} Note that  w-knots are strictly weaker than v-knots --- a notorious
example is the Kishino knot (e.g.~\cite{Dye:Kishinos}) which is non-trivial
as a v-knot; yet both it and its mirror are trivial as w-knots. Yet ordinary
knots inject even into w-knots, as the Wirtinger presentation makes sense
for w-knots and therefore w-knots have a ``fundamental quandle'' which
generalizes the fundamental quandle of ordinary
knots~\cite{Kauffman:VirtualKnotTheory}, and as the fundamental
quandle of ordinary knots separates ordinary
knots~\cite[Corollary~16.3]{Joyce:TheKnotQuandle}.
\end{remark}

\subsubsection{A topological construction of Satoh's tubing map}
\label{subsubsec:TopTube}
Following Satoh~\cite{Satoh:RibbonTorusKnots}
and using the same constructions as in Section~\ref{subsubsec:ribbon}, we
can map w-knots to (``long'') ribbon tubes in $\bbR^4$ (and the relations
in Figure~\ref{fig:VKnotRels} still hold). It is natural to expect that
this ``tubing'' map is an isomorphism; in other words, that the theory
of w-knots provides a ``Reidemeister framework'' for long ribbon tubes
in $\bbR^4$ --- that every long ribbon tube is in the image of this map
and that two ``w-knot diagrams'' represent the same long ribbon tube iff
they differ by a sequence of moves as in Figure~\ref{fig:VKnotRels}. This
remains unproven.

Let $\glos{\delta}\colon\{\text{v-knots}\} \to \{\text{Ribbon
tori in } \bbR^4\}$ denote the tubing map. In Satoh's~\cite{Satoh:RibbonTorusKnots}
$\delta$ is called ``Tube''.  It is worthwhile to give a completely
``topological'' definition of $\delta$. To do this we must start with
a topological interpretation of v-knots.

The standard topological interpretation of v-knots
(see e.g.~\cite{Kuperberg:VirtualLink}) is that they are oriented framed knots
drawn\footnote{Here and below, ``drawn on $\Sigma$'' means ``embedded in
$\Sigma\times[-\epsilon,\epsilon]$''.} on an oriented surface $\Sigma$,
modulo ``stabilization'', which is the addition and/or removal of empty
handles (handles that do not intersect with the knot). We prefer an
equivalent, yet even more bare-bones approach. For us, a virtual knot is an
oriented framed knot $\gamma$ drawn on a ``virtual surface $\glos{\Sigma}$ for
$\gamma$''. More precisely, $\Sigma$ is an oriented surface that may have
a boundary, $\gamma$ is drawn on $\Sigma$, and the pair $(\Sigma,\gamma)$
is taken modulo the following relations:
\begin{itemize}
\item Isotopies of $\gamma$ on $\Sigma$ (meaning, in
  $\Sigma\times[-\epsilon,\epsilon]$).
\item Tearing and puncturing parts of $\Sigma$ away from $\gamma$:
\end{itemize}
\[ \input{figs/TearingAndPuncturing.pstex_t} \]
(We call $\Sigma$ a ``virtual surface'' because tearing and puncturing
imply that we only care about it in the immediate vicinity of $\gamma$).

We can now define\footnote{Following a private discussion with Dylan
Thurston.} a map $\delta$, defined on v-knots and taking values
in ribbon tori in $\bbR^4$: given $(\Sigma,\gamma)$, embed $\Sigma$
arbitrarily in $\bbR^3_{xzt}\subset\bbR^4$. Note that the unit normal
bundle of $\Sigma$ in $\bbR^4$ is a trivial circle bundle and it has a
distinguished trivialization, constructed using its positive-$y$-direction
section and the orientation that gives each fibre a linking number
$+1$ with the base $\Sigma$.  We say that a normal vector to $\Sigma$
in $\bbR^4$ is ``near unit'' if its norm is between $1-\epsilon$ and
$1+\epsilon$. The near-unit normal bundle of $\Sigma$ has as fibre
an annulus that can be identified with $[-\epsilon,\epsilon]\times
S^1$ (identifying the radial direction $[1-\epsilon,1+\epsilon]$
with $[-\epsilon,\epsilon]$ in an orientation-preserving manner), and
hence, the near-unit normal bundle of $\Sigma$ defines an embedding
of $\Sigma\times[-\epsilon,\epsilon]\times S^1$ into $\bbR^4$. On the
other hand, $\gamma$ is embedded in $\Sigma\times[-\epsilon,\epsilon]$ so
$\gamma\times S^1$ is embedded in $\Sigma\times[-\epsilon,\epsilon]\times
S^1$, and we can let $\delta(\gamma)$ be the composition
\[ \gamma\times S^1
  \hookrightarrow\Sigma\times[-\epsilon,\epsilon]\times S^1
  \hookrightarrow\bbR^4,
\]
which is a torus in $\bbR^4$, oriented using the given orientation of
$\gamma$ and the standard orientation of $S^1$.

\begin{figure}
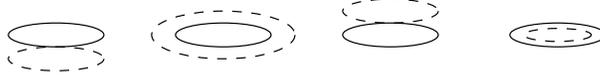

 \input figs/CompanionRing.pstex_t
\caption{Framing as companion rings.}\label{fig:CompanionRing}
\end{figure}

A framing of a knot (or a v-knot) $\gamma$ can be thought of as a
``nearby companion'' to $\gamma$. Applying the above procedure to a knot
and a nearby companion simultaneously, we find that $\delta$ takes framed
v-knots to framed ribbon tori in $\bbR^4$, where a framing of a tube in
$\bbR^4$ is a continuous up-to-homotopy choice of unit normal vector at
every point of the tube. Note that from the perspective of flying rings as
in Section~\ref{subsubsec:FlyingRings}, a framing is a ``companion ring''
to a flying ring. In the framing of $\delta(\gamma)$ the companion ring
is never linked with the main ring, but can fly parallel inside, outside,
above or below it and change these positions, as shown in Figure~\ref{fig:CompanionRing}.

We leave it to the reader to verify that $\delta(\gamma)$ is ribbon, that
it is independent of the choices made within its construction, that it is
invariant under isotopies of $\gamma$ and under tearing and puncturing of
$\Sigma$, that it is also invariant under the OC
relation of Figure~\ref{fig:VKnotRels} and hence, the true domain of
$\delta$ is w-knots, and that it is equivalent to Satoh's tubing map.

\draftcut
\subsection{Finite Type Invariants of v-Knots and w-Knots}
\label{subsec:FTforvwKnots}

Much as for v-braids and w-braids (see Section~\ref{subsec:FT4Braids}) and
much as for ordinary knots (e.g.~\cite{Bar-Natan:OnVassiliev}) we define
finite type invariants for v-knots and for w-knots using an alternation
scheme with $\semivirtualover\to\overcrossing-\virtualcrossing$
and $\semivirtualunder\to\virtualcrossing-\undercrossing$. That is,
given any invariant of v- or w-knots taking values in an abelian group, 
we extend the invariant to v- or
w-knots also containing ``semi-virtual crossings'' like $\semivirtualover$
and $\semivirtualunder$ using the above assignments, and we declare an
invariant to be ``of type $m$'' if it vanishes on v- or w-knots with more
than $m$ semi-virtuals. As for v- and w-braids and as for ordinary knots,
such invariants have an ``$m$th derivative'', their ``weight system'',
which is a linear functional on the space $\calA^{sv}(\uparrow)$ (for
v-knots) or $\calA^{sw}(\uparrow)$ (for w-knots). We turn to the definitions
of these spaces, following~\cite{GoussarovPolyakViro:VirtualKnots,
Bar-NatanHalachevaLeungRoukema:v-Dims}:

\begin{definition} \label{def:ArrowDiagrams} An ``arrow diagram''
is a chord diagram along a long line (called ``the skeleton''),
in which the chords are oriented (hence ``arrows''). An example is
given in Figure~\ref{fig:ADand6T}. Let $\glos{\calD^v(\uparrow)}$ be the space of
formal linear combinations of arrow diagrams.  Let $\glos{\calA^v(\uparrow)}$
be $\calD^v(\uparrow)$ modulo all ``6T relations''. Here a 6T relation is
any (signed) combination of arrow diagrams obtained from the diagrams in
Figure~\ref{fig:6T} by placing the 3 vertical strands there along a long
line skeleton in any order, and possibly adding some further arrows in between, 
as shown in Figure~\ref{fig:ADand6T}. Let $\glos{\calA^{sv}(\uparrow)}$
be the further quotient of $\calA^v(\uparrow)$ by the \glost{RI} relation,
where the RI (for rotation number independence) relation asserts that an
isolated arrow pointing to the right equals an isolated arrow pointing
to the left\footnote{
  The XII relation of~\cite{Bar-NatanHalachevaLeungRoukema:v-Dims} follows
  from RI and need not be imposed.
}, as shown in Figure~\ref{fig:ADand6T}.

\begin{figure}
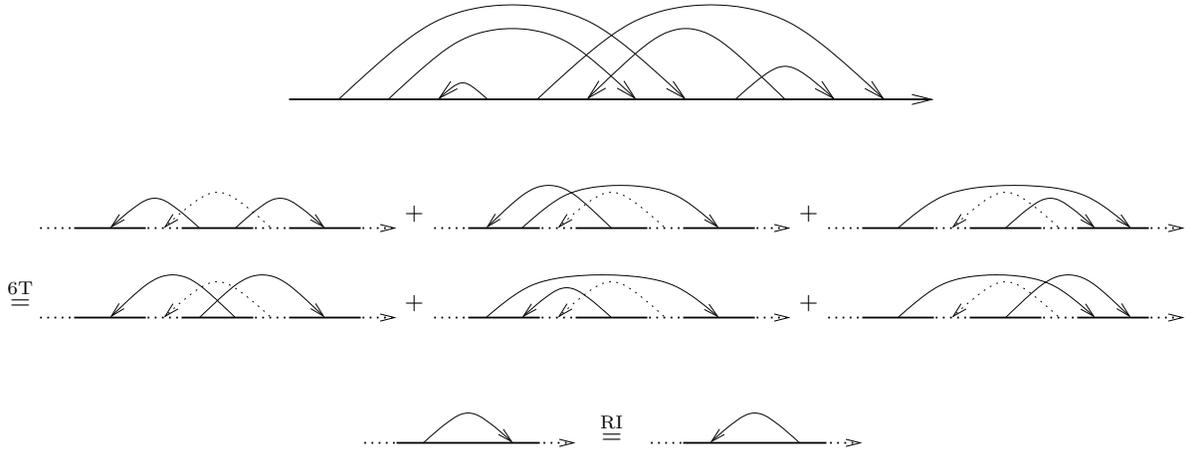

\[ \input figs/ADand6T.pstex_t \]
\caption{
  An arrow diagram of degree 6, a 6T relation, and an RI relation. The dotted parts indicate 
  that there may be more arrows on other parts of the skeleton, however these remain the 
  same throughout the relation.
} \label{fig:ADand6T}
\end{figure}

Let $\glos{\calA^w(\uparrow)}$ be the further quotient of $\calA^v(\uparrow)$
by the TC relation, first displayed
in Figure~\ref{fig:TCand4T} and reproduced for the case of a
long line skeleton in Figure~\ref{fig:TCand4TForKnots}. Likewise, let
$\glos{\calA^{sw}(\uparrow)}:=\calA^{sv}(\uparrow)/TC=\calA^w(\uparrow)/RI$.
Alternatively, noting that given TC two of the terms in 6T drop out,
$\calA^w(\uparrow)$ is the space of formal linear combinations
of arrow diagrams modulo TC and $\aft$ relations, displayed in
Figures~\ref{fig:TCand4T} and~\ref{fig:TCand4TForKnots}. Likewise,
$\calA^{sw}=\calD^v/TC,\aft,RI$. Finally, grade $\calD^v(\uparrow)$ and
all of its quotients by declaring that the degree of an arrow diagram
is the number of arrows in it.

\begin{figure}
\[ \input figs/TCand4TForKnots.pstex_t \]
\caption{The TC and the $\protect\aft$ relations for knots.}
\label{fig:TCand4TForKnots}
\end{figure}

\end{definition}

As an example, the spaces $\calA^{v,sv,w,sw}(\uparrow)$ (that is, any of the spaces above)
restricted to degrees up to 2 are studied in detail in
Section~\ref{subsec:ToTwo}.

In the same manner as in the theory of finite type invariants of ordinary
knots (see especially~\cite[Section~3]{Bar-Natan:OnVassiliev}), the spaces
$\calA^{\star}(\uparrow)$ (meaning, all of the spaces above) carry much algebraic structure.  The
juxtaposition product makes them into graded algebras. The product of two
finite type invariants is a finite type invariant (whose type is the sum
of the types of the factors); this induces a product on weight systems,
and therefore a co-product $\Delta$ on arrow diagrams. In brief (and much
the same as in the usual finite type story), the co-product $\Delta D$
of an arrow diagram $D$ is the sum of all ways of dividing the arrows
in $D$ between a ``left co-factor'' and a ``right co-factor''. In summary:

\begin{proposition} \label{prop:CoarseStructure} $\calA^v(\uparrow)$,
$\calA^{sv}(\uparrow)$, $\calA^w(\uparrow)$, and
$\calA^{sw}(\uparrow)$ are co-commutative graded bi-algebras.
\end{proposition}

By the Milnor-Moore theorem~\cite[Theorem 6.11]{MilnorMoore:Hopf} we find that
$\calA^{v,sv,w,sw}(\uparrow)$ are the
universal enveloping algebras of their Lie algebras of primitive
elements (that is, elements $D$ such that $\Delta(D)=1\otimes D + D \otimes 1$). Denote these (graded) Lie algebras by
$\glos{\calP^{v,sv,w,sw}(\uparrow)}$, respectively.

When we grow up we'd like to understand $\calA^v(\uparrow)$ and
$\calA^{sv}(\uparrow)$. At the moment we know only very little about these
spaces beyond the generalities of Proposition~\ref{prop:CoarseStructure}.
In Section \ref{subsec:SomeDimensions} some dimensions of low degree parts of
$\calA^{v,sv}(\uparrow)$ are discussed. Also, given a finite dimensional
Lie bialgebra and a finite dimensional representation thereof, we know
how to construct linear functionals on $\calA^v(\uparrow)$ (one in each
degree~\cite{Haviv:DiagrammaticAnalogue, Leung:CombinatorialFormulas}),
but not on $\calA^{sv}(\uparrow)$. But we don't even know which degree
$m$ linear functionals on $\calA^{sv}(\uparrow)$ are the weight systems
of degree $m$ invariants of v-knots (that is, we have not solved the
``Fundamental Problem''~\cite{Bar-NatanStoimenow:Fundamental} for
v-knots).

As we shall see below, the situation is much brighter for
$\calA^{w,sw}(\uparrow)$.

\draftcut
\subsection{Expansions for w-Knots} \label{subsec:Z4Knots}
The notion of ``an expansion'' (or ``a universal finite type invariant'')
for w-knots (or v-knots) is defined in complete analogy with the
parallel notion for usual knots (see e.g.~\cite{Bar-Natan:OnVassiliev}),
except replacing double points $\doublepoint$ with semi-virtual
crossings $\semivirtualover$ and $\semivirtualunder$, and replacing
chord diagrams by arrow diagrams. Alternatively, it is the same as
an expansion for w-braids (as in Definition~\ref{def:vwbraidexpansion}),
simply replacing w-braids by w-knots. Just as in the
cases of u-knots (i.e., ordinary knots) and/or w-braids, the existence of an expansion
$Z\colon \{\text{w-knots}\}\to\calA^{sw}(\uparrow)$ is equivalent to the
statement ``every weight system integrates'', i.e., ``every degree $m$
linear functional on $\calA^{sw}(\uparrow)$ is the $m$th derivative of
a type $m$ invariant of long w-knots''.

\begin{theorem} \label{thm:ExpansionForKnots}
There exists an expansion $Z\colon \{\text{w-knots}\}\to\calA^{sw}(\uparrow)$.
\end{theorem}

\begin{proof} It is best to define $Z$ by an example, and it is best to
display the example only as a picture:
\[ \input figs/ZwKnotsExample.pstex_t \]
It is clear how to define $Z(K)$ in the general case --- for every crossing
in $K$ place an exponential reservoir of arrows (compare
with Equation~\eqref{eq:reservoir}) next to that crossing, with
the arrows heading from the upper strand to the lower strand, taking
positive reservoirs ($e^a$, with $a$ symbolizing the arrow) for positive
crossings and negative reservoirs ($e^{-a}$) for negative crossings, and
then tug the skeleton until it looks like a straight line. Note that the
TC relation in $\calA^{sw}$ is used to show that all reasonable
ways of placing an arrow reservoir at a crossing (with its heading and sign
fixed) are equivalent:
\[ \input figs/FourWays.pstex_t \]

The same proof that shows the invariance of $Z$ in the braid case
(see Theorem~\ref{thm:RInvariance}) works here as well\footnote{A tiny bit of
extra care is required for invariance under \Rs: it easily follows from
RI.}, and the same argument as in the braid case shows the universality
of $Z$. \qed
\end{proof}

\begin{remark} \label{rem:ZwForGD} Using the language of Gauss diagrams
(Remark~\ref{rem:GD}) the definition of $Z$ is even simpler. Simply map
every positive arrow in a Gauss diagram to a positive ($e^a$) reservoir,
and every negative one to a negative ($e^{-a}$) reservoir:
\[ \input figs/ZwForGD.pstex_t \]
\end{remark}

An expansion (a universal finite type invariant) is as interesting as its
target space, for it is just a tool that takes linear functionals on the
target space to finite type invariants on its domain space. The purpose
of the next section is to find out how interesting are our present target
space, $\calA^{sw}(\uparrow)$, and its ``parent'', $\calA^w(\uparrow)$.

\draftcut
\subsection{Jacobi Diagrams, Trees and Wheels} \label{subsec:Jacobi}

In studying $\calA^w(\uparrow)$ we again follow the model set
by usual knots: we introduce the space $\calA^{wt}$ of ``w-Jacobi diagrams'' and show that
it is isomorphic to $\calA^w$. Major advantages of working with $\calA^{wt}$
are that the co-product, the primitives, and the relationship with Lie algebras are
much more natural and easy to describe. Compare the following definitions and theorem
with~\cite[Section~3]{Bar-Natan:OnVassiliev}.

\begin{figure}
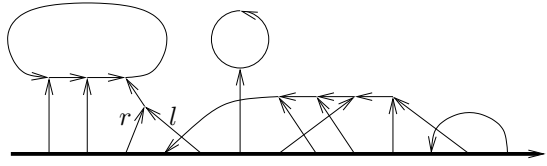

\[ \input figs/wJacDiag.pstex_t \]
\caption{A degree 11 w-Jacobi diagram on a long line skeleton. It has a
skeleton line at the bottom, 13 vertices along the skeleton (of which 2 are
incoming and 11 are outgoing), 9 internal vertices (with only one
explicitly marked with ``left'' ($l$) and ``right'' ($r$)) and one
bubble. The five quadrivalent vertices that seem to appear in the diagram
are just projection artifacts and graph-theoretically, they don't exist.}
\label{fig:wJacDiag}
\end{figure}

\begin{definition} \label{def:wJac} A ``w-Jacobi diagram on a long line
skeleton''\footnote{What a mouthful! We usually short this to
``w-Jacobi diagram'', or sometimes ``arrow diagram'' or just ``diagram''.}
is a connected graph made of the following ingredients:
\begin{itemize}
\item A ``long'' oriented ``skeleton'' line. We usually draw the skeleton
  line a bit thicker for emphasis.
\item Other directed edges, usually called ``arrows''.
\item Trivalent ``skeleton vertices'' in which an arrow starts or ends on
  the skeleton line.
\item Trivalent ``internal vertices'' in which {\em two arrows end and one arrow
  begins} (this will be important in Section \ref{subsec:LieAlgebras} where
  we relate these diagrams to Lie algebras). 
  The internal vertices are ``oriented'' --- of the two arrows that
  end in an internal vertices, one is marked as ``left'' and the other is
  marked as ``right''. In reality when a diagram is drawn in the plane, we
  almost never mark ``left'' and ``right'', but instead assume the
  ``left'' and ``right'' inherited from the plane, as seen from the
  outgoing arrow from the given vertex.
\end{itemize}
Note that we allow multiple arrows connecting the same two vertices
(though at most two are possible, given connectedness and trivalence)
and we allow ``bubbles'' --- arrows that begin and end in the same
vertex. Also keep in mind that for the purpose of determining equality of diagrams
the skeleton line is distinguished.
The ``degree'' of a w-Jacobi diagram is half the number of
trivalent vertices in it, including both internal and skeleton vertices.
An example of a w-Jacobi diagram is in Figure~\ref{fig:wJacDiag}.
\end{definition}

\begin{figure}
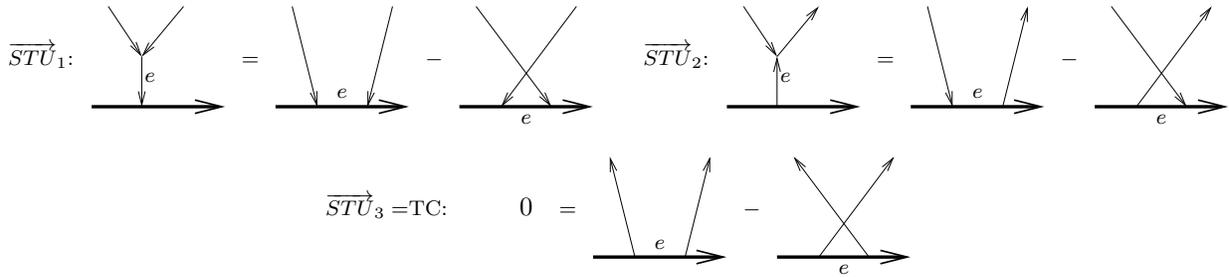

\[ \input figs/aSTU.pstex_t \]
\caption{The $\protect\aSTU_{1,2}$ and TC relations with
their ``central edges'' marked $e$.}
\label{fig:aSTU}
\end{figure}

\begin{definition}
Let $\glos{\calD^{wt}}(\uparrow)$ be the graded vector space of formal
linear combinations\footnote{$\bbQ$-linear, or any other field of 
characteristic 0.} of w-Jacobi diagrams on a long line skeleton,
and let $\glos{\calA^{wt}}(\uparrow)$ be $\calD^{wt}(\uparrow)$
modulo the $\glos{\aSTU_1}$, $\glos{\aSTU_2}$, and TC relations of
Figure~\ref{fig:aSTU}. Note that that each diagram appearing in each
$\aSTU$ relation has a ``central edge'' $e$ which can serve as an
``identifying name'' for that $\aSTU$. Thus, given a diagram $D$ with
a marked edge $e$ which is either on the skeleton or which contacts
the skeleton, there is an unambiguous $\aSTU$ relation ``around'' or
``along'' the edge $e$.
\end{definition}

\begin{figure}
\[ \input figs/aIHX.pstex_t \]
\caption{The $\protect\aAS$ and $\protect\aIHX$ relations.}
\label{fig:aIHX}
\end{figure}

We like to call the following theorem ``the bracket-rise theorem'',
for it justifies the introduction of internal vertices, and as
should be clear from the $\aSTU$ relations and as will become even
clearer in Section~\ref{subsec:LieAlgebras}, internal vertices can be
viewed as ``brackets''. Two other bracket-rise theorems are Theorem~6
of~\cite{Bar-Natan:OnVassiliev} and Ohtsuki's theorem, i.e., Theorem~4.9
of~\cite{Polyak:ArrowDiagrams}.

\begin{theorem}[Bracket-rise] \label{thm:BracketRise} The obvious inclusion
$\iota\colon \calD^v(\uparrow)\to\calD^{wt}(\uparrow)$ of arrow diagrams
(see Definition~\ref{def:ArrowDiagrams}) into w-Jacobi diagrams descends
to the quotient $\calA^w(\uparrow)$ and induces an isomorphism\footnote{At 
this point a vector space isomorphism, but we'll soon define a bi-algebra 
structure on $\calA^{wt}$ to make it into an isomorphism of bi-algebras.}
$\bar\iota\colon \calA^w(\uparrow)\stackrel{\sim}{\longrightarrow}
\calA^{wt}(\uparrow)$.  Furthermore, the $\glos{\aAS}$ and $\glos{\aIHX}$
relations of Figure~\ref{fig:aIHX} hold in $\calA^{wt}(\uparrow)$.
\end{theorem}

\begin{proof} The proof, joint with D.~Thurston, is modelled after
the proof of Theorem~6 of~\cite{Bar-Natan:OnVassiliev}. To show that
$\iota$ descends to $\calA^w(\uparrow)$ we just need to show that in
$\calA^{wt}(\uparrow)$, $\aft$ follows from $\aSTU_{1,2}$. Indeed,
applying $\aSTU_1$ along the edge $e_1$ and $\aSTU_2$ along $e_2$ in
the picture below, we get the two sides of $\aft$:
\begin{equation} \label{eq:STUto4T}
  \input figs/STUto4T.pstex_t
\end{equation}

The fact that $\bar\iota$ is surjective is easy: indeed, for diagrams
in $\calA^{wt}(\uparrow)$ that have no internal vertices there is
nothing to show, for they are really in $\calA^w(\uparrow)$. Further,
by repeated use of $\aSTU_{1,2}$ relations, all internal vertices in
any diagram in $\calA^{wt}(\uparrow)$ can be removed (remember that
the diagrams in $\calA^{wt}(\uparrow)$ are always connected, and in
particular, if they have an internal vertex they must have an internal
vertex connected by an edge to the long line skeleton, and the latter 
vertex can be removed first).

To complete the proof that $\bar\iota$ is an isomorphism it is enough
to show that the ``elimination of internal vertices'' procedure of
the last paragraph is well-defined -- that its output is independent
of the order in which $\aSTU_{1,2}$ relations are applied in order to
eliminate internal vertices. Indeed, this done, the elimination map would
by definition satisfy the $\aSTU_{1,2}$ relations and thus, descend to
a well-defined inverse for $\bar\iota$.

On diagrams with just one internal vertex, Equation~\eqref{eq:STUto4T}
shows that all ways of eliminating that vertex are equivalent modulo $\aft$
relations, and hence, the elimination map is well-defined on such diagrams.

We proceed by induction on the number of internal vertices. We have shown 
that $\bar\iota$ is well-defined if there is only one internal vertex.
Now assume that we have shown that the elimination map is well defined on
all diagrams with at most $k$ internal vertices for some $k\geq 1$ positive integer, 
and let $D$ be a diagram
with $(k+1)$ internal vertices. Let
$e$ and $e'$ be edges in $D$ that connect the skeleton of $D$ to an
internal vertex. We need to show that any elimination process that
begins with eliminating $e$ yields the same answer, modulo $\aft$, as
any elimination process that begins with eliminating $e'$. There are
several cases to consider.

\parpic[r]{$\input figs/CaseI.pstex_t$}
{\bf Case I.} $e$ and $e'$ connect the skeleton to {\em different} internal
vertices of $D$. In this case, after eliminating $e$ we get a signed sum
of two diagrams with exactly 7 internal vertices, and since the elimination
process is well-defined on such diagrams, we may as well continue by
eliminating $e'$ in each of those, getting a signed sum of 4 diagrams with
6 internal vertices each. On the other hand, if we start by eliminating
$e'$ we can continue by eliminating $e$, and we get the {\em same} signed
sum of 4 diagrams with 6 internal vertices.

\parpic[r]{$\input figs/CaseII.pstex_t$}
{\bf Case II.} $e$ and $e'$ are connected to the same internal vertex $v$
of $D$, yet some other edge $e''$ exists in $D$ that connects the skeleton
of $D$ to some other internal vertex $v'$ in $D$. In that case, use the
previous case and the transitivity of equality: (elimination starting with
$e$)=(elimination starting with $e''$)=(elimination starting with $e'$).

\parpic[r]{$\input figs/CaseIII.pstex_t$}
{\bf Case III.} This is what remains if neither Case I nor Case II
hold. In that case, $D$ must have a schematic form as on the right,
with the ``blob'' not connected to the skeleton other than via $e$
or $e'$, yet further arrows may exist outside of the blob. Let $f$
denote the edge connecting the blob to $e$ and $e'$. The ``two in one
out'' rule for vertices implies that any part of a diagram must have
an excess of incoming edges over outgoing edges, equal to the total
number of vertices in that diagram part. Applying this principle to
the blob, we find that it must contain exactly one vertex, as shown on
the right below. Then by the ``two in one out'' rule $f$ must be oriented 
upwards, and hence, by the ``two in one out'' rule again, 
$e$ and $e'$ must be oriented upwards as well.

\parpic[r]{$\input figs/CaseIIIa.pstex_t$}
We leave it to the reader to verify that in this case the two ways of
applying the elimination procedure, $e$ and then $f$ or $e'$ and then $f$,
yield the same answer modulo $\aft$ (in fact, that answer is $0$).

We also leave it to the reader to verify that $\aSTU_1$ implies $\aAS$
and $\aIHX$.  In Section \ref{subsec:LieAlgebras} we'll describe
the relationship between $\calA^{wt}$ and Lie algebras. 
Algebraically, the relations $\aSTU_1$, $\aAS$
and $\aIHX$ are restatements of the anti-symmetry
of the bracket and of Jacobi's identity: if $[x,y]:=xy-yx$, then
$0=[x,y]+[y,x]$ and $[x,[y,z]]=[[x,y],z]-[[x,z],y]$. \qed
\end{proof}

Note that $\calA^{wt}(\uparrow)$ inherits algebraic structure from
$\calA^w(\uparrow)$: it is an algebra by concatenation of diagrams,
and a co-algebra with $\Delta(D)$, for $D\in\calD^{wt}(\uparrow)$,
being the sum of all ways of dividing $D$ between a ``left co-factor''
and a ``right co-factor'' so that connected components of $D-S$
are kept intact, where $S$ is the skeleton line of $D$ (compare
with~\cite[Definition~3.7]{Bar-Natan:OnVassiliev}).

\begin{figure}
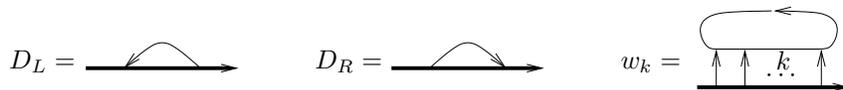

\[ \input figs/AwGenerators.pstex_t \]
\caption{The left-arrow diagram $D_L$, the right-arrow diagram $D_R$ and
  the $k$-wheel $w_k$.}
\label{fig:AwGenerators}
\end{figure}

As $\calA^w(\uparrow)$ and $\calA^{wt}(\uparrow)$ are canonically
isomorphic, from this point on we will not keep the distinction between
the two spaces. 
One may add the RI relation to the definition of $\calA^{wt}(\uparrow)$
to get a space $\calA^{swt}(\uparrow)$. For an unframed version one may 
add the stronger framing independence (FI) relation, setting $D_L=D_R=0$, with $D_L$
and $D_R$ the single arrows as in Figure~\ref{fig:AwGenerators}. The resulting space is 
called $\calA^{rwt}(\uparrow)$. The statement and proof of the
bracket rise theorem adapt with no difficulty, and we find that
$\calA^{sw}(\uparrow)\cong\calA^{swt}(\uparrow)$ and
$\calA^{rw}(\uparrow)\cong\calA^{rwt}(\uparrow)$. In the
future we'll drop the $t$ from all superscripts.

The advantages of allowing internal trivalent vertices are
already apparent (for example, note that there is a nice description of
primitive elements: they are the arrow diagrams which remain connected 
if the skeleton is removed). Further advantages will emerge in Section
\ref{subsec:LieAlgebras}.

\begin{theorem} \label{thm:Aw}
The bi-algebra $\calA^w(\uparrow)$ is the bi-algebra of polynomials
in the diagrams $\glos{D_L}$, $\glos{D_R}$ and $\glos{w_k}$
(for $k\geq 1$) shown in Figure~\ref{fig:AwGenerators}, where
$\deg D_L=\deg D_R=1$ and $\deg w_k=k$, subject to the one relation
$w_1=D_L-D_R$. Thus, $\calA^w(\uparrow)$ has two generators in degree
1 and one generator in every degree greater than 1, as stated in
Section~\ref{subsec:SomeDimensions}.
\end{theorem}

\begin{proof} (sketch). Readers familiar with the diagrammatic PBW
theorem~\cite[Theorem~8]{Bar-Natan:OnVassiliev} will note that it has
a direct analogue for the $\calA^w(\uparrow)$ case, and that the proof
in~\cite{Bar-Natan:OnVassiliev} carries through almost verbatim. Namely,
the space $\calA^w(\uparrow)$ is isomorphic to a space $\glos{\calB^w}$
of ``unitrivalent diagrams'' with symmetrized univalent ends modulo
$\aAS$ and $\aIHX$. Given the ``two in one out'' rule for arrow
diagrams in $\calA^w(\uparrow)$ (and hence, in $\calB^w$)
the connected components of diagrams in $\calB^w$ can only be
``trees'' or ``wheels''. A tree is a unitrivalent diagram with no
cycles (oriented or not). A wheel is a single oriented cycle with some
number of incoming ``spokes'' (see $w_k$ in Figure \ref{fig:AwGenerators} and remove the skeleton line).
The reader might object that there are also ``wheels of trees'': trees 
attached to an oriented cycle, but these can be reduced to linear
combinations of wheels using the $\aIHX$ relation.

Trees vanish if they have more than one leaf, as their
leafs are symmetric while their internal vertices are anti-symmetric,
so $\calB^w$ is generated by wheels and by the one-leaf-one-root tree, which is simply
a single arrow. Wheels map to the $w_k$'s in
$\calA^w(\uparrow)$ under the isomorphism, and  the arrow maps to the average of 
$D_L$ and $D_R$. The relation $w_1=D_L-D_R$ is then easily verified using $\aSTU_2$.

One may also argue directly, without using $\calB^w$. In
short, let $D$ be a diagram in $\calA^w(\uparrow)$ and $S$ is its
skeleton. Then $D-S$ may have several connected components, whose ``legs''
are intermingled along $S$. Using the $\aSTU$ relations these legs can
be sorted (at a cost of diagrams with fewer connected components, which
could have been treated earlier in an inductive proof). At the end of the
sorting procedure one can see that the only diagrams that remain are our
declared generators. It remains to show that our generators are linearly
independent (apart for the relation $w_1=D_L-D_R$). For the generators
in degree 1, simply write everything out explicitly in the spirit of
Section~\ref{subsubsec:DegreeOne}. In higher degrees there is only one
primitive diagram in each degree, so it is enough to show that $w_k\neq
0$ for every $k$. This can be done ``by hand'', but it is more easily
done using Lie algebraic tools in Section~\ref{subsec:LieAlgebras}. \qed
\end{proof}

\begin{exercise} \label{exe:Asw} Show that the bi-algebra
$\calA^{rw}(\uparrow)$ (see Section~\ref{subsec:SomeDimensions}) is
the bi-algebra of polynomials in the wheel diagrams $w_k$ ($k\geq 2$),
and that $\calA^{sw}(\uparrow)$ is the bi-algebra of polynomials in the
same wheel diagrams and an additional generator $\glos{D_A}:=D_L=D_R$.
\end{exercise}

\begin{proposition} \label{prop:AwCirc} In $\calA^w(\bigcirc)$ all wheels
vanish, and hence, the bi-algebra $\calA^w(\bigcirc)$ is the bi-algebra
of polynomials in a single variable $D_L=D_R$.
\end{proposition}

\begin{figure}
 \input{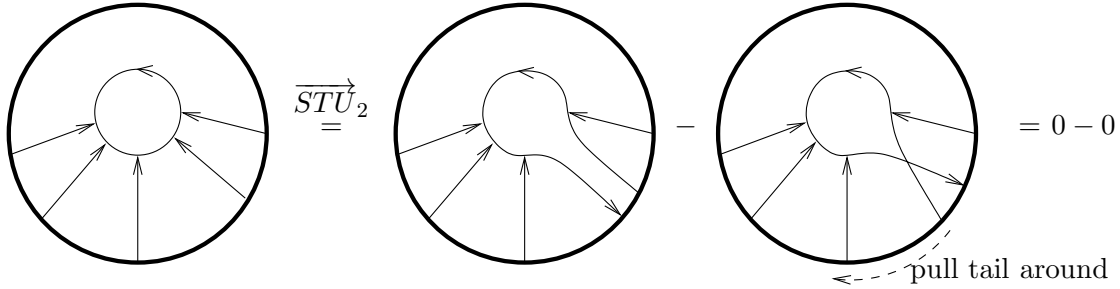}
\caption{Wheels in a circle vanish.}\label{fig:WheelInCircle}
 \end{figure}

\begin{proof} This is Lemma~2.7 of~\cite{Naot:BF}. In short, a wheel in
$\calA^w(\bigcirc)$ can be reduced using $\aSTU_2$ to a difference of
trees, as shown in Figure \ref{fig:WheelInCircle}. 
One of these trees has two adjoining leafs, and hence, it is 0 by TC and
$\aAS$. In the other two of the leafs can be commuted ``around the circle''
using TC until they are adjoining and hence vanish by TC and $\aAS$.
\qed
\end{proof}

\begin{exercise} Show that $\calA^{sw}(\bigcirc)\cong\calA^w(\bigcirc)$
yet $\calA^{rw}(\bigcirc)$ vanishes except in degree $0$.
\end{exercise}

The following two exercises may help the reader to develop a better
``feel'' for $\calA^w(\uparrow)$ and will be needed, within the discussion
of the Alexander polynomial (especially within
Definition~\ref{def:InterpretationMap}).

\parpic[r]{\raisebox{-12mm}{$\input figs/CC.pstex_t$}}
\begin{exercise} Show
that the ``commutators commute'' (\glost{CC}) relation,
shown on the right, holds in $\calA^w(\uparrow)$. (Interpreted in
Lie algebras as in the next section, this relation becomes $[[x,y],
[z,w]]=0$, and hence the name ``commutators commute''). Note that the
proof of CC depends on the skeleton having a single component; later,
when we will work with $\calA^w$-spaces with more complicated skeleta,
the CC relation will not hold.
\end{exercise}

\parpic[r]{\raisebox{-2mm}{$\input figs/Hair.pstex_t$}}
\begin{exercise} \label{ex:Hair} Show that ``detached wheels'' and
``hairy $Y$'s'' make sense in $\calA^w(\uparrow)$. As on the right, a
detached wheel is a wheel with a number of spokes, and a hairy $Y$ is a
combinatorial $Y$ shape (three arrows meeting at a single internal vertex) 
with further ``hair'' on its trunk (its outgoing
arrow). It is specified where the trunk and the leafs of the $Y$ connect
to the skeleton, but it is not specified where the spokes of the wheel
and where the hair on the $Y$ connect to the skeleton. The content of the
exercise is to show that modulo the relations of $\calA^w(\uparrow)$,
it is not necessary to specify this further information: all ways of
connecting the spokes and the hair to the skeleton are equivalent. Like
the previous exercise, this result depends on the skeleton having a
single component.
\end{exercise}

\begin{remark} In the case of usual knots and usual chord diagrams,
Jacobi diagrams have a topological interpretation using the
Goussarov-Habiro calculus of claspers~\cite{Goussarov:3Manifolds,
Habiro:Claspers}. In the w case a similar such calculus was developed by 
Watanabe in~\cite{Watanabe:ClasperMoves}. Various related results are 
at~\cite{HabiroKanenobuShima:R2K, HabiroShima:R2KII}.
\end{remark}

\draftcut
\subsection{The Relation with Lie Algebras} \label{subsec:LieAlgebras}
The
theory of finite type invariants of knots is related to the theory
of metrized Lie algebras via the space $\calA$ of chord diagrams, as
explained in~\cite[Theorem~4 and Exercise~5.1]{Bar-Natan:OnVassiliev}. In
a similar manner the theory of finite type invariants of w-knots is
related to arbitrary finite-dimensional Lie algebras (or equivalently, to
doubles of co-commutative Lie bialgebras, as explained below) via the space $\calA^w(\uparrow)$
of arrow diagrams.

\subsubsection{Preliminaries} Given a finite dimensional Lie 
algebra\footnote{Over $\bbQ$, or another field of characteristic 0.}
$\glos{\frakg}$ let $\glos{I\frakg}:=\frakg^\ast\rtimes\frakg$ be the
semi-direct product of the dual $\frakg^\ast$ of $\frakg$ with $\frakg$,
with $\frakg^\ast$ taken as an abelian algebra and with $\frakg$ acting
on $\frakg^\ast$ by the usual coadjoint action. In formulae,
\[ I\frakg=\{(\varphi, x)\colon \,\varphi\in\frakg^\ast,\,x\in\frakg\}, \]
\[ [(\varphi_1,x_1), (\varphi_2,x_2)]
  = (x_1\varphi_2-x_2\varphi_1, [x_1,x_2]).
\]

In the case where $\frakg$ is the algebra $\mathfrak{so}(3)$ of infinitesimal
symmetries of $\bbR^3$, its dual $\frakg^\ast$ is $\bbR^3$ itself with the
usual action of $\mathfrak{so}(3)$ on it, and $I\frakg$ is the algebra $\bbR^3\rtimes
\mathfrak{so}(3)$ of infinitesimal affine isometries of $\bbR^3$. This is the
Lie algebra of the Euclidean group of isometries of $\bbR^3$, which is
often denoted $ISO(3)$. This explains our choice of the name $I\frakg$.

Note that, if $\frakg$ is a co-commutative Lie bialgebra, then $I\frakg$
is the ``double'' of $\frakg$~\cite{Drinfeld:QuantumGroups}. This
is a significant observation, for it is a part of the relationship
between this paper and the Etingof-Kazhdan theory of quantization of
Lie bialgebras~\cite{EtingofKazhdan:BialgebrasI}. Yet we will make no
explicit use of this observation below.

In the construction that follows we are going to define a map
from $\calA^w$ to $\glos{\calU}(I\frakg)$, the universal enveloping algebra of
$I\frakg$.  Note that a map ${\calA}^w
\to \calU(I\frakg)$ is ``almost the same'' as a map $\calA^{sw} \to
\calU(I\frakg)$, in the following sense.  The quotient
map $p\colon {\calA}^w \to \calA^{sw}$ has a one-sided inverse
$F\colon  \calA^{sw} \to {\calA}^w$ defined by
\[ F(D)= \sum_{k=0}^\infty \frac{(-1)^k}{k!} S_L^k(D)\cdot w_1^k. \]
Here $S_L$ denotes the map that sends an arrow diagram to the sum of
all ways of deleting a left-going arrow, $S_L^k$ is $S_L$ applied $k$ times, 
and $w_1$ denotes the 1-wheel,
as shown in Figure~\ref{fig:AwGenerators}.  The reader can verify that
$F$ is well-defined, an algebra- and co-algebra homomorphism, and that
$p\circ F= id_{\calA^{sw}}$.

\subsubsection{The Construction} Fixing a finite dimensional Lie algebra
$\frakg$, we construct a map $\glos{\calT^w_\frakg}\colon
{\calA}^w\to\calU(I\frakg)$ which assigns to every arrow diagram $D$
an element of the universal enveloping algebra $\calU(I\frakg)$. As is
often the case in our subject, a picture of a typical example is worth
more than a formal definition:
\[
  \def\I{{$I$}}
  \def\B{{$B$}}
  \def\g{{$\frakg$}}
  \def\d{{$\frakg^\ast$}}
  \def\p{{$\frakg^\ast\otimes\frakg^\ast\otimes\frakg\otimes\frakg
    \otimes\frakg^\ast\otimes\frakg^\ast$}}
  \def\u{{$\calU(I\frakg)$}}
  \input figs/Twg.pstex_t
\]

In short, we break up the diagram $D$ into its constituent
pieces and assign a copy of the structure constants tensor
$B\in\frakg^\ast\otimes\frakg^\ast\otimes\frakg$ to each internal
vertex $v$ of $D$ (keeping an association between the tensor factors
in $\frakg^\ast\otimes\frakg^\ast\otimes\frakg$ and the edges emanating
from $v$, as dictated by the orientations of the edges and of the vertex
$v$ itself). We assign the identity tensor in
$\frakg^\ast\otimes\frakg$ to every arrow in $D$ that is not connected to an
internal vertex, and contract any pair of factors connected by a fully
internal arrow. The remaining tensor factors
($\frakg^\ast\otimes\frakg^\ast\otimes\frakg\otimes\frakg
\otimes\frakg^\ast\otimes\frakg^\ast$ in our examples) are all along the
skeleton and can thus be ordered by the skeleton. We then multiply these
factors to get an output $\calT^w_\frakg(D)$ in $\calU(I\frakg)$.

It is also useful to restate this construction given a choice of a basis.
Let $\glos{(x_j)}$ be a basis of $\frakg$ and let $\glos{(\varphi^i)}$
be the dual basis of $\frakg^\ast$, so that $\varphi^i(x_j)=\delta^i_j$,
and let $\glos{b_{ij}^k}$ denote the structure constants of $\frakg$ in
the chosen basis: $[x_i,x_j]=\sum b_{ij}^kx_k$. Mark every arrow in $D$
with lower case Latin letter from within\footnote{The
supply of these can be made inexhaustible by the addition of numerical
subscripts.} $\{i,j,k,\dots\}$. Form a product $P_D$ by taking one $b_{\alpha\beta}^\gamma$
factor for each internal vertex $v$ of $D$ using the letters marking the
edges around $v$ for $\alpha$, $\beta$ and $\gamma$ and by taking one
$x_\alpha$ or $\varphi^\beta$ factor for each skeleton vertex of $D$,
taken in the order that they appear along the long line skeleton, with the indices
$\alpha$ and $\beta$ dictated by the edge markings and with the choice
between factors in $\frakg$ and factors in $\frakg^\ast$ dictated by the
orientations of the edges. Finally let $\calT^w_\frakg(D)$ be the sum
of $P_D$ over the indices $i,j,k,\dots$ running from $1$ to $\dim\frakg$:

\begin{equation} \label{eq:Twb}
  \def\P{{$\displaystyle
    \sum_{i,j,k,l,m,n=1}^{\dim\frakg}
    \hspace{-4mm} b_{ij}^kb_{kl}^m
    \varphi^i\varphi^jx_nx_m\varphi^l\in\calU(I\frakg)
  $}}
  \input figs/Twb.pstex_t
\end{equation}

The following is easy to verify (compare with~\cite[Theorem~4 and
Exercise~5.1]{Bar-Natan:OnVassiliev}):

\begin{proposition} The above two definitions of $T^w_\frakg$ agree, are
independent of the choices made within them, and respect all the relations
defining ${\calA}^w$. \qed
\end{proposition}

While we do not provide a proof of this proposition here, it is worthwhile
to state the correspondence between the relations defining ${\calA}^w$
and the Lie algebraic information in $\calU(I\frakg)$: $\aAS$ is
the antisymmetry of the bracket of $\frakg$, $\aIHX$ is the Jacobi
identity of $\frakg$, $\aSTU_1$ and  $\aSTU_2$ are the relations
$[x_i,x_j]=x_ix_j-x_jx_i$ and $[\varphi^i,x_j]=\varphi^ix_j-x_j\varphi^i$
in $\calU(I\frakg)$, $TC$ is the fact that $\frakg^\star$ is taken as an
abelian algebra, and $\aft$ is the fact that the identity tensor in
$\frakg^\ast\otimes\frakg$ is $\frakg$-invariant.

\subsubsection{Example: The 2 Dimensional Non-Abelian Lie Algebra}
Let $\frakg$ be the Lie algebra with two generators $x_{1,2}$
satisfying $[x_1,x_2]=x_2$, so that the only non-vanishing structure
constants $b_{ij}^k$ of $\frakg$ are $b_{12}^2=-b_{21}^2=1$. Let
$\varphi^i\in\frakg^\ast$ be the dual basis of $x_i$; by an easy
calculation, we find that in $I\frakg$ the element $\varphi^1$ is
central, while $[x_1,\varphi^2]=-\varphi^2$ and
$[x_2,\varphi^2]=\varphi^1$. We calculate $\calT^w_\frakg(D_L)$,
$\calT^w_\frakg(D_R)$ and $\calT^w_\frakg(w_k)$ using the ``in basis''
technique of Equation~\eqref{eq:Twb}. The outputs of these calculations lie
in $\calU(I\frakg)$; we display these results in a PBW basis in which the
elements of $\frakg^\ast$ precede the elements of $\frakg$:

\begin{eqnarray} 
  \calT^w_\frakg(D_L)
    &=& x_1\varphi^1+x_2\varphi^2 =
      \varphi^1x_1+\varphi^2x_2+[x_2,\varphi^2]
      = \varphi^1x_1+\varphi^2x_2+\varphi_1, \notag \\
  \calT^w_\frakg(D_R) &=& \varphi^1x_1+\varphi^2x_2, \label{eq:2DExample} \\
  \calT^w_\frakg(w_k) &=& (\varphi^1)^k. \notag
\end{eqnarray}

\parpic[r]{$\input figs/4wheel.pstex_t$}
For the last assertion above, note that all non-vanishing structure
constants $b_{ij}^k$ in our case have $k=2$, and therefore all indices
corresponding to edges that exit an internal vertex must be set equal to
$2$. This forces the ``hub'' of $w_k$ to be marked $2$ and therefore the
legs to be marked $1$, and therefore $w_k$ is mapped to $(\varphi^1)^k$.

Note that the calculations in~\eqref{eq:2DExample} are consistent with the
relation $D_L-D_R=w_1$ of Theorem~\ref{thm:Aw} and that they show that
other than that relation, the generators of ${\calA}^w$ are linearly
independent.

\draftcut \subsection{The Alexander Polynomial} \label{subsec:Alexander}
Let
$K$ be a long w-knot, let $Z(K)$ be the invariant of
Theorem~\ref{thm:ExpansionForKnots}. Theorem~\ref{thm:Alexander} below
asserts that apart from self-linking, $Z(K)$ contains precisely the
same information as the Alexander polynomial $A(K)$ of $K$ (recalled
below). But we have to start with some definitions.

\begin{figure}
\begin{center}
  $\input figs/8-17.pstex_t$
  \qquad
  \raisebox{-18mm}{\includegraphics[height=40mm]{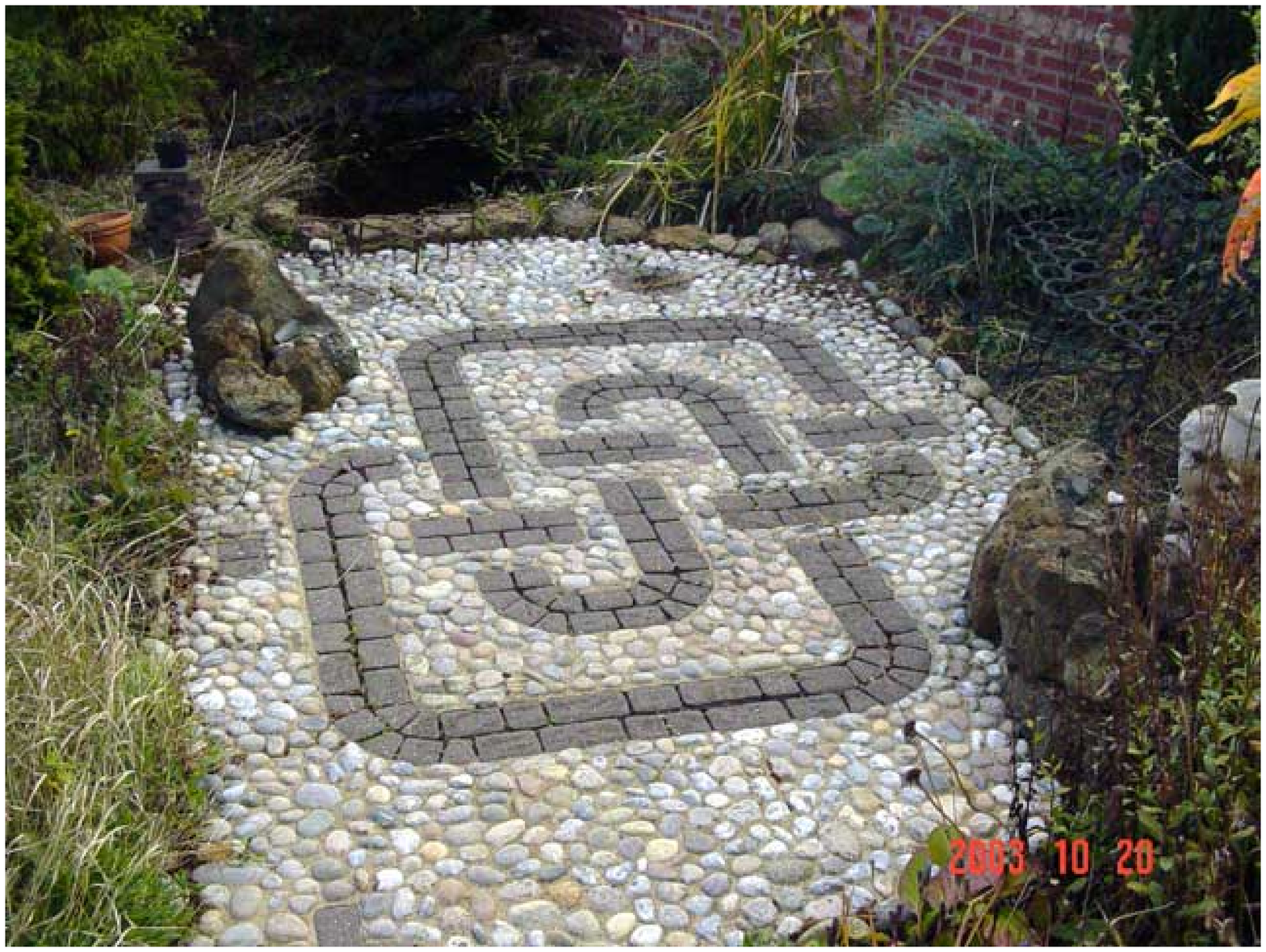}}
\end{center}
\caption{
  A long $8_{17}$, with the span of crossing $\#3$
  marked.  The projection is as in Brian Sanderson's garden.
  See~\cite{WKO}/\href{http://www.math.toronto.edu/~drorbn/papers/WKO/SandersonsGarden.html}{\tt
  SandersonsGarden.html}.
} \label{fig:817}
\end{figure}

\begin{definition} \label{def:STA} Enumerate the crossings of $K$
from $1$ to $n$ in some arbitrary order. For \linebreak $1\leq j\leq n$, the
``span'' of crossing $\#i$ is the connected open interval along the line
parametrizing $K$ between the two times $K$ ``visits'' crossing $\#i$
(see Figure~\ref{fig:817}). Form a matrix $T=\glos{T(K)}$ with $T_{ij}$
the indicator function of ``the lower strand of crossing $\#j$ is within
the span of crossing $\#i$'' (so $T_{ij}$ is $1$ if for a given $i,j$
the quoted statement is true, and $0$ otherwise). Let $\glos{s_i}$ be the
sign of crossing $\#i$ (recall that $\overcrossing$ is positive, 
$\undercrossing$ is negative; $(-,-,-,-,+,+,+,+)$ for Figure~\ref{fig:817}),
let $\glos{d_i}$ be $+1$ if $K$ visits the ``over'' strand of crossing
$\#i$ before visiting the ``under'' strand of that crossing, and let
$d_i=-1$ otherwise ($(-,+,-,+,-,+,-,+$) for Figure~\ref{fig:817}). Let
$S=\glos{S(K)}$ be the diagonal matrix with $S_{ii}=s_id_i$, and for
an indeterminate $\glos{X}$, let $X^{-S}$ denote the diagonal matrix
with diagonal entries $X^{-s_id_i}$.  Finally, let $\glos{A(K)}$ be the
Laurent polynomial in $\bbZ[X,X^{-1}]$ given by
\begin{equation} \label{eq:AKDef}
  A(K)(X) := \det\left(I+T(I-X^{-S})\right).
\end{equation}
\end{definition}

\begin{example} For the knot diagram in Figure~\ref{fig:817},
\[ \scriptstyle
  T = \left(\begin{smallmatrix}
    0 & 1 & 1 & 1 & 1 & 0 & 1 & 0 \\
    0 & 0 & 1 & 0 & 1 & 0 & 0 & 0 \\
    0 & 1 & 0 & 0 & 1 & 0 & 0 & 0 \\
    0 & 1 & 0 & 0 & 1 & 0 & 1 & 0 \\
    0 & 1 & 0 & 1 & 0 & 1 & 1 & 1 \\
    0 & 1 & 0 & 1 & 0 & 0 & 1 & 0 \\
    0 & 0 & 0 & 1 & 0 & 1 & 0 & 0 \\
    0 & 0 & 0 & 1 & 0 & 1 & 0 & 0
  \end{smallmatrix}\right),
  \quad
  S = \left(\begin{smallmatrix}
    1 & 0 & 0 & 0 & 0 & 0 & 0 & 0 \\
    0 & -1 & 0 & 0 & 0 & 0 & 0 & 0 \\
    0 & 0 & 1 & 0 & 0 & 0 & 0 & 0 \\
    0 & 0 & 0 & -1 & 0 & 0 & 0 & 0 \\
    0 & 0 & 0 & 0 & -1 & 0 & 0 & 0 \\
    0 & 0 & 0 & 0 & 0 & 1 & 0 & 0 \\
    0 & 0 & 0 & 0 & 0 & 0 & -1 & 0 \\
    0 & 0 & 0 & 0 & 0 & 0 & 0 & 1
  \end{smallmatrix}\right),
  \quad\text{and}\quad
  A = \left|\begin{smallmatrix}
    1 & 1-X & 1-X^{-1} & 1-X & 1-X & 0 & 1-X & 0 \\
    0 & 1 & 1-X^{-1} & 0 & 1-X & 0 & 0 & 0 \\
    0 & 1-X & 1 & 0 & 1-X & 0 & 0 & 0 \\
    0 & 1-X & 0 & 1 & 1-X & 0 & 1-X & 0 \\
    0 & 1-X & 0 & 1-X & 1 & 1-X^{-1} & 1-X & 1-X^{-1} \\
    0 & 1-X & 0 & 1-X & 0 & 1 & 1-X & 0 \\
    0 & 0 & 0 & 1-X & 0 & 1-X^{-1} & 1 & 0 \\
    0 & 0 & 0 & 1-X & 0 & 1-X^{-1} & 0 & 1
  \end{smallmatrix}\right|.
\]
The last determinant equals $-X^3+4X^2-8X+11-8X^{-1}+4X^{-2}-X^{-3}$,
the Alexander polynomial of the knot $8_{17}$
(see e.g.~\cite{Rolfsen:KnotsAndLinks}).
\end{example}

\begin{theorem} \label{thm:AlexanderFormula}
(Lee,~\cite[Theorem~1]{Lee:AlexanderInvariant}) For any (classical)
knot $K$, $A(K)$ is equal to the normalized Alexander
polynomial~\cite{Rolfsen:KnotsAndLinks} of $K$. \qed
\end{theorem}

The Mathematica notebook~\cite[``wA'']{WKO} verifies
Theorem~\ref{thm:AlexanderFormula} for all prime knots with up to 11
crossings.

The following theorem asserts that $Z(K)$ can be computed from $A(K)$
(see Equation~\eqref{eq:AtoZ}) and that modulo a certain additional relation
and with the appropriate identifications in place, $Z(K)$ {\em is} $A(K)$
(see Equation~\eqref{eq:ZisA}).

\begin{theorem} \label{thm:Alexander} (Proof in
Section~\ref{subsec:AlexanderProof}). Let $x$ be an indeterminate, let $\sl$
be self-linking as in Exercise~\ref{ex:sl}, let $D_A:=D_L=D_R$ and $w_k$
be as in Figure~\ref{fig:AwGenerators}, and let $\glos{w}\colon
\bbQ\llbracket x\rrbracket \to\calA^w$ be the linear map defined by
$x^k\mapsto w_k$. Then for a long w-knot $K$,
\begin{equation} \label{eq:AtoZ}
  Z(K) = 
    \underbrace{
      \exp_{\calA^{sw}}\left(\sl(K)D_A\right)
    }_\text{$\sl$ coded in arrows} \cdot
    \underbrace{
      \exp_{\calA^{sw}}\left(-w\left(\log_{\bbQ\llbracket x\rrbracket}
        A(K)(e^x)
      \right)\right)
    }_\text{main part: Alexander coded in wheels},
\end{equation}
where the logarithm and inner exponentiation are computed by formal power
series in $\bbQ\llbracket x\rrbracket$ and the outer exponentiations
are likewise computed in $\calA^{sw}$.
\end{theorem}

\parpic[r]{$\input figs/wkl.pstex_t$}
Let $\calA^\text{reduced}$ be $\calA^{sw}$ modulo the additional
relations $D_A=0$ and $w_kw_l=w_{k+l}$ for $k,l\neq 1$. The quotient
$\calA^\text{reduced}$ can be identified with the vector space of (infinite)
linear combinations of $w_k$'s (with $k\neq 1$).  Identifying the
$k$-wheel $w_k$ with $x^k$, we see that $\calA^\text{reduced}$ is
the space of power series in $x$ having no linear terms. Note by
inspecting Equation~\eqref{eq:AKDef} that $A(K)(e^x)$ never has a term linear
in $x$, and that modulo $w_kw_l=w_{k+l}$, the exponential and the
logarithm in Equation~\eqref{eq:AtoZ} cancel each other out. Hence, within
$\calA^\text{reduced}$,

\begin{equation} \label{eq:ZisA} Z(K) = A^{-1}(K)(e^x). \end{equation}

\begin{remark} In~\cite{HabiroKanenobuShima:R2K} Habiro, Kanenobu,
and Shima show that all coefficients of the Alexander polynomial are
finite type invariants of w-knots, and in~\cite{HabiroShima:R2KII}
Habiro and Shima show that all finite type invariants of w-knots are
polynomials in the coefficients of the Alexander polynomial. Thus,
Theorem~\ref{thm:Alexander} is merely an ``explicit form'' of these earlier
results.
\end{remark}

\draftcut
\subsection{Proof of Theorem~\ref{thm:Alexander}}
\label{subsec:AlexanderProof}

We start with a sketch. The proof of Theorem~\ref{thm:Alexander} can be
divided in three parts: differentiation, bulk management, and computation.

\noindent{\bf Differentiation.} 
Both sides of our goal, that is,
Equation~\eqref{eq:AtoZ}, are exponential in nature. When seeking to
show an equality of exponentials it is often beneficial to compare
their derivatives\footnote{Thanks, Dylan.}. In our case the useful
``derivatives'' to use are the ``Euler operator'' $\glos{E}$ (``multiply
every term by its degree'', an analogue of $f\mapsto xf'$, defined
in Section~\ref{subsubsec:Euler}), and the ``normalized Euler
operator'' $Z\mapsto\glos{\tilE} Z:=Z^{-1}EZ$, which is a variant of the
logarithmic derivative $f\mapsto x(\log f)'=xf'/f$. Since $\tilE$
is one to one (see Section~\ref{subsubsec:Euler}) and since we know how
to apply $\tilE$ to the right hand side of Equation~\eqref{eq:AtoZ}
(see Section~\ref{subsubsec:Euler}), it is enough to show that with
$\glos{B}:=T(\exp(-xS)-I)$ and suppressing the fixed w-knot $K$ from the
notation,
\begin{equation} \label{eq:EofAtoZ}
  EZ = Z\cdot\left(
    \sl\cdot D_A-w\!\left[x\tr\left( (I-B)^{-1}TS\exp(-xS) \right)\right]
  \right) \qquad \text{ in }\calA^{sw}.
\end{equation}

\noindent{\bf Bulk Management.}
Next we seek to understand the left hand
side of Equation~\eqref{eq:EofAtoZ}. $Z$ is made up of ``quantities in bulk'':
arrows that come in exponential ``reservoirs''. As it turns out,
$EZ$ is made up of the same bulk quantities, but also allowing for a
single non-bulk ``excitation'', which we often highlight in red. (compare
with $Ee^x={\red x}\cdot e^x$; the
``bulk'' $e^x$ remains, and single ``excited red'' $\red x$ gets created). We
wish manipulate and simplify that red excitation. This is best done by
introducing a certain module, $\glos{\IAM_K}$, the ``Infinitesimal Alexander
Module'' of $K$ (see Section~\ref{subsubsec:IAM}). The elements of $\IAM_K$
can be thought of as names for ``bulk objects with a red excitation'',
and hence, there is an ``interpretation map'' $\glos{\iota}\colon
\IAM_K\to\calA^{sw}$, which maps every ``name'' into the object it
represents. There are three special elements in $\IAM_K$: an element
$\glos{\lambda}$, which is the name of $EZ$ (that is, $\iota(\lambda)=EZ$),
the element $\glos{\delta_A}$ which is the name of $D_A\cdot Z$
(so $\iota(\delta_A)=D_A\cdot Z$), and an element $\glos{\omega_1}$
which is the name of a ``detached'' 1-wheel that is appended to
$Z$. The latter can take a coefficient which is a power of $x$,
with $\iota(x^k\omega_1)=w(x^{k+1})\cdot Z=(Z\text{ times a }
(k+1)\text{-wheel})$. Thus, it is enough to show that in $\IAM_K$,
\begin{equation} \label{eq:GoalInIAM}
  \lambda = \sl\cdot\delta_A
    - \tr\left((I-B)^{-1}TSX^{-S}\right)\omega_1,
  \quad\text{with}\quad X=e^x.
\end{equation}
Indeed, applying $\iota$ to both sides of the above equation, we get
Equation~\eqref{eq:EofAtoZ} back again.

\noindent{\bf Computation.} Last, we show in
Section~\ref{sec:ComputeLambda} that Equation~\eqref{eq:GoalInIAM} holds true. This
is a computation that happens entirely in $\IAM_K$ and does not mention
finite type invariants, expansions or arrow diagrams in any way.

\subsubsection{The Euler Operator} \label{subsubsec:Euler} Let $A$ be
a completed, graded algebra with unit, in which all degrees are $\geq
0$. Define a continuous linear operator $E\colon A\to A$ by setting $Ea=(\deg
a)a$ for homogeneous $a\in A$. In the case $A=\bbQ\llbracket x\rrbracket$,
we have $Ef=xf'$, the standard ``Euler operator'': indeed, for each $n$, 
$Ex^n=nx^n=x\cdot(x^n)'$. Hence, we adopt
the name $E$ for this operator in general.

We say that $Z\in A$ is a ``perturbation of the identity'' if its
degree 0 piece is 1. Such a $Z$ is always invertible. For such a $Z$,
set $\tilE Z:=Z^{-1}\cdot EZ$, and call the thus (partially) defined
operator $\tilE \colon A\to A$ the ``normalized Euler operator''. From this
point on when we write $\tilE Z$ for some $Z\in A$, we automatically
assume that $Z$ is a perturbation of the identity or that it is trivial
to show that $Z$ is a perturbation of the identity. Note that for
$f\in\bbQ\llbracket x\rrbracket$, we have $\tilE f=x(\log f)'$,
so $\tilE$ is a variant of the logarithmic derivative.

\begin{claim} $\tilE$ is one to one.
\end{claim}

\begin{proof} Assume $Z_1\neq Z_2$ and let $d$ be the smallest degree
in which they differ. Then $d>0$ and in degree $d$ the difference
$\tilE Z_1-\tilE Z_2$ is $d$ times the difference $Z_1-Z_2$, and
hence, $\tilE Z_1\neq\tilE Z_2$. \qed
\end{proof}

Thus, in order to prove our goal, that is, Equation~\eqref{eq:AtoZ}, it is enough to
compute $\tilE$ of both sides and to show the equality then. We start
with the right hand side of Equation~\eqref{eq:AtoZ}; but first, we need some
simple properties of $E$ and $\tilE$. The proofs of these properties are
routine, and hence, they are omitted.

\begin{proposition} The following hold true:
\begin{enumerate}
\item $E$ is a derivation: $E(fg)=(Ef)g+f(Eg)$.
\item If $Z_1$ commutes with $Z_2$, then $\tilE(Z_1Z_2)=\tilE Z_1+\tilE Z_2$.
\item If $z$ commutes with $Ez$, then $Ee^z=e^z(Ez)$ and $\tilE e^z=Ez$.
\item If $w\colon A\to\calA$ is a morphism of graded algebras,
then it commutes with $E$ and $\tilE$. \qed
\end{enumerate}
\end{proposition}

Let us denote the right hand side of Equation~\eqref{eq:AtoZ} by $Z_1(K)$. Then, by
the above proposition, remembering (see Theorem~\ref{thm:Aw}) that $\calA^{sw}$ is
commutative and that $\deg D_A=1$, we have
\[ \tilE Z_1(K) = \sl\cdot D_A-w(E\log A(K)(e^x))
  = \sl\cdot D_A-w\left(x\frac{d}{dx}\log A(K)(e^x)\right).
\]
The rest is an exercise in matrices and
differentiation. $A(K)$ is a determinant, see Equation~\eqref{eq:AKDef}, and in general,
$\frac{d}{dx}\log\det(M) = \tr\left(M^{-1}\frac{d}{dx}M\right)$. So with
$B=T(e^{-xS}-I)$ (so $M=I-B$), we have
\[ \tilE Z_1(K) =
  \sl\cdot D_A + w\left(x\tr\left((I-B)^{-1}\frac{d}{dx}B\right)\right)
  = \sl\cdot D_A - w\left(x\tr\left((I-B)^{-1}TSe^{-xS}\right)\right),
\]
as promised in Equation~\eqref{eq:EofAtoZ}.

\subsubsection{The Infinitesimal Alexander Module} \label{subsubsec:IAM}
Let $K$ be a w-knot diagram. The ``Infinitesimal Alexander Module'' $\IAM_K$
of $K$, which is defined in detail below, is a certain module made from
a certain space $\glos{\IAM^0_K}$ of pictures ``annotating'' $K$ with
``red excitations'' modulo some pictorial relations that indicate how
the red excitations can be moved around. The space $\IAM^0_K$ in itself
is made of three pieces, or ``sectors''. The ``A sector'' in which the
excitations are red arrows, the ``Y sector'' in which the excitations
are ``red hairy Y-diagrams'', and a rank 1 ``W sector'' for ``red
hairy wheels''. There is an ``interpretation map'' $\glos{\iota}\colon
\IAM^0_K\to\calA^w$ which descends to a well-defined (and homonymous)
$\iota\colon \IAM_K\to\calA^w$. Finally, there are some special elements
$\lambda$ and $\delta_A$ that live in the A sector of $\IAM^0_K$ and
$\omega_1$ that lives in the W sector.

In principle, the description of $\IAM^0_K$ and of $\IAM_K$ can be given
independently of the interpretation map $\iota$, and there are some good
questions to ask about $\IAM_K$ (and the special elements in it) that are
completely independent of the interpretation of the elements of $\IAM_K$ as
``perturbed bulk quantities'' within $\calA^{sw}$. Yet $\IAM_K$ is a
complicated object and we fear its definition will appear completely
artificial without its interpretation. Hence, below the two definitions will
be woven together.

\begin{figure}
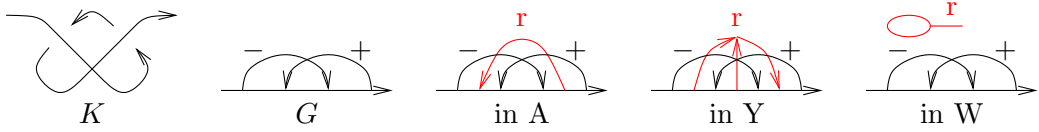

\[ \input figs/IAM0Def.pstex_t \]
\caption{
  A sample w-knot $K$, it's Gauss diagram $G$, and one generator from
  each of the A, Y, and W sectors of $\IAM^0_K$. Red parts are marked
  with the letter ``r''.
} \label{fig:IAM0Def}
\end{figure}

$\IAM_K$ and $\iota$ may equally well be described in terms of $K$ or in
terms of the Gauss diagram of $K$ (see Remark~\ref{rem:GD}). For pictorial
simplicity, we choose to use the latter; so let $G=G(K)$ be the Gauss
diagram of $K$. It is best to read the following definition while at the
same time studying Figure~\ref{fig:IAM0Def}.

\begin{definition} Let $\glos{R}$ be the ring $\bbZ[X,X^{-1}]$ of Laurent
polynomials in a variable $X$ with integer coefficients\footnote{Later,
$X$ is interpreted in $\calA^w$ as a formal exponential $e^x$. So within $\IAM$
we can restrict to coefficients in $\bbZ$, yet in $\calA^w$ we
must allow coefficients in $\bbQ$.}, and let $\glos{R_1}$ be the subring
of polynomials that vanish at $X=1$ (i.e., whose sum of coefficients
is $0$)\footnote{$R_1$ is only very lightly needed, and only within
Definition~\ref{def:InterpretationMap}. In particular, all that we say
about $\IAM_K$ that does not concern the interpretation map $\iota$ is
equally valid with $R$ replacing $R_1$.}.  Let $\IAM^0_K$ be the direct
sum of the following three modules (which for the purpose of taking the
direct sum, are all regarded as $\bbZ$-modules):
\begin{enumerate}
\item The ``A sector'' is the free $\bbZ$-module generated by all diagrams
made from $G$ by the addition of a single unmarked ``red excitation''
arrow, whose endpoints are on the long line skeleton of $G$ and are distinct from
each other and from all other endpoints of arrows in $G$. Such diagrams
are considered combinatorially --- so two are equivalent iff they differ
only by an orientation preserving diffeomorphism of the skeleton. Let
us count: if $K$ has $n$ crossings, then $G$ has $n$ arrows and the
skeleton of $G$ get subdivided into $m:=2n+1$ arcs. An A sector diagram
is specified by the choice of an arc for the tail of the red arrow and
an arc for the head ($m^2$ choices), except if the head and the tail
fall within the same arc, their relative ordering has to be specified
as well ($m$ further choices). So the rank of the A sector over $\bbZ$
is $m(m+1)$.
\item The ``Y sector'' is the free $R_1$-module generated by all
diagrams made from $G$ by the addition of a single ``red excitation''
$Y$-shape single-vertex graph, with two incoming edges (``tails'') and
one outgoing (``head''), modulo anti-symmetry for the two incoming edges
(again, considered combinatorially). Counting is more elaborate: when
the three edges of the $Y$ end in distinct arcs in the skeleton of $G$,
we have $\frac12m(m-1)(m-2)$ possibilities ($\frac12$ for the
antisymmetry). When the two tails of the $Y$ lie on the same arc, we get $0$
by anti-symmetry. The remaining possibility is to have the head and
one tail on one arc (order matters!) and the other tail on another,
at $2m(m-1)$ possibilities. So the rank of the Y sector over $R_1$
is $m(m-1)(\frac12m+1)$.
\item The ``W sector'' is the rank 1 free $R$-module with a single
generator $w_1$. It is not necessary for $w_1$ to have a pictorial
representation, yet one, involving a single ``red'' 1-wheel, is shown in
Figure~\ref{fig:IAM0Def}. This pictorial representation is consistent
with the interpretation in the definition below of $\omega_1$ as a
detached 1-wheel.
\end{enumerate}
\end{definition}

\begin{definition} \label{def:InterpretationMap}
The ``interpretation map'' $\iota\colon \IAM^0_K\to\calA^w$
is defined by sending the arrows (marked $+$ or $-$) of a diagram in
$\IAM^0_K$ to $(e^{\pm a})$-exponential reservoirs of arrows, as in the
definition of $Z$ (see Remark~\ref{rem:ZwForGD}). In addition, the red
excitations of diagrams in $\IAM^0_K$ are interpreted as follows:
\begin{enumerate}
\item In the A sector, the red arrow is simply mapped to itself, with the
colour red suppressed.
\item In the $Y$ sector diagrams have red $Y$'s and coefficients $f\in
R_1$. Substitute $X=e^x$ in $f$, expand in powers of $x$,
and interpret $x^kY$ as a ``hairy $Y$ with $k-1$ hairs'' as in
Exercise~\ref{ex:Hair}. Note that $f(1)=0$, so only positive powers of $x$
occur. So we never need to worry about ``$Y$'s with $-1$ hairs''. This is
the only point where the condition $f\in R_1$ (as opposed to $f\in R$) is
needed.
\item In the $W$ sector treat the coefficients as above, but interpret
$x^kw_1$ as a detached $w_{k+1}$. I.e., as a detached wheel with $k+1$
spokes, as in Exercise~\ref{ex:Hair}.
\end{enumerate}
\end{definition}

As stated above, $\IAM_K$ is the quotient of $\IAM^0_K$ by some set of
relations. The best way to think of this set of relations is as
``everything that's obviously annihilated by $\iota$''. Here's the same
thing, in a more formal language:

\begin{figure}
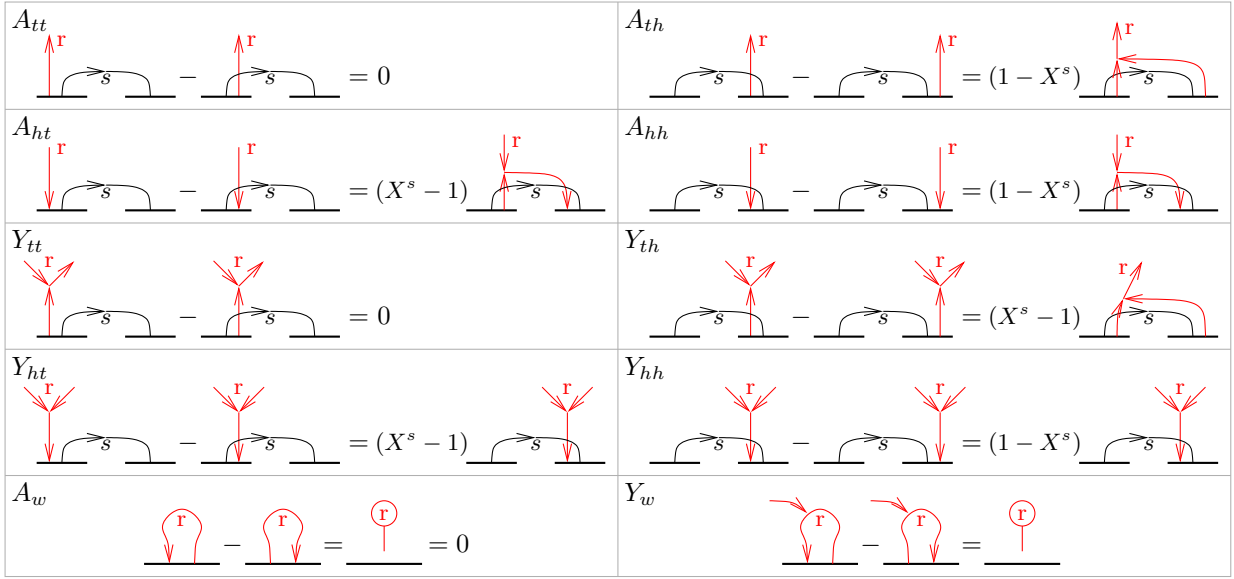

\[ \input figs/IAMRelations.pstex_t \]
\caption{The relations $\calR$ making $\IAM_K$.} \label{fig:IAMRelations}
\end{figure}

\begin{definition} Let $\glos{\IAM_K}:=\IAM^0_K/\calR$, where
$\glos{\calR}$ is the linear span of the relations depicted in
Figure~\ref{fig:IAMRelations}. The top 8 relations are about moving
a leg of the red excitation across an arrow head or an arrow tail in
$G$. Since the red excitation may be either an arrow $A$ or a $Y$,
its leg in motion may be either a tail or a head, and it may be moving
either past a tail or past a head, there are 8 relations of that type. The
$A_w$ relation corresponds to $D_L-D_R=w_1=0$. The $Y_w$ relation indicates
the ``price'' (always a red $w_1$) of commuting a red head across a red
tail. As per custom, in each case only the changing part of the diagrams
involved is shown. Further, the red excitations are marked with the
letter ``r'' and the sign of an arrow in $G$ is marked $s$; so always
$s\in\{\pm 1\}$. The relations in the left column may be multiplied
by a scalar in $\bbZ$, while the relations in the right column may be
multiplied by a scalar in $R$. Hence, for example, $x^0w_1=0$ by $A_w$,
yet $x^kw_1\neq 0$ for $k>0$.
\end{definition}

\begin{proposition} The interpretation map $\iota$ indeed annihilates all
the relations in $\calR$.
\end{proposition}

\begin{proof} $\iota A_{tt}$ and $\iota Y_{tt}$ follow immediately from
the TC relation. The formal identity $e^{\ad b}(a)=e^bae^{-b}$ (here $\ad$ denotes 
the adjoint representation) implies
$e^{\ad b}(a)e^b=e^ba$, and hence, $ae^b-e^ba=(1-e^{\ad b})(a)e^b$. With
$a$ interpreted as ``red head'', $b$ as ``black head'', and $\ad b$
as ``hair'' (justified by the $\iota$-meaning of hair and by the
$\aSTU_1$ relation, see Figure~\ref{fig:aSTU}), the last equality becomes
a proof of $\iota Y_{hh}$.  Further pushing that same equality, we get
$ae^b-e^ba=\frac{1-e^{\ad b}}{\ad b}([b,a])$, where $\frac{1-e^{\ad
b}}{\ad b}$ is first interpreted as a power series $\frac{1-e^y}{y}$
involving only non-negative powers of $y$, and then the substitution
$y=\ad b$ is made. But that's $\iota A_{hh}$, when one remembers
that $\iota$ on the Y sector automatically contains a single
``$\frac{1}{\text{hair}}$'' factor. Similar arguments, though using
$\aSTU_2$ instead of $\aSTU_1$, prove that $Y_{ht}$, $Y_{th}$, $A_{ht}$,
and $A_{th}$ are all in $\ker\iota$. Finally, $\iota A_w$ is RI,
and $\iota Y_w$ is a direct consequence of $\aSTU_2$. \qed
\end{proof}

Finally, we come to the special elements $\lambda$, $\delta_A$, and $\omega_1$.

\begin{figure}
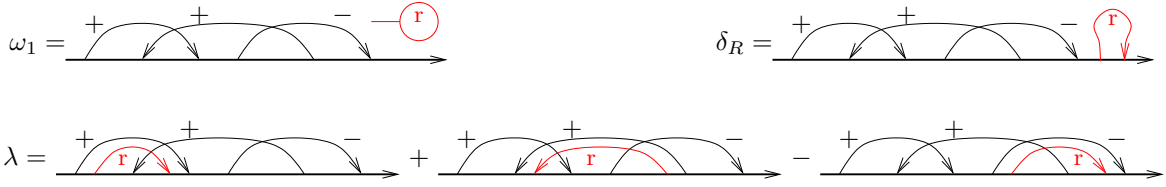

\[ \input figs/SpecialElements.pstex_t \]
\caption{The special elements $\omega_1$, $\delta_A$, and $\lambda$
  in $\IAM_G$, for a sample 3-arrow Gauss diagram $G$.
} \label{fig:SpecialElements}
\end{figure}

\begin{definition} Within $\IAM_G$, let $\omega_1$ be, as before, the
generator of the W sector. Let $\delta_A$ be a ``short''
red arrow, as in the $A_w$ relation (exercise:
modulo $\calR$, this is independent of the placement of the short
arrows within $G$). Finally, let $\lambda$ be the signed sum of exciting
each of the (black) arrows in $G$ in turn. The picture says all, and it is
Figure~\ref{fig:SpecialElements}.
\end{definition}

\begin{lemma} In $\calA^{sw}(\uparrow)$, the special elements of
$\IAM_G$ are interpreted as follows: \linebreak $\iota(\omega_1)=Zw_1$,
$\iota(\delta_A)=ZD_A$, and most interesting, $\iota(\lambda)=EZ$.
Therefore, Equation~\eqref{eq:GoalInIAM} (if true) implies
Equation~\eqref{eq:EofAtoZ} and hence, it implies our goal,
Theorem~\ref{thm:Alexander}.
\end{lemma}

\begin{proof} For the proof of this lemma, the only thing that isn't
done yet and isn't trivial is the assertion $\iota(\lambda)=EZ$. But this
assertion is a consequence of $Ee^{\pm a}=\pm ae^{\pm a}$ and of a
Leibniz law for the derivation $E$, appropriately generalized to a
context where $Z$ can be thought of as a ``product'' of ``arrow
reservoirs''. The details are left to the reader. \qed
\end{proof}

\subsubsection{The Computation of $\lambda$} \label{sec:ComputeLambda}

Naturally, our next task is to prove Equation~\eqref{eq:GoalInIAM}. This is
done entirely algebraically within the finite rank module $\IAM_G$. To
read this section one need not know about $\calA^{sw}(\uparrow)$, or $\iota$,
or $Z$, but we do need to lay out some notation. Start by marking the arrows
of $G$ with $a_1$ through $a_n$ in some order.

\begin{figure}
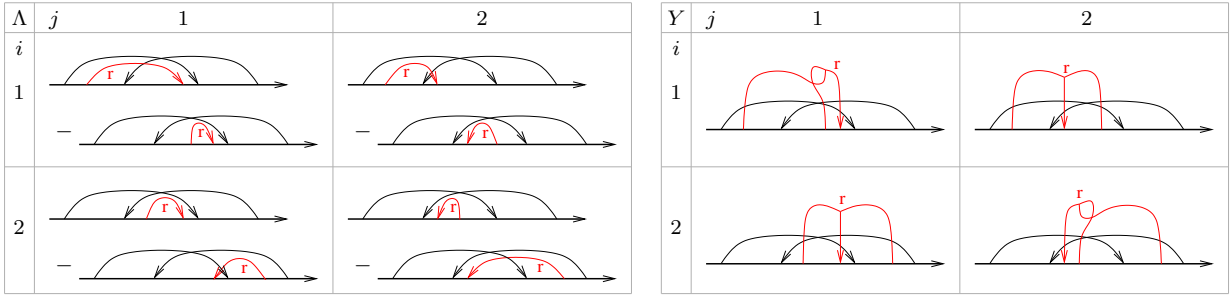

\[ \input figs/LambdaAndY.pstex_t \]
\caption{The matrices $\Lambda$ and $Y$ for a sample 2-arrow Gauss
  diagram (the signs on $a_1$ and $a_2$ are suppressed, and so are the $r$
  marks). The twists in $y_{11}$ and $y_{22}$ may be replaced by minus
  signs.
} \label{fig:LambdaAndY}
\end{figure}

Let $\epsilon$ stand for the informal yet useful quantity
``a little''. Let $\lambda_{ij}$ denote the difference
$\lambda'_{ij}-\lambda''_{ij}$ of red excitations in the A sector of
$\IAM_G$, where $\lambda'_{ij}$ is the diagram with a red arrow whose
tail is $\epsilon$ to the right of the left end of $a_i$ and whose head
is $\frac12\epsilon$ away from head of $a_j$ in the direction of the
tail of $a_j$, and where $\lambda''_{ij}$ has a red arrow whose tail
is $\epsilon$ to the left of the right end of $a_i$ and whose head is
as before, $\frac12\epsilon$ away from head of $a_j$ in the direction
of the tail of $a_j$.  Let $\Lambda=(\lambda_{ij})$ be the matrix whose
entries are the $\lambda_{ij}$'s, as shown in Figure~\ref{fig:LambdaAndY}.

Similarly, let $y_{ij}$ denote the element in the Y sector of $\IAM_G$
whose red Y has its head $\frac12\epsilon$ away from head of $a_j$
in the direction of the tail of $a_j$, its right tail (as seen
from the head) $\epsilon$ to the left of the right end of $a_i$ and
its left tail $\epsilon$ to the right of the left end of $a_i$. Let
$Y=(y_{ij})$ be the matrix whose entries are the $y_{ij}$'s, as shown
in Figure~\ref{fig:LambdaAndY}.

\begin{lemma} \label{lem:IAMStructure}
With $S$ and $T$ as in Definition~\ref{def:STA}, and with $B=T(X^{-S}-I)$
and $\lambda$ as above, the following identities between elements
of $\IAM_G$ and matrices with entries in $\IAM_G$ hold true:
\begin{eqnarray}
  \lambda-\sl\cdot D_A &=& \tr S\Lambda \label{eq:lambda}, \\
  \Lambda &=& -BY-TX^{-S}w_1 \label{eq:Lambda}, \\
  Y &=& BY + TX^{-S}w_1 \label{eq:Y}.
\end{eqnarray}
\end{lemma}

\noindent{\em Proof of Equation~\eqref{eq:GoalInIAM} given
Lemma~\ref{lem:IAMStructure}.} The last of the equalities above
implies that $Y=(I-B)^{-1}TX^{-S}w_1$. Thus,
\begin{align*}
  \lambda-\sl\cdot D_A = \tr S\Lambda = -\tr S(BY+TX^{-S}w_1) &=
    -\tr S(B(I-B)^{-1}TX^{-S}+TX^{-S})w_1 \\
  &= -\tr\left((I-B)^{-1}TSX^{-S}\right)w_1,
\end{align*}
and this is exactly Equation~\eqref{eq:GoalInIAM}. \qed

\noindent{\em Proof of Lemma~\ref{lem:IAMStructure}.}
Equation~\eqref{eq:lambda} is trivial. The proofs of
Equations~\eqref{eq:Lambda} and~\eqref{eq:Y} both have the same simple
cores, that have to be supplemented by highly unpleasant tracking of signs
and conventions and powers of $X$. Let us start from the cores.

To prove Equation~\eqref{eq:Lambda} we wish to ``compute''
$\lambda_{ik}=\lambda'_{ik}-\lambda''_{ik}$. As $\lambda'_{ik}$ and
$\lambda''_{ik}$ have their heads in the same place, we can compute their
difference by gradually sliding the tail of $\lambda'_{ik}$ from its
original position near the left end of $a_i$ towards the right end of
$a_i$, where it would be cancelled by $\lambda''_{ik}$. As the tail slides
we pick up a $y_{jk}$ term each time it crosses a head of an $a_j$ (relation
$A_{th}$), we pick up a vanishing term each time it crosses a tail
(relation $A_{tt}$), and we pick up a $w_1$ term if the tail needs to
cross over its own head (relation $A_w$). Ignoring signs and $(X^{\pm
1}-1)$ factors, the sum of the $y_{jk}$-terms should be proportional
to $TY$, for indeed, the matrix $T$ has non-zero entries precisely when
the head of an $a_j$ falls within the span of an $a_i$. Unignoring these
signs and factors, we get $-BY$ (recall that $B=T(X^{-S}-I)$ is just $T$
with added $(X^{\pm 1}-1)$ factors). Similarly, a $w_1$ term arises
in this process when a tail has to cross over its own head, that is,
when the head of $a_k$ is within the span of $a_i$. Thus, the $w_1$
term should be proportional to $Tw_1$, and we claim it is $-TX^{-S}w_1$.

The core of the proof of Equation~\eqref{eq:Y} is more or less the
same. We wish to ``compute'' $y_{ik}$ by sliding its left leg, starting
near the left end of $a_i$, towards its right leg, which is stationary
near the right end of $a_i$. When the two legs come together, we get 0
because of the anti-symmetry of Y excitations. Along the way we pick up
further Y terms from the $Y_{th}$ relations, and sometimes a $w_1$ term
from the $Y_w$ relation. When all signs and $(X^{\pm 1}-1)$ factors are
accounted for, we get Equation~\eqref{eq:Y}.

We leave it to the reader to complete the details in the above proofs. It
is a major headache, and we would not have trusted ourselves had we not
written a computer program to manipulate quantities in $\IAM_G$ by a
brute force application of the relations in $\calR$. Everything checks;
see~\cite[``The Infinitesimal Alexander Module'']{WKO}. \qed

This concludes the proof of Theorem~\ref{thm:Alexander}. \qed

\begin{remark} We chose the name ``Infinitesimal Alexander Module'' as in
our mind there is some similarity between $\IAM_K$ and the ``Alexander
Module'' of $K$. Yet beyond the above, we did not embark on any serious
study of $\IAM_K$. In particular, we do not know if $\IAM_K$ in itself
is an invariant of $K$ (though we suspect it wouldn't be hard to show
that it is), we do not know if $\IAM_K$ contains any further information
beyond $\sl$ and the Alexander polynomial, and we do not know if there is any
formal relationship between $\IAM_K$ and the Alexander module of $K$.
\end{remark}

\begin{remark} The logarithmic derivative of the Alexander polynomial
also appears in Lescop's work, see~\cite{Lescop:EquivariantLinking, Lescop:Cube}. We
don't know if its appearances there are related to its appearance here.
\end{remark}

\draftcut
\subsection{The Relationship with u-Knots} \label{subsec:RelWithKont}
Unlike in the case of braids, there is a canonical universal finite type
invariant of $u$-knots: the Kontsevich integral $\glos{Z^u}$. So it
makes sense to ask how it is related to the expansion $Z^w$.

\parpic[l]{$\xymatrix{
  \calK^u(\uparrow) \ar[r]^{Z^u} \ar[d]^a
    & \glos{\calA^u}(\uparrow) \ar[d]^\alpha \\
  \calK^w(\uparrow) \ar[r]^{Z^{w}}
    & \calA^{sw}(\uparrow)
}$}
We claim that the square on the left commutes, where
$\glos{\calK^u}(\uparrow)$ stands for long $u$-knots (knottings of
an oriented line), and similarly $\calK^w(\uparrow)$ denotes long
$w$-knots. As before, $a$ is the composition of the maps $u$-knots
$\to$ $v$-knots $\to$ $w$-knots, and $\alpha$ is the induced map on
the projectivizations, mapping each chord to the sum of the two ways to
direct it.

Recall that $\alpha$ kills everything but wheels and arrows (hence $Z^w$ is much
weaker, but also easier to handle, than the Kontsevich integral).
We are going to use the formula for the ``wheel part'' of the Kontsevich integral
as stated in \cite{Kricker:Kontsevich}. 
Let $K$ be a 0-framed long knot, and let $A(K)$ denote the Alexander polynomial. Then by \cite{Kricker:Kontsevich},
$$Z^u(K)= \exp_{\calA^u}\left(-\frac{1}{2} \log A(K)(e^h)|_{h^{2n}\to w^u_{2n}}\right)+\text{ ``loopy terms''},$$
where $w^u_{2n}$ stands for the unoriented wheel with $2n$ spokes; and ``loopy terms'' means terms that contain 
diagrams with more than one loop, which
are killed by $\alpha$. Note that by the symmetry $A(z)=A(z^{-1})$ of the Alexander polynomial,
$A(K)(e^h)$ contains only even powers of $h$, as suggested by the formula.

We need to understand how $\alpha$ acts on wheels. Due to the two-in-one-out
rule, a wheel is zero unless all the ``spokes'' are oriented inward, and the cycle oriented in
one direction. In other words, there are two ways to orient an unoriented wheel:
clockwise or counterclockwise. Due to 
the anti-symmetry of chord vertices, we get that for odd wheels $\alpha(w^u_{2h+1})=0$ and
for even wheels $\alpha(w^u_{2h})=2w^w_{2h}$. As a result,
$$\alpha Z^u(K)=\exp_{\calA^{sw}}\left(-\frac{1}{2} \log A(K)(e^h)|_{h^{2n}\to 2w_{2n}}\right)
=\exp_{\calA^{sw}}\left(-\log A(K)(e^h)|_{h^{2n}\to w_{2n}}\right)$$
which agrees with the Formula (\ref{eq:AtoZ}) of Theorem
\ref{thm:Alexander}. Note that since $K$ is 0-framed, the first part
(``$\sl$ coded in arrows'') of Formula~(\ref{eq:AtoZ}) is trivial.

\draftcut
\section{Odds and Ends} \label{sec:OddsAndEnds}

\draftcut
\subsection{Some Dimensions} \label{subsec:SomeDimensions}

The table below lists what we could find about $\calA^v$ and $\calA^w$ by
crude brute force computations in low degrees. We list degrees 0 through
7. The spaces we study are $\calA^-(\uparrow)$, $\calA^{s-}(\uparrow)$ 
(the $-$ in the subscript means ``$v$ and $w$''), and
$\calA^{r-}(\uparrow)$ which is $\calA^-(\uparrow)$ moded out by
``isolated'' arrows\footnote{That is, $\calA^{r-}(\uparrow)$
is $\calA^-(\uparrow)$ modulo ``framing independence'' (\glost{FI})
relations (see Section~\ref{subsec:Jacobi}, cf. ~\cite{Bar-Natan:OnVassiliev}, with the isolated arrow
taken with either orientation). It is the space related to finite type
invariants of unframed knots, on which the R1 move
is also imposed, in the same way as $\calA^-(\uparrow)$ is related
to framed knots.}, $\calP^-(\uparrow)$ which is the space of
primitives in $\calA^-(\uparrow)$, and $\glos{\calA^-(\bigcirc)}$,
$\glos{\calA^{s-}(\bigcirc)}$, and $\glos{\calA^{r-}(\bigcirc)}$,
which are the same as $\calA^-(\uparrow)$, $\calA^{s-}(\uparrow)$,
and $\calA^{r-}(\uparrow)$ except with closed knots (knots with
a circle skeleton) replacing long knots. Each of these spaces we
study in three variants: the ``v'' and the ``w'' variants, as well
as the \underline{u}sual knots ``u'' variant which is here just for
comparison. We also include a row ``$\dim\calG_m\calL ie^-(\uparrow)$''
for the dimensions of ``Lie-algebraic weight systems''. Those
are explained in the u and v cases in~\cite{Bar-Natan:OnVassiliev,
Haviv:DiagrammaticAnalogue, Leung:CombinatorialFormulas}, and in the w
case in Section~\ref{subsec:LieAlgebras}.

{
\def\uvw#1#2#3{{\hspace{-2.5mm}\text{\small
  $\begin{array}{c}#1\mid#2\\#3\end{array}$}\hspace{-3mm}
}}
\begin{center}\begin{tabular}{||c|c||c|c|c|c|c|c|c|c|c||}
\hline \hline
&& \multicolumn{3}{c|}{\footnotesize See Section~\ref{subsec:ToTwo}} &&&&&& \\
$m$ && 0 & 1 & 2 & 3 & 4 & 5 & 6 & 7 & \footnotesize Comments \\
\hline
$\dim\calG_m\calA^-(\uparrow)$ & \uvw{u}{v}{w} &
  \uvw{1}{1}{1} & \uvw{1}{2}{2} & \uvw{2}{7}{4} & \uvw{3}{27}{7} &
  \uvw{6}{139}{12} & \uvw{10}{813}{19} & \uvw{19}{?}{30} & \uvw{33}{?}{45}
  & \uvw{\ref{com:uknots}}{\ref{com:longv}}{\ref{com:wknots},
    \ref{com:longw}, \ref{com:nextfew}} \\
\hline
$\dim\calG_m\calL ie^-(\uparrow)$ & \uvw{u}{v}{w} &
  \uvw{1}{1}{1} & \uvw{1}{2}{2} & \uvw{2}{7}{4} & \uvw{3}{27}{7} &
  \uvw{6}{\,\geq\!128}{12} & \uvw{10}{?}{19} & \uvw{19}{?}{30} &
  \uvw{33}{?}{45}
  & \uvw{\ref{com:uknots}}{\ref{com:Lie}}{\ref{com:nextfew}} \\
\hline
$\dim\calG_m\calA^{s-}(\uparrow)$ & \uvw{u}{v}{w} &
  \uvw{-}{1}{1} & \uvw{-}{1}{1} & \uvw{-}{3}{2} & \uvw{-}{10}{3} &
  \uvw{-}{52}{5} & \uvw{-}{298}{7} & \uvw{-}{?}{11} & \uvw{-}{?}{15}
  & \uvw{\ref{com:su}}{\ref{com:longv}}{\ref{com:wknots}, \ref{com:nextfews}} \\
\hline
$\dim\calG_m\calA^{r-}(\uparrow)$ & \uvw{u}{v}{w} &
  \uvw{1}{1}{1} & \uvw{0}{0}{0} & \uvw{1}{2}{1} & \uvw{1}{7}{1} &
  \uvw{3}{42}{2} & \uvw{4}{246}{2} & \uvw{9}{?}{4} & \uvw{14}{?}{4}
  & \uvw{\ref{com:uknots}}{\ref{com:fiwarning}}{\ref{com:wknots},
    \ref{com:nextfewr}} \\
\hline
$\dim\calG_m\calP^-(\uparrow)$ & \uvw{u}{v}{w} &
  \uvw{0}{0}{0} & \uvw{1}{2}{2} & \uvw{1}{4}{1} & \uvw{1}{15}{1} &
  \uvw{2}{82}{1} & \uvw{3}{502}{1} & \uvw{5}{?}{1} & \uvw{8}{?}{1}
  & \uvw{\ref{com:uknots}}{\ref{com:Pv}}{\ref{com:wknots}} \\
\hline
$\dim\calG_m\calA^-(\bigcirc)$ & \uvw{u}{v}{w} &
  \uvw{1}{1}{1} & \uvw{1}{1}{1} & \uvw{2}{2}{1} & \uvw{3}{5}{1} &
  \uvw{6}{19}{1} & \uvw{10}{77}{1} & \uvw{19}{?}{1} & \uvw{33}{?}{1}
  & \uvw{\ref{com:uknots}}{\ref{com:closedv}}{\ref{com:wknots}} \\
\hline
$\dim\calG_m\calA^{s-}(\bigcirc)$ & \uvw{u}{v}{w} &
  \uvw{-}{1}{1} & \uvw{-}{1}{1} & \uvw{-}{1}{1} & \uvw{-}{2}{1} &
  \uvw{-}{6}{1} & \uvw{-}{23}{1} & \uvw{-}{?}{1} & \uvw{-}{?}{1}
  & \uvw{\ref{com:su}}{\ref{com:longv}}{\ref{com:wknots}} \\
\hline
$\dim\calG_m\calA^{r-}(\bigcirc)$ & \uvw{u}{v}{w} &
  \uvw{1}{1}{1} & \uvw{0}{0}{0} & \uvw{1}{0}{0} & \uvw{1}{1}{0} &
  \uvw{3}{4}{0} & \uvw{4}{17}{0} & \uvw{9}{?}{0} & \uvw{14}{?}{0}
  & \uvw{\ref{com:uknots}}{\ref{com:closedv}}{\ref{com:wknots}} \\
\hline \hline
\end{tabular}\end{center}
}

\begin{comments} \begin{enumerate}
\item \label{com:uknots} Much more is known computationally on the u-knots
  case.  See especially~\cite{Bar-Natan:OnVassiliev, Bar-Natan:Computations,
  Kneissler:Twelve, Amir-KhosraviSankaran:VasCalc}.
\item \label{com:longv} These dimensions were computed by Louis Leung and
  DBN using a program available at~\cite[``Dimensions'']{WKO}.
\item \label{com:wknots} As we have seen in Section~\ref{subsec:Jacobi},
  the spaces associated with w-knots are understood to all degrees.
\item \label{com:longw} To degree 4, these numbers were also verified
  by~\cite[``Dimensions'']{WKO}.
\item \label{com:nextfew} The next few numbers in these sequences are 67, 97,
  139, 195, 272.
\item \label{com:Lie} These dimensions were computed by Louis Leung and
  DBN using a program available at~\cite[``Arrow Diagrams and
  $\mathfrak{gl}(N)$'']{WKO}. Note the match with the row above.
\item \label{com:su} There is no ``s'' quotient in the ``u'' case.
\item \label{com:nextfews} The next few numbers in this sequence are 22, 30,
  42, 56, 77.
\item \label{com:fiwarning} These numbers were computed
  by~\cite[``Dimensions'']{WKO}. Contrary to the $\calA^u$
  case, $\calA^{rv}$ is {\em not} the quotient of $\calA^{v}$ by the
  ideal generated by degree 1 elements, and therefore the dimensions
  of the graded pieces of these two spaces cannot be deduced from each
  other using the Milnor-Moore theorem.
\item \label{com:nextfewr} The next few numbers in this sequence are
  7,8,12,14,21.
\item \label{com:Pv} These dimensions were deduced from the dimensions of
  $\calG_m\calA^v(\uparrow)$ using the Milnor-Moore theorem.
\item \label{com:closedv} Computed
  by~\cite[``Dimensions'']{WKO}. Contrary to the $\calA^u$
  case, $\calA^v(\bigcirc)$, $\calA^{sv}(\bigcirc)$, and
  $\calA^{rv}(\bigcirc)$ are {\em not} isomorphic to $\calA^v(\uparrow)$,
  $\calA^{sv}(\uparrow)$, and $\calA^{rv}(\uparrow)$ and separate
  computations are required.
\end{enumerate}
\end{comments}

\subsection{What Means ``Closed Form''?} \label{subsec:ClosedForm}

As stated earlier, one of our hopes for this sequence of papers is that it will
lead to closed-form formulae for tree-level associators. The
notion ``closed-form'' in itself requires an explanation. Is $e^x$ a closed form expression for
$\sum_{n=0}^\infty\frac{x^n}{n!}$, or is it just an artificial name given
for a transcendental function we cannot otherwise reduce? Likewise,
why not call some tree-level associator $\Phi^\text{tree}$ and now it is
``in closed form''?

For us, ``closed-form'' should mean ``useful for computations''. More
precisely, it means that the quantity in question is an element of some
space $\calAcf$ of ``useful closed-form thingies'' whose elements have
finite descriptions (hopefully, finite and short) and on which some
operations are defined by algorithms which terminate in finite time
(hopefully, finite and short). Furthermore, there should be a finite-time
algorithm to decide whether two descriptions of elements of $\calAcf$
describe the same element\footnote{In our context, if it is hard to
decide within the target space of an invariant whether two elements
are equal or not, the invariant is not too useful in deciding whether
two knotted objects are equal or not.}. It is even better if the said
decision algorithm takes the form ``bring each of the two elements in question
to a canonical form by means of some finite (and hopefully short)
procedure, and then compare the canonical forms verbatim''; if this is the
case, then many algorithms that involve managing a large number of elements
become simpler and faster.

Thus, for example, polynomials in a variable $x$ are always of closed form,
for they are simply described by finite sequences of integers (which in
themselves are finite sequences of digits), the standard operations on
polynomials ($+$, $\times$, and, say, $\frac{d}{dx}$) are algorithmically
computable, and it is easy to write the ``polynomial equality'' computer
program. Likewise for rational functions and even for rational functions
of $x$ and $e^x$.

On the other hand, general elements $\Phi$ of the space
$\calA^\text{tree}(\uparrow_3)$ of potential tree-level associators
are not closed-form, for they are determined by infinitely many
coefficients. Thus, iterative constructions of associators, such
as the one in~\cite{Bar-Natan:NAT} are computationally useful only
within bounded-degree quotients of $\calA^\text{tree}(\uparrow_3)$
and not as all-degree closed-form formulae. Likewise, ``explicit''
formulae for an associator $\Phi$ in terms of multiple $\zeta$-values
(e.g.~\cite{LeMurakami:HOMFLY}) are not useful for computations as it
is not clear how to apply tangle-theoretic operations to $\Phi$ (such as
$\Phi\mapsto\Phi^{1342}$ or $\Phi\mapsto(1\otimes\Delta\otimes 1)\Phi$)
while staying within some space of ``objects with finite description in
terms of multiple $\zeta$-values''. And even if a reasonable space of such
objects could be defined, it remains an open problem to decide whether
a given rational linear combination of multiple $\zeta$-values is equal
to $0$.

\draftcut
\subsection{Arrow Diagrams up to Degree 2} \label{subsec:ToTwo} Just as
an example, in this section we study the spaces $\calA^-(\uparrow)$,
$\calA^{s-}(\uparrow)$, $\calA^{r-}(\uparrow)$, $\calP^-(\uparrow)$,
$\calA^-(\bigcirc)$, $\calA^{s-}(\bigcirc)$, and $\calA^{r-}(\bigcirc)$ in
degrees $m\leq 2$ in detail, both in the ``v'' case and in the ``w'' case
(the ``u'' case has been known since long \cite{Bar-Natan:OnVassiliev,
Kneissler:Twelve, Bar-Natan:Computations}).

\subsubsection{Arrow Diagrams in Degree 0} There is only
one degree 0 arrow diagram, the empty diagram $D_0$ (see
Figure~\ref{fig:Deg0-2Diagrams}).  There are no relations, and thus,
$\{D_0\}$ is the basis of all $\calG_0\calA^-(\uparrow)$ spaces
and its closure, the empty circle, is the basis of all
$\calG_0\calA^-(\bigcirc)$ spaces. $D_0$ is the unit $1$, yet $\Delta
D_0=D_0\otimes D_0=1\otimes 1\neq D_0\otimes 1+1\otimes D_0$, so $D_0$
is not primitive and $\dim\calG_0\calP^-(\uparrow)=0$.

\subsubsection{Arrow Diagrams in Degree 1} \label{subsubsec:DegreeOne}
\begin{figure}
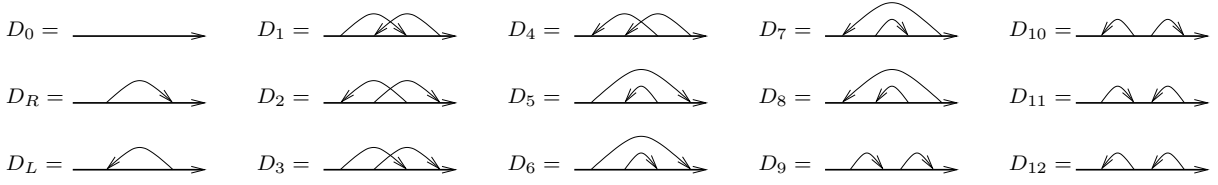

\[ \input figs/Deg0-2Diagrams.pstex_t \]
\caption{The 15 arrow diagrams of degree at most 2.}
\label{fig:Deg0-2Diagrams}
\end{figure}

There is only two degree 1 arrow diagrams, the ``right
arrow'' diagram $D_R$ and the ``left arrow'' diagram $D_L$ (see
Figure~\ref{fig:Deg0-2Diagrams}).  There are no $6T$ relations, and
thus, $\{D_R, D_L\}$ is the basis of $\calG_1\calA^-(\uparrow)$. Modulo
RI, $D_L=D_R$ and hence, $D_A:=D_L=D_R$ is the single basis element of
$\calG_1\calA^{s-}(\uparrow)$. Both $D_R$ and $D_L$ vanish modulo FI, so
$\dim\calG_1\calA^{r-}(\uparrow)=\dim\calG_1\calA^{r-}(\bigcirc)=0$. Both
$D_R$ and $D_L$ are primitive, so $\dim\calG_1\calP^-(\uparrow)=2$.
Finally, the closures ${\bar D}_R$ and ${\bar D}_L$ of $D_R$ and $D_L$ are
equal, so $$\calG_1\calA^{s-}(\bigcirc)=\calG_1\calA^-(\bigcirc)=\langle
{\bar D}_R\rangle=\langle {\bar D}_L\rangle=\langle {\bar D}_A\rangle.$$

\subsubsection{Arrow Diagrams in Degree 2} There are 12 degree
2 arrow diagrams, which we denote $D_1,\dots,D_{12}$ (see
Figure~\ref{fig:Deg0-2Diagrams}). There are six $6T$ relations,
corresponding to the 6 ways of ordering the 3 vertical strands that
appear in a $6T$ relation (see Figure~\ref{fig:6T}) along a long line. The
ordering $(ijk)$ becomes the relation $D_3+D_9+D_3=D_6+D_3+D_6$. Likewise,
$(ikj)\mapsto D_6+D_1+D_{11}=D_3+D_5+D_1$, $(jik)\mapsto
D_{10}+D_2+D_6=D_2+D_5+D_3$, $(jki)\mapsto D_4+D_7+D_1=D_8+D_1+D_{11}$,
$(kij)\mapsto D_2+D_7+D_4=D_{10}+D_2+D_8$, and $(kji)\mapsto
D_8+D_4+D_8=D_4+D_{12}+D_4$. After some linear algebra, we find
that $\{D_1, D_2, D_6, D_8, D_9, D_{11}, D_{12}\}$ form a basis
of $\calG_2\calA^v(\uparrow)$, and that the remaining diagrams
reduce to the basis as follows: $D_3=2D_6-D_9$, $D_4=2D_8-D_{12}$,
$D_5=D_9+D_{11}-D_6$, $D_7=D_{11}+D_{12}-D_8$, and $D_{10}=D_{11}$.  In
$\calG_2\calA^{sv}(\uparrow)$ we further have that $D_5=D_6$, $D_7=D_8$,
and $D_9=D_{10}=D_{11}=D_{12}$, and so $\calG_2\calA^{sv}(\uparrow)$
is 3-dimensional with basis $D_1$, $D_2$, and $D_3=\ldots=D_{12}$.
In $\calG_2\calA^{rv}(\uparrow)$ we further have that $D_{5-12}=0$.
Thus $\{D_1, D_2\}$ is a basis of $\calG_2\calA^{rv}(\uparrow)$.

There are 3 OC relations to write for $\calG_2\calA^w(\uparrow)$:
$D_2=D_{10}$, $D_3=D_6$, and $D_4=D_8$. Along with the $6T$ relations,
we find that $$\{D_1, D_3=D_6=D_9, D_2=D_5=D_7=D_{10}=D_{11},
D_4=D_8=D_{12}\}$$ is a basis of $\calG_2\calA^w(\uparrow)$. Similarly
$\{D_1,D_2=\ldots=D_{12}\}$ is a basis of the two-dimensional
$\calG_2\calA^{sw}(\uparrow)$.  When we mod out by FI, only one
diagram remains non-zero in $\calG_2\calA^{rw}(\uparrow)$ and it is $D_1$.

We leave the determination of the primitives and the spaces with a circle
skeleton as an exercise to the reader.

\draftcut \section{Glossary of Notation} \label{sec:glossary}
Greek letters, then Latin, then symbols:

\noindent
{\small \begin{multicols}{2}
\begin{list}{}{
  \renewcommand{\makelabel}[1]{#1\hfil}
}

% \alpha
\item[{$\alpha$}] maps $\calA^u\to \calA^v$ or $\calA^u\to \calA^w$
  \hfill\ \ref{subsubsec:RelWithu}
% \beta
% \gamma
% \delta
\item[{$\Delta$}] cloning, co-product\hfill\ \ref{par:Delta}
\item[{$\delta$}] Satoh's tube map\hfill\ \ref{subsubsec:TopTube}
\item[{$\delta_A$}] a formal $D_A$\hfill\ \ref{subsec:AlexanderProof}
% \epsilon
% \zeta
% \eta
% \theta
\item[{$\theta$}] inversion, antipode\hfill\ \ref{par:theta}
% \iota
\item[{$\iota$}] an inclusion $\wB_n\to\wB_{n+1}$\hfill\
  \ref{subsubsec:McCool}
\item[{$\iota$}] interpretation map\hfill\ \ref{subsec:AlexanderProof},
  \ref{subsubsec:IAM}
% \kappa
% \lambda
\item[{$\lambda$}] a formal $EZ$\hfill\ \ref{subsec:AlexanderProof}
% \mu
% \nu
% \xi
\item[{$\xi_i$}] the generators of $F_n$\hfill\ \ref{subsubsec:McCool}
% \omicron
% \pi
% \rho
% \sigma
\item[{$\Sigma$}] a virtual surface\hfill\ \ref{subsubsec:TopTube}
\item[{$\sigma_i$}] a crossing between adjacent strands\hfill\ 
  \ref{subsubsec:Planar}
\item[{$\sigma_{ij}$}] strand $i$ crosses over strand $j$\hfill\ 
  \ref{subsubsec:Abstract}
\item[{$\varsigma$}] the skeleton morphism\hfill\ \ref{subsubsec:Planar}
% \tau
% \upsilon
% \phi
\item[{$(\varphi^i)$}] a basis of $\frakg^\ast$\hfill\ \ref{subsec:LieAlgebras}
% \chi
% \psi
% \omega
\item[{$\omega_1$}] a formal 1-wheel\hfill\ \ref{subsec:AlexanderProof}

\item
%a
\item[{$\calA(G)$}] associated graded of $G$\hfill\
  \ref{subsubsec:FTAlgebraic}
\item[{$\calA^-_n$}] $\calD^v_n$ mod relations\hfill\
  \ref{subsubsec:FTPictorial}
\item[{$\calA^{-t}$}] $\calA^-$ allowing trivalent vertices\hfill\
  \ref{subsec:Jacobi}
\item[{$\calA^-(\uparrow)$}] $\calD^v(\uparrow)$ mod relations\hfill\
  \ref{subsec:FTforvwKnots}
\item[{$\calA^-(\bigcirc)$}] $\calA^-(\uparrow)$ for round skeletons\hfill\
  \ref{subsec:SomeDimensions}
\item[{$\calA^u$}] usual chord diagrams\hfill\ \ref{subsec:RelWithKont}
\item[{$A(K)$}] the Alexander polynomial\hfill\ \ref{subsec:Alexander}
\item[{$\aAS$}] arrow-AS relations\hfill\ \ref{subsec:Jacobi}
\item[{$a$}] maps $u\to v$ or $u\to w$\hfill\ \ref{subsubsec:RelWithu}
\item[{$a_{ij}$}] an arrow from $i$ to $j$\hfill\
  \ref{subsubsec:FTPictorial}
%b
\item[{$\calB^w$}] unitrivalent arrow diagrams\hfill\ \ref{subsec:Jacobi}
\item[{$B$}] the matrix $T(\exp(-xS)-I)$\hfill\ \ref{subsec:AlexanderProof}
\item[{$b_{ij}^k$}] structure constants of $\frakg^\ast$\hfill\ 
  \ref{subsec:LieAlgebras}
%c
\item[{CC}] the Commutators Commute relation\hfill\ \ref{subsec:Jacobi}
%d
\item[{$\calD^v_n$}] arrow diagrams for braids\hfill\
  \ref{subsubsec:FTPictorial}
\item[{$\calD^{-t}$}] $\calD^-$ allowing trivalent vertices\hfill\
  \ref{subsec:Jacobi}
\item[{$\calD^v(\uparrow)$}] arrow diagrams long knots\hfill\ 
  \ref{subsec:FTforvwKnots}
\item[{$D_A$}] either $D_L$ or $D_R$\hfill\ \ref{subsec:Jacobi}
\item[{$D_L$}] left-going isolated arrow\hfill\ \ref{subsec:Jacobi}
\item[{$D_R$}] right-going isolated arrow\hfill\ \ref{subsec:Jacobi}
\item[{$d_k$}] strand deletion\hfill\ \ref{par:deletions}
\item[{$d_i$}] the direction of a crossing\hfill\ \ref{subsec:Alexander}
%e
\item[{$E$}] the Euler operator\hfill\ \ref{subsec:AlexanderProof}
\item[{$\tilE$}] the normalized Euler operator\hfill\ 
  \ref{subsec:AlexanderProof}
%f
\item[{FI}] Framing Independence\hfill\ \ref{subsec:SomeDimensions}
\item[{$F_n$}] the free group\hfill\ \ref{subsubsec:McCool}
\item[{$\FA_n$}] the free associative algebra\hfill\ \ref{par:action}
%g
\item[{$\frakg$}] a finite-dimensional Lie algebra\hfill\
  \ref{subsec:LieAlgebras}
\item[{$\calG_m$}] degree $m$ piece\hfill\ \ref{subsubsec:FTPictorial}
%h
%i
\item[{$\calI$}] augmentation ideal\hfill\ \ref{subsubsec:FTAlgebraic},
\item[{$I\frakg$}] $\frakg^\ast\rtimes\frakg$\hfill\ \ref{subsec:LieAlgebras}
\item[{$\IAM$}] Infinitesimal Alexander\newline Module\hfill\ 
  \ref{subsec:AlexanderProof}, \ref{subsubsec:IAM}
\item[{$\IAM^0$}] $\IAM$, before relations\hfill\ \ref{subsubsec:IAM}
\item[{$\aIHX$}] arrow-IHX relations\hfill\ \ref{subsec:Jacobi}
\item[{$i_u$}] an inclusion $F_n\to\wB_{n+1}$\hfill\ \ref{subsubsec:McCool}
%j
%k
\item[{$\calK^u$}] usual knots\hfill\ \ref{subsec:RelWithKont}
%l
%m
\item[{M}] the ``mixed'' move\hfill\ \ref{subsec:VirtualKnots}
%n
%o
\item[{OC}] the Overcrossings Commute relation\hfill\ \ref{subsec:wBraids}
%p
\item[{$\calP^-(\uparrow)$}] primitives of $\calA^-(\uparrow)$\hfill\ 
  \ref{subsec:FTforvwKnots}
\item[{$\PvB_n$}] the group of pure v-braids\hfill\ \ref{subsubsec:Planar}
\item[{$\PwB_n$}] the group of pure w-braids\hfill\ \ref{subsec:wBraids}
%q
%r
\item[{$\calR$}] the relations in $\IAM$\hfill\ \ref{subsubsec:IAM}
\item[{$R$}] $Z(\overcrossing)$\hfill\ \ref{subsec:wBraidExpansion}
\item[{$R$}] the ring $\bbZ[X,X^{-1}]$\hfill\ \ref{subsubsec:IAM}
\item[{$R_1$}] the augmentation ideal of $R$\hfill\ \ref{subsubsec:IAM}
\item[{RI}] Rotation number Independence\hfill\ \ref{subsec:FTforvwKnots}
\item[{R123}] Reidemeister moves\hfill\ \ref{subsec:VirtualKnots}
%\item[{R4}] a Reidemeister move for foams/graphs\hfill\ \ref{subsubsec:wrels}
\item[{\Rs}] the ``spun'' R1 move\hfill\ \ref{subsec:VirtualKnots}
%s
\item[{$S(K)$}] a matrix of signs\hfill\ \ref{subsec:Alexander}
\item[{$S_n$}] the symmetric group\hfill\ \ref{subsubsec:Planar}
\item[{$\aSTU$}] arrow-STU relations\hfill\ \ref{subsec:Jacobi}
\item[{$s_i$}] a virtual crossing between adjacent strands\hfill\ 
  \ref{subsubsec:Planar}
\item[{$s_i$}] the sign of a crossing\hfill\ \ref{subsec:Alexander}
\item[{$\sl$}] self-linking\hfill\ \ref{subsec:VirtualKnots}
%t
\item[{$\calT^w_\frakg$}] a map ${\calA}^w\to\calU(I\frakg)$\hfill\ 
  \ref{subsec:LieAlgebras}
\item[{TC}] Tails Commute\hfill\ \ref{subsubsec:FTPictorial}
\item[{$T(K)$}] the ``trapping'' matrix\hfill\ \ref{subsec:Alexander}
%u
\item[{$\calU$}] universal enveloping algebra\hfill\ \ref{subsec:LieAlgebras}
\item[{UC}] Undercrossings Commute\hfill\ \ref{subsec:wBraids}
\item[{$u_k$}] strand unzips\hfill\ \ref{par:unzip}
\item[{$\uB_n$}] the (usual) braid group\hfill\ \ref{subsubsec:Planar}
%v
\item[{$V$}] a finite-type invariant\hfill\ \ref{subsubsec:FTPictorial}
\item[{VR123}] virtual Reidemeister moves\hfill\ \ref{subsec:VirtualKnots}
\item[{$\vB_n$}] the virtual braid group\hfill\ \ref{subsubsec:Planar}
%w
\item[{$W_m$}] weight system\hfill\ \ref{subsubsec:FTPictorial}
\item[{$w$}] the map $x^k\mapsto w_k$\hfill\ \ref{subsec:Alexander}
\item[{$w_i$}] flip ring $\#i$\hfill\ \ref{subsubsec:FlyingRings}
\item[{$w_k$}] the $k$-wheel\hfill\ \ref{subsec:Jacobi}
\item[{$\wB_n$}] the group of w-braids\hfill\ \ref{subsec:wBraids}
%x
\item[{$X$}] an indeterminate\hfill\ \ref{subsec:Alexander}
\item[{$X_n,\, \tilde{X}_n$}] moduli of horizontal rings\hfill\ 
  \ref{subsubsec:FlyingRings}
\item[{$x_i$}] the generators of $FA_n$\hfill\ \ref{par:action}
\item[{$(x_j)$}] a basis of $\frakg$\hfill\ \ref{subsec:LieAlgebras}
%y
\item[{$Y_n,\, \tilde{Y}_n$}] moduli of rings\hfill\ \ref{subsubsec:NonHorRings}
%z
\item[{$Z$}] expansions \hfill\ throughout
\item[{$Z^u$}] the Kontsevich integral\hfill\ \ref{subsec:RelWithKont}

\item
\item[{$\aft$}] $\aft$ relations\hfill\ \ref{subsubsec:FTPictorial}
\item[{$6T$}] $6T$ relations\hfill\ \ref{subsubsec:FTPictorial}
\item[{$\semivirtualover,\,\semivirtualunder$}] semi-virtual
  crossings\hfill\ \ref{subsubsec:FTPictorial}
\item[{$\sslash$}] right action\hfill\ \ref{subsubsec:McCool}
\item[{$\uparrow$}] a ``long'' strand\hfill\ throughout

\end{list}
\end{multicols}}

\draftcut

\if\draft y
  \clearpage
  Everything below is to be blanked out before the completion of this paper.
  \section*{Recycling}

\begin{exercise} Do the same for the obviously-defined ``w-links'',
excluding the material about the Alexander polynomial. Note that the wheels
that are obtained in the case of w-links have legs coloured by the
components of the w-link in question. Hence if there is more than one
component, the number of such wheels grows exponentially in the degree and
thus $Z$ contains more information than can be coded in a polynomial of
even a multi-variable polynomial.
\end{exercise}

\vskip -5mm
\parpic[r]{\raisebox{-14mm}{$\input figs/CC.pstex_t$}}
\begin{conjecture} In the case of ordinary links seen as w-links,
if we mod out the target space of $Z$ by the ``Commutators Commute''
relation shown on the right, what remains of the wheels part of $Z$
is precisely the multi-variable Alexander polynomial.
\end{conjecture}

Note that $D \in \tder_n$ is never an arrow on a single strand (these are
elements of $\fraka_n$),
and hence $\operatorname{div}D$ is never a 1-wheel, more precisely it never
has a degree 1 component. Thus, even
though the target space of div is $\attr_n/(\text{deg }1)$, we can just as
well think of it as a
map to $\attr_n$ itself.

\subsection{The Injectivity of $i_u\colon F_n\to\wB_{n+1}$}
\label{subsec:FreeInW}

Just for completeness, we sketch here an algebraic proof of the
injectivity of the map $i_u\colon F_n\to\wB_{n+1}$ discussed in
Section~\ref{subsubsec:McCool}. There's some circularity in our argument
--- we need this injectivity in order to motivate the definition of the
map $\Psi\colon \wB_n\to\Aut(F_n)$, and in the proof below we use $\Psi$
to prove the injectivity of $i_u$. But $\Psi$ exists regardless of how
its definition is motivated, and it can be shown to be well defined by
explicitly verifying that it respects the relations defining $\wB_n$. So
our proof is logically valid.

\begin{claim} The map $i_u\colon F_n\to\wB_{n+1}$ is injective.
\end{claim}

\begin{proof} (sketch). Let $H$ be the subgroup of $\wB_{n+1}$ MORE
\end{proof}

\subsection{Finite Type Invariants of v-Braids and w-Braids, in some
Detail} \label{subsec:FTDetails}

As mentioned in Section~\ref{subsec:wBraids}, w-braids are v-braids modulo
an additional relation. So we start with a discussion of finite type
invariants of v-braids. For simplicity we take our base ring to be $\bbQ$;
everywhere we could replace it by an arbitrary field of characteristic
$0$\footnote{Or using the variation of constants method, we can simply
declare that $\bbQ$ is an arbitrary field of characteristic $0$.},
and many definitions make sense also over $\bbZ$ or even with $\bbQ$
replaced by an arbitrary Abelian group.

\subsubsection{Basic Definitions.}
Let $\bbQ\vB_n$ denote group ring of $\vB_n$, the algebra of formal linear
combinations of elements of $\vB_n$, and let $\bbQ S_n$ be the group
ring of $S_n$. The skeleton homomorphism of Remark~\ref{rem:Skeleton}
extends to a homomorphism $\varsigma\colon \bbQ\vB_n\to\bbQ S_n$. Let $\calI$
(or $\calI_n$ when we need to be more precise) denote the kernel of
the skeleton homomorphism; it is the ideal in $\bbQ\vB_n$ generated
by formal differences of v-braids having the same skeleton. One may
easily check that $\calI$ is generated by differences of the form
$\overcrossing-\virtualcrossing$ and $\virtualcrossing-\undercrossing$.
Following~\cite{GoussarovPolyakViro:VirtualKnots} we call such differences
``semi-virtual crossings'' and denote them by $\semivirtualover$
and $\semivirtualunder$, respectively\footnote{The signs in
$\semivirtualover\leftrightarrow\overcrossing-\virtualcrossing$ and
$\semivirtualunder\leftrightarrow\virtualcrossing-\undercrossing$ are
``crossings come with their sign and their virtual counterparts come with
the opposite sign''.}. In a similar manner, for any natural number $m$
the $m$th power $\calI^m$ of $\calI$ is generated by ``$m$-fold iterated
differences'' of v-braids, or equally well, by ``$m$-singular v-braids'',
which are v-braids that also have exactly $m$ semi-virtual crossings
(subject to relations which we don't need to specify).

Let $V\colon \vB_n\to A$ be an invariant of v-braids with values in some vector
space $A$. We say that $V$ is ``of type $m$'' (for some $m\in\bbZ_{\geq
0}$) if its linear extension to $\bbQ\vB_n$ vanishes on $\calI^{m+1}$
(alternatively, on all $(m+1)$-singular v-braids, in clear analogy with the
standard definition of finite type invariants). If $V$ is of type $m$
for some unspecified $m$, we say that $V$ is ``of finite type''. Given
a type $m$ invariant $V$, we can restrict it to $\calI^m$ and as it
vanishes on $\calI^{m+1}$, this restriction can be regarded as an element
of $\left(\calI^m/\calI^{m+1}\right)^\star$. If two type $m$ invariants
define the same element of $\left(\calI^m/\calI^{m+1}\right)^\star$ then
their difference vanishes on $\calI^m$, and so it is an invariant of type
$m-1$. Thus it is clear that an understanding of $\calI^m/\calI^{m+1}$
will be instrumental to an inductive understanding of finite type
invariants. Hence the following definition.

\begin{definition} The projectivization\footnote{Why ``projectivization''?
See Section~\ref{subsec:Projectivization}.} $\proj\vB_n$ is the direct sum
\[ \proj\vB_n:=\bigoplus_{m\geq 0} \calI^m/\calI^{m+1}. \]
Note that throughout this paper, whenever we write an infinite direct sum,
we automatically complete it. Therefore an element in $\proj\vB_n$ is an
infinite sequence $D=(D_0, D_1,\dots)$, where $D_m\in\calI^m/\calI^{m+1}$.
The projectivization $\proj\vB_n$ is a graded space, with the degree $m$
piece being $\calI^m/\calI^{m+1}$.
\end{definition}

We proceed with the study of $\proj\vB_n$ (and thus of finite type
invariants of v-braids) in three steps. In
Section~\ref{subsubsec:ArrowDiagrams} we introduce a space $\calD^v_n$ and
a surjection $\rho_0\colon \calD^v_n\to\proj\vB_n$. In Section~\ref{subsubsec:6T}
we find some relations in $\ker\rho_0$, most notably the $6T$ relation, and
introduce the quotient $\calA^v_n:=\calD^v_n/6T$. And then in
Section~\ref{subsubsec:UFTI} we introduce the notion of a ``universal
finite type invariant'' and explain how the existence of such a gadget
proves that $\proj\vB_n$ is isomorphic to $\calA^v_n$ (in a more
traditional language this is the statement that every weight system
integrates to an invariant).

Unfortunately, we do not know if there is a universal finite type invariant
of v-braids. Thus in Section~\ref{subsec:wbraids} we return to the subject
of w-braids and prove the weaker statement that there exists a universal
finite type invariant of w-braids.

\subsubsection{Arrow Diagrams.} \label{subsubsec:ArrowDiagrams}

We are looking for a space that will surject on $\calI^m/\calI^{m+1}$. In
other words, we are looking for a set of generators for $\calI^m$, and
we are willing to call two such generators the same if their difference
is in $\calI^{m+1}$. But that's easy. Left and right multiples of the
formal differences $\semivirtualover=\overcrossing-\virtualcrossing$
and $\semivirtualunder=\virtualcrossing-\undercrossing$
generate $I$, so products of the schematic form
\begin{equation} \label{eq:GeneratingProduct}
  B_0 (\semivirtualover|\semivirtualunder) B_1
  (\semivirtualover|\semivirtualunder) B_2 \cdots
  B_{m-1} (\semivirtualover|\semivirtualunder) B_m
\end{equation}
\parpic[r]{\input figs/SemiVirtRels.pstex_t }
\noindent generate $\calI^m$ (here $(\semivirtualover|\semivirtualunder)$
means ``either a $\semivirtualover$ or a $\semivirtualunder$'', and there
are exactly $m$ of those in any product). Furthermore, inside such a
product any $B_k$ can be replaced by any other v-braid $B'_k$ having the
same skeleton (e.g., with $\varsigma(B_k)$), for then $B_k-B'_k\in\calI$
and the whole product changes by something in $\calI^{m+1}$.  Also,
the relations in~\eqref{eq:R3} and in~\eqref{eq:MixedRelations} imply
the relations shown on the right for $\semivirtualover$, and similar
relations for $\semivirtualunder$. With this freedom, a product as
in~\eqref{eq:GeneratingProduct} is determined by the strand-placements
of the $\semivirtualover$'s and the $\semivirtualunder$'s. That is,
for each semi-virtual crossing in such a product, we only need to
know which strand number is the ``over'' strand, which strand number
is the ``under'' strand, and a sign that determines whether it is the
positive semi-virtual $\semivirtualover$ or the negative semi-virtual
$\semivirtualunder$. This motivates the following definition.

\begin{definition} A ``horizontal $m$-arrow diagrams'' (analogues to the
``chord diagrams'' of, say, \cite{Bar-Natan:OnVassiliev}) is an ordered
pair $(D,\beta)$ in which $D$ is a word of length $m$ in the alphabet
$\{a^+_{ij},a^-_{ij}\colon i,j\in\{1,\ldots,n\},\,i\neq j\}$ and $\beta$ is a
permutation in $S_n$. Let $\calD_m^{vh}$ be the space of formal linear
combinations of horizontal $m$-arrow diagrams. We usually use a pictorial
notation for horizontal arrow diagram, as demonstrated in
Figure~\ref{fig:Dvh}.
\end{definition}

\begin{figure}
\parpic[r]{\hspace{-5mm}\raisebox{-29mm}{$\input figs/Dvh.pstex_t$}}
\caption{
  The horizontal $3$-arrow diagram
  $(D,\beta)=$ $(a^+_{12}a^-_{41}a^+_{23},\,3421)$ and its image via 
  $\rho_0$. The first arrow, $a^+_{12}$ starts at strand $1$, ends 
  at strand $2$ and carries a $+$ sign, so it is mapped to a positive
  semi-virtual crossing of strand $1$ over strand $2$. Likewise the second
  arrow $a^-_{41}$ maps to a negative semi-virtual crossing of strand 
  $4$ over strand $1$, and $a^+_{23}$ to a positive semi-virtual
  crossing of strand $2$ over strand $3$. We also show one possible
  choice for a representative of the image of $\rho_0(D,\beta)$ in
  $\calI^m/\calI^{m+1}$: it is a v-braid with semi-virtual crossings as
  specified by $D$ and whose overall skeleton is $3421$.
} \label{fig:Dvh}
\end{figure}

There is a surjection $\rho_0\colon \calD_m^{vh}\to\calI^m/\calI^{m+1}$. The
definition of $\rho_0$ is suggested by the first paragraph of this section
and an example is shown in Figure~\ref{fig:Dvh}; we will skip the formal
definition here. We also skip the formal proof of the surjectivity of
$\rho_0$.

Finally, consider the product $\semivirtualover\cdot\semivirtualunder$ and
use the second Reidemeister move for both virtual and non-virtual
crossings:
\[
  \semivirtualover\semivirtualunder
  = (\overcrossing-\virtualcrossing)(\virtualcrossing-\undercrossing)
  = \overcrossing\virtualcrossing+\virtualcrossing\undercrossing
    - \overcrossing\undercrossing - \virtualcrossing\virtualcrossing
  = (\overcrossing\virtualcrossing-1) + (\virtualcrossing\undercrossing)
  = \semivirtualover\virtualcrossing - \virtualcrossing\semivirtualunder.
\]
If a total of $m-1$ further semi-virtual crossings are multiplied into this
equality on the left and on the right, along with arbitrary further
crossings and virtual crossings, the left hand side of the equality becomes
a member of $\calI^{m+1}$, and therefore, as a member of
$\calI^m/\calI^{m+1}$, it is $0$. Thus with ``$\ldots$'' standing for
extras added on the left and on the right, we have that in
$\calI^m/\calI^{m+1}$,
\[ 0 = 
  \ldots(
    \semivirtualover\virtualcrossing-\virtualcrossing\semivirtualunder
  )\ldots
  = \rho_0(\ldots??\ldots)
\]

MORE.

\subsubsection{The $6T$ Relations.} \label{subsubsec:6T}

MORE.

\subsubsection{The Notion of a Universal Finite Type Invariant.}
\label{subsubsec:UFTI}

MORE.

\subsection{Finite type invariants of w-braids} \label{subsec:wbraids}

MORE.
 
  \section*{To Do}

\par\noindent{\bf Sorted.}
\begin{itemize}
\item Finish the paper.
\item Freeze Mathematica notebooks.
\end{itemize}

\fi

\end{document}
\endinput